\documentclass[hidelinks,12pt]{amsart}

\usepackage{amsmath}
\usepackage[psamsfonts]{amssymb}
\usepackage{amsfonts}
\usepackage{graphicx}
\usepackage{placeins}

\usepackage{verbatim}
\IfFileExists{dsfont.sty}{\usepackage{dsfont}}{}
\usepackage{tikz-cd}
\usepackage{mathrsfs}
\usepackage{color}
\usepackage[english]{babel}
\usepackage{indentfirst}
\usepackage{tikz}
\usetikzlibrary{matrix,arrows,decorations.pathmorphing}
\IfFileExists{algorithm2e.sty}{\usepackage{algorithm2e}}{}
\usepackage[all]{xy}
\usepackage[T1]{fontenc}
\usepackage{hyperref}
\usepackage{mathtools}
\usepackage{multicol}
\usepackage{enumerate}

\newtheorem{theorem}{Theorem}[section]
\newtheorem{prop}[theorem]{Proposition}
\newtheorem{lemma}[theorem]{Lemma}
\newtheorem{cor}[theorem]{Corollary}

\theoremstyle{definition}
\newtheorem{example}[theorem]{Example}
\newtheorem{definition}[theorem]{Definition}
\theoremstyle{plain}

\newtheorem{remark}[theorem]{Remark}

\newcommand{\norm}[1]{\left\lVert#1\right\rVert}

\numberwithin{equation}{section}

\DeclareMathOperator{\R}{\mathbb{R}}

\def\ep{\epsilon}

\def\R{\mathbb{R}}
\def\Z{\mathbb{Z}}
\def\N{\mathbb{N}}

\def\C{\mathbb{C}}

\def\T{\mathbb T}

\def\cal R{\mathcal R}

\newtheorem{mainthm}{Theorem}

\allowdisplaybreaks[1]
\usepackage[letterpaper,textwidth=6.68in,textheight=8.5in,centering]{geometry}

\newcommand{\Palpha}{\mathrm P_\alpha}
\DeclareMathOperator{\Cont}{Cont}
\DeclareMathOperator{\supp}{supp}
\hypersetup{pdftitle={Contact rigidity and Poisson bracket invariants},pdfauthor={Jun Zhang}}

\newif\ifdraftnotes
\draftnotestrue

\author{Jun Zhang}
\email{jzhang4518@ustc.edu.cn}
\address{The Institute of Geometry and Physics, University of Science and Technology of China, 96 Jinzhai Road, Hefei Anhui, 230026, China}
 
\title{Contact rigidity and Poisson bracket invariants}

\DeclareRobustCommand{\rev}[1]{{#1}}
\DeclareRobustCommand{\keepmath}[1]{{#1}}
\DeclareRobustCommand{\mrev}[1]{{#1}}
\DeclareRobustCommand{\rrev}[1]{{#1}}
\begin{document}

\begin{abstract}
\rev{We study \keepmath{$C^0$}-rigidity and Poisson bracket invariants on closed contact manifolds. We prove that Poisson commutation is closed under uniform convergence, establish weak rigidity for pairs without common zeros, and obtain lower semicontinuity under quantitative control of Reeb derivatives. A three-dimensional counterexample shows that unrestricted weak rigidity fails. For partitions of unity, the contact bracket itself controls the Reeb derivatives. Using monotone homogeneous quasi-morphisms and the Cant--Kilgore--Zhang spectral invariant, we derive positive lower bounds for invariants of open covers and configurations of compact subsets. Applications include a quantitative contact big fiber theorem, bounds on Hamiltonian chord durations and Reeb travel times, and deformations of contact one-forms.}
\end{abstract}

\maketitle

\vspace*{-18pt}
\tableofcontents

\section{Introduction}\label{sec:introduction}

The Poisson bracket relates the algebra of smooth functions to the geometry of Hamiltonian flows. In symplectic geometry, the bracket measures the rate of change of one Hamiltonian along the flow of the other. In particular, when the bracket vanishes, the two flows commute. Since the bracket involves first derivatives, $C^0$-convergence of its inputs does not necessarily imply $C^0$-convergence of the bracket. However, symplectic topology imposes strong restrictions on how much a small uniform perturbation can decrease its norm \cite{CV,EP-C0,Buhovsky-C0}. These restrictions extend beyond individual pairs of functions: when the inputs satisfy conditions prescribed by subsets or an open cover, certain minimization problems can also reveal geometric rigidity \cite{EPZ,Polterovich-noise,BEP}. \rev{Chapters 2 and 7--9 of Polterovich--Rosen's book \cite{Polterovich-Rosen} develop this point of view in detail and provide much of the machinery adapted in this paper.}

\rev{Motivated by these symplectic results, we develop the same point of view in contact geometry.} \rev{The Poisson bracket is also defined in this setting (see (\ref{dfn-CBP})), where various forms of rigidity have been established using generating function methods \cite{Givental,Sandon,GKPS} and Floer-theoretic methods \cite{EKP,Albers-Merry,DUZ,CKZ}.} \rev{Building on these developments, we establish several forms of rigidity for the contact Poisson bracket analogous to those in symplectic geometry. Some proofs use the recently developed contact invariants in \cite{DUZ,CKZ}, and some results are closely related to the contact rigidity phenomena in \cite{SUV,DUZ}.} 

\subsection{\rev{Symplectic background and contact notation}}\label{subsec:background}

Let $(M,\omega)$ be a closed symplectic manifold. We use $\iota_{X_F}\omega=-dF$ and write $\{F,G\}_\omega=\rrev{X_F G}$. A basic form of symplectic $C^0$-rigidity is the lower-semicontinuity statement
\begin{equation} \label{eq:symplectic-lsc}
\norm{\{F,G\}_\omega}_\infty \le \liminf_{j\to\infty}\norm{\{F_j,G_j\}_\omega}_\infty
\end{equation}
\rev{for sequences of smooth Hamiltonian functions \keepmath{$F_j, G_j$} converging uniformly to the respective components of a smooth pair \keepmath{$(F,G)$};} see \cite{Zapolsky-surfaces,EP-C0} and the exposition in \cite[Section 2.1]{Polterovich-Rosen}. In particular, a non-commuting pair cannot be uniformly approximated by pairs whose brackets tend uniformly to zero. 

\rev{These rigidity questions also concern the geometry of subsets.} For a finite open cover $\mathcal U=\{U_1,\ldots,U_N\}$ of $M$, \rev{define the Poisson bracket invariant of the cover \keepmath{$\mathcal U$} by} 
\begin{equation} \label{eq:symplectic-cover}
pb_\omega(\mathcal U):=\inf_{\vec F}\max_{a,b\in[-1,1]^N}\rrev{\norm{\left\{\sum_i a_i F_i,\sum_j b_j F_j\right\}_\omega}_\infty},
\end{equation}
where the infimum is taken over smooth partitions of unity $\vec F=(F_1,\ldots,F_N)$ subordinate to $\mathcal U$. The maximum tests all linear combinations with coefficients in $[-1,1]$ at once\rev{. In this way, the invariant fits into the function-theoretic approach} initiated in \cite{EPZ} and developed in \cite{Polterovich-noise}. \rev{A positive lower bound for \keepmath{$pb_\omega(\mathcal U)$} means that the pairwise Poisson brackets cannot all be made arbitrarily small in the uniform norm by varying the subordinate partition.} \rev{Such a lower bound can, for instance, be obtained in terms of symplectic spectral width; see \cite[Section 9.2]{Polterovich-Rosen}.}

\medskip

Fix a closed co-oriented contact manifold $(Y,\xi=\ker\alpha)$ of dimension $2n+1$, together with the contact form $\alpha$. Its Reeb vector field is characterized by $\alpha(R_\alpha)=1$ and $\iota_{R_\alpha}d\alpha=0$. For a contact Hamiltonian function $f\in C^\infty(Y)$, the contact Hamiltonian vector field $X_f$ is uniquely determined by the following equations:
\begin{equation} \label{eq:contact-hamiltonian}
\alpha(X_f)=f,\qquad \iota_{X_f}d\alpha=\rrev{(R_\alpha f)}\alpha-df.
\end{equation}
Its flow is denoted by $\varphi_f^t$. The constant Hamiltonian $1$ generates the Reeb vector field, so $\varphi_1^t$ denotes its flow. The contact Poisson bracket (with respect to the contact 1-form $\alpha$) is defined as follows\footnote{Our sign convention is the negative of the convention in Rosen--Zhang \cite[Section 2]{Rosen-Zhang}, \rev{but the vanishing of brackets and their \keepmath{$C^0$}-norms are unaffected.}}:
\begin{equation} \label{dfn-CBP}
\{f,g\}_\alpha:=\rrev{X_f g}-g\rrev{R_\alpha f}.
\end{equation}
It is skew-symmetric and satisfies $[X_f,X_g]=X_{\{f,g\}_\alpha}$. \rev{For the definition in (\ref{dfn-CBP}) and its relation to the symplectic bracket on the symplectization,} see McDuff--Salamon \cite[Exercise 3.57(iv)]{McDuff-Salamon}. All functions in the statements below are smooth, and $\rrev{\norm{h}_\infty=\max_Y|h|}$.

\rev{Two differences from the symplectic bracket are central to this paper.} First, constant functions cannot be overlooked. In fact, 
\begin{equation} \label{eq:constant-bracket}
\{1,h\}_\alpha=R_\alpha h, 
\end{equation}
\rev{because the constant function one generates the Reeb flow.} Second, contact Hamiltonian flows preserve $\xi$ but in general change $\alpha$. For a contactomorphism $\varphi$, write
\begin{equation} \label{dfn-CF}
\varphi^*\alpha=e^{\kappa(\varphi)}\alpha\text{\rev{,}}
\end{equation}
where $\kappa(\varphi)$ is the (logarithmic) conformal factor of $\varphi$ (the same notation applies to contact isotopies). Along an autonomous flow, $\frac{d}{dt}\kappa(\varphi_f^t)=\rrev{(R_\alpha f)}\circ\varphi_f^t$. Therefore, $R_\alpha f=0$ precisely when the flow is strict, meaning that it preserves $\alpha$.

Throughout the paper, displaceability means displaceability by an element of $\Cont_0(Y,\xi)$. Thus a subset $K\subset Y$ is contact displaceable if some $\varphi\in\Cont_0(Y,\xi)$ satisfies $\varphi(K)\cap K=\varnothing$. The contact form is fixed in all related norms and bracket invariants. 

\subsection{Three forms of contact \texorpdfstring{$C^0$}{C0}-rigidity}\label{subsec:rigidity}

\rev{Our first result states that contact Poisson commutation is preserved under uniform convergence.} It requires no control of derivatives of the approximating functions.

\begin{mainthm} \label{prop-1} Let $(Y,\xi=\ker\alpha)$ be a closed contact manifold, and let $f_j,g_j,f,g\in C^\infty(Y)$ satisfy $f_j\to f$ and $g_j\to g$ uniformly. If $\{f_j,g_j\}_\alpha=0$ for every $j$, then $\{f,g\}_\alpha=0$. \end{mainthm}

The proof of Theorem~\ref{prop-1} is given in Section~\ref{sec:proof-A}, and \rev{the main input from contact geometry is the non-degeneracy of Shelukhin's norm} \cite[Theorem A]{Shelukhin-cont}. \rev{Next, we relax exact commutation and consider approximating functions that form almost commuting pairs. In the symplectic setting, the rigidity result in \cite[Section 3]{CV} allows the approximating brackets to tend uniformly to zero (see \cite[Definition 1.1 and Theorem 1.2]{CV}). However, a direct analogue in contact geometry fails, as the following example shows.} 

\begin{example} \label{ex-counter} \rev{Equip \keepmath{$\T^3=(\R/2\pi\Z)^3$} with the standard contact 1-form \keepmath{$\alpha=\cos x\,dz+\sin x\,dy$}. Consider} 
\begin{equation} \label{ce:limit-torus}
f=\cos x\cos z,\qquad g=\cos x\sin z.
\end{equation}
\rev{To construct the approximating sequences, introduce the auxiliary functions \keepmath{$A_j(x) :=\sqrt{4j(\cos x)^2-2\sin x \cos x+j^{-1}}$} for every \keepmath{$j \in \N$}, and set}
\[ f_j:=f+u_j \,\,\,\,\mbox{where} \,\,\,\, u_j =j^{-2}A_j(x)\cos(j^4y) \]
and similarly, 
\[ g_j:=g+v_j \,\,\,\,\mbox{where} \,\,\,\, v_j=j^{-2}A_j(x)\sin(j^4y). \]
Both $f_j$ and $g_j$ are well-defined on $\T^3$ for any $j \in \N$, and $f_j \to f$, $g_j \to g$ uniformly. \rev{A direct computation gives} 
\[ \norm{\{f,g\}_\alpha}_\infty=1\,\,\,\,\mbox{but $\norm{\{f_j,g_j\}_\alpha}_\infty \to 0$ as $j \to \infty$}. \]
The computation is given in Section~\ref{sec:counterexample}. 
\end{example} 

\rev{Thus, unlike Theorem \ref{prop-1}, this weaker form of contact Poisson bracket rigidity requires an additional hypothesis. In our next result, we assume that the two functions \keepmath{$f, g$} have no common zeros on the ambient contact manifold \keepmath{$(Y, \ker \alpha)$}; equivalently, \keepmath{$f^2 + g^2 >0$} on \keepmath{$Y$}. By comparison, the functions \keepmath{$f, g$} in (\ref{ce:limit-torus}) in Example \ref{ex-counter} do have common zeros in \keepmath{$\T^3$}.} \rev{The absence of common zeros is useful because it allows a change of variables.} \rev{The two-term expression for the contact Poisson bracket \keepmath{$\{f,g\}_\alpha$} in (\ref{dfn-CBP}) can then be related to the derivative of the circle-valued function \keepmath{$(f+\sqrt{-1}g)/\sqrt{f^2+g^2}$} along the contact flow generated by \keepmath{$\sqrt{f^2+g^2}$} (see Lemma~\ref{lem-1}). The resulting expression is easier to control in the \keepmath{$C^0$}-norm, leading to the following weak rigidity statement.}

\begin{mainthm} \label{prop-2} Let $(Y,\xi=\ker\alpha)$ be a closed contact manifold, and let $f_j,g_j,f,g\in C^\infty(Y)$ satisfy $f_j\to f$ and $g_j\to g$ uniformly. If $f^2+g^2>0$ on $Y$ and $\norm{\{f_j,g_j\}_\alpha}_\infty\to0$, then $\{f,g\}_\alpha=0$. \end{mainthm}

The proof of Theorem \ref{prop-2} is given in Section~\ref{sec:proof-B}. \rev{Besides the polar identity in Lemma~\ref{lem-1}, a crucial ingredient is the lower semicontinuity of \keepmath{$\norm{\rrev{X_f g}}_\infty$}, established by adapting the proof of \cite[Theorem 1.1.2]{Buhovsky-C0}. The key point is that a decrease in the directional derivative would force a flow to displace a fixed open set, contradicting the contact energy-capacity inequality in \cite[Corollary 15]{Shelukhin-cont}.} 

\medskip

\rev{Our third result strengthens the conclusion to a lower bound analogous to \eqref{eq:symplectic-lsc}. It allows common zeros of the two input contact Hamiltonian functions, but in exchange requires control of the interaction} between the uniform approximation error and derivatives along the Reeb vector field. \rev{Indeed, Example \ref{ex-counter} shows that an unconditional contact analogue of (\ref{eq:symplectic-lsc}) fails. More specifically, the computations in Section \ref{sec:counterexample} show that its approximating sequences violate condition (\ref{Reeb-cond}) in the following result.} 

\begin{mainthm} \label{prop-3} Let $(Y,\xi=\ker\alpha)$ be a closed contact manifold, and let $f_j,g_j,f,g\in C^\infty(Y)$ satisfy $f_j\to f$ and $g_j\to g$ uniformly. If
\begin{equation} \label{Reeb-cond}
\max\{\norm{f_j-f}_\infty,\norm{g_j-g}_\infty\} \cdot \norm{R_\alpha g_j}_\infty\longrightarrow0,
\end{equation}
then $\norm{\{f,g\}_\alpha}_\infty \le\liminf_{j\to\infty}\norm{\{f_j,g_j\}_\alpha}_\infty$.
\end{mainthm}

Note that the condition~\eqref{Reeb-cond} holds whenever the functions $R_\alpha g_j$ are uniformly bounded, and in particular when every $g_j$ is Reeb-invariant. \rev{It also allows \keepmath{$\norm{R_\alpha g_j}_\infty$} to grow, provided that it grows more slowly than the reciprocal of the approximation error.} Moreover, by skew-symmetry of the contact Poisson bracket, the condition with $f_j$ in place of $g_j$ is equally sufficient.

\rev{For comparison, the symplectic lower-semicontinuity statement \eqref{eq:symplectic-lsc} was established in \cite{Zapolsky-surfaces} for surfaces and in \cite{EP-C0} for general symplectic manifolds. The latter proof uses the length-minimizing property of short Hamiltonian paths in Hofer geometry. An alternative proof in \cite{Buhovsky-C0} uses displacement energy and also yields quantitative} $2/3$\rev{-rate estimates.} Our proof of Theorem \ref{prop-3}, given in Section~\ref{sec:proof-C}, follows this local displacement method, \rev{using the energy-capacity inequality in \cite[Corollary 15]{Shelukhin-cont} to obtain a positive lower bound for the relevant displacement energy.} 

\begin{remark}[\rev{Quantitative control in (\ref{Reeb-cond})}]  \label{rmk-conformal-factor} \rev{The condition (\ref{Reeb-cond}) provides a form of control over conformal factors.} This also appears in topological contact dynamics \cite{Muller-Spaeth} and in contact rigidity for uniform limits of ambient maps \cite{Usher-contact}. Recall that rigidity of smooth images of coisotropic and Legendrian submanifolds under such limits is studied in \cite{HLS,DRS,Rosen-Zhang,Usher-contact}. \rev{The author began investigating an unconditional version of the \keepmath{$C^0$}-rigidity statement in Theorem \ref{prop-3} shortly after completing \cite{Rosen-Zhang}. However, this investigation was set aside because the conformal factor arising from terms such as \keepmath{$\norm{R_\alpha g_j}_\infty$} seemed difficult to eliminate. More recently, the contact rigidity results in \cite{DUZ}, particularly the Reeb-invariance hypothesis in the contact big fiber theorem, suggested imposing quantitative control on \keepmath{$\norm{R_\alpha g_j}_\infty$} instead of eliminating it. In fact, Example \ref{ex-counter} confirms that some additional control is necessary in general. Thus this observation motivates the present paper.}    \end{remark}

\begin{remark}\label{rmk-change-contact-form}
Theorems~A, B, and~C hold for every fixed contact 1-form. Explicitly, if $\widehat\alpha=e^\sigma\alpha$ for some fixed smooth function $\sigma: Y \to \R$, then
\[
\begin{aligned}
\{e^\sigma f,e^\sigma g\}_{\widehat\alpha}&=e^\sigma\{f,g\}_\alpha, \qquad R_{\widehat\alpha}(e^\sigma g)&=R_\alpha g+X_g^\alpha\sigma
\end{aligned}
\]
where $X_g^{\alpha}$ denotes the contact Hamiltonian vector field of $g$ with respect to $\alpha$. Since $Y$ is compact, the rescaling $f\mapsto e^\sigma f$ preserves uniform convergence and the absence of common zeros. Therefore, Theorems~A and~B transfer directly. If the original sequence satisfies the hypothesis of Theorem~C, then the rescaled sequence satisfies the corresponding hypothesis precisely when
\[
\max\{\|f_j-f\|_\infty,\|g_j-g\|_\infty\}
\|X_{g_j}^\alpha\sigma\|_\infty\longrightarrow0.
\]
This is automatic for constant $\sigma$. More generally, it suffices that the same approximation error multiplied by $\|g_j\|_{C^1}$ tends to zero.
\end{remark}

\subsection{The contact Poisson bracket invariant of a cover}\label{subsec:covers}

Motivated by the symplectic cover invariant \eqref{eq:symplectic-cover}, we introduce a contact bracket functional for partitions of unity and the associated invariant of an open cover. The definition uses arbitrary smooth partition functions, without imposing Reeb-invariance.

\begin{definition}[Contact Poisson bracket invariant of a cover]\label{def:contact-pb-cover} Let $(Y,\xi=\ker\alpha)$ be a closed contact manifold. For a smooth partition of unity $\vec f=(f_1,\ldots,f_N)$ on $Y$, with $f_i\ge0$ and $\sum_i f_i=1$, define
\begin{equation} \label{dfn-Pf}
\Palpha(\vec f):= \max_{a,b\in[-1,1]^N} \rrev{\norm{\left\{\sum_i a_i f_i,\sum_j b_j f_j\right\}_\alpha}_\infty}.
\end{equation}
For a finite open cover $\mathcal U=\{U_1,\ldots,U_N\}$ of $Y$, its contact Poisson bracket invariant is
\begin{equation} \label{eq:contact-pb-cover}
pb_\alpha(\mathcal U):=\inf_{\vec f}\Palpha(\vec f),
\end{equation}
where the infimum is taken over smooth partitions of unity subordinate to $\mathcal U$, meaning $\supp(f_i)\subset U_i$ for every $i$. \end{definition}

Since $\sum_i f_i=1$ and $\{1,h\}_\alpha=R_\alpha h$ as in (\ref{eq:constant-bracket}), choosing $a=(1,\ldots,1)$ in \eqref{dfn-Pf} gives, for every $b\in[-1,1]^N$,
\begin{equation} \label{eq:partition-reeb-control}
\rrev{\norm{R_\alpha\left(\sum_j b_jf_j\right)}_\infty} \le\Palpha(\vec f).
\end{equation}
Thus the same functional measures the brackets and controls the Reeb derivatives of the partition functions. \rev{Unlike in symplectic geometry, even for a partition with two members, \keepmath{$\{f,1-f\}_\alpha=-R_\alpha f$} need not vanish. This explains the modification introduced later in Definition \ref{def:contact-pb-subsets}.} Along any subsequence for which $\Palpha(\vec f_j)$ is bounded, \eqref{eq:partition-reeb-control} supplies the Reeb bound needed in Theorem~\ref{prop-3}. This gives the following lower-semicontinuity statement.

\begin{mainthm}\label{cor-1} Let $\vec f_j=(f_{j,1},\ldots,f_{j,N})$ and $\vec f=(f_1,\ldots,f_N)$ be smooth partitions of unity on a closed contact manifold $(Y,\xi=\ker\alpha)$. If $f_{j,i}\to f_i$ uniformly for every $i$, then
\[ \Palpha(\vec f)\le\liminf_{j\to\infty}\Palpha(\vec f_j). \]
\end{mainthm}
The proof of Theorem \ref{cor-1} is given in Section~\ref{sec:proof-D}. \rev{Having established lower semicontinuity, we turn to two positive lower bounds for the cover invariant} in Definition~\ref{def:contact-pb-cover}, using quasi-morphisms and contact spectral invariants, respectively. Their hypotheses will be stated separately.

\medskip

Write $\widetilde\varphi_h$ for the class of the contact isotopy generated by $h$ in $\widetilde{\Cont}_0(Y,\xi)$. For an autonomous $f$, we also write $\widetilde\varphi_f^{\,t}=[s\mapsto\varphi_f^{ts}]$, $0\le s\le1$. In particular, $\widetilde\varphi_1^{\,t}$ is the class of the time-$t$ Reeb isotopy. Recall that a homogeneous quasi-morphism $\mu:\widetilde{\Cont}_0(Y,\xi)\to\R$ satisfies $\mu(a^m)=m\mu(a)$ for every $m\in\Z$ and admits a defect bound $D(\mu)>0$ such that
\begin{equation} \label{eq:mu-defect}
|\mu(ab)-\mu(a)-\mu(b)|\le D(\mu)\qquad\text{for all }a,b\in\widetilde{\Cont}_0(Y,\xi).
\end{equation}
It is monotone if $\mu(a)\le\mu(b)$ whenever $a^{-1}b$ has a representative generated by a nonnegative contact Hamiltonian. It has the vanishing property if $\mu(\widetilde\varphi_h)=0$ whenever $h$ is supported in $[0,1]\times U$ for a contact displaceable open set $U$. A nonzero monotone homogeneous quasi-morphism has $\mu(\widetilde\varphi_1)>0$. Therefore, for the estimates involving $\mu$, we may use the normalization
\begin{equation} \label{eq:mu-normalization}
\mu(\widetilde\varphi_1)=1,
\end{equation}
and $D(\mu)$ always denotes a defect bound after this normalization. \rev{Remark \ref{rmk-normalization} explains in detail how to achieve the normalization in (\ref{eq:mu-normalization}).} 

\rev{The construction in \cite{Givental}, as further developed by Borman--Zapolsky \cite{BZ}, supplies} such quasi-morphisms on standard contact real projective spaces and on further prequantizations obtained by contact reduction. In particular, \cite[Proposition 1.2]{BZ} establishes monotonicity and the vanishing property for the projective-space quasi-morphism. \rev{Related constructions on lens spaces are given in \cite{GKPS}.}

\begin{mainthm}\label{prop-4} Suppose that $(Y,\xi=\ker\alpha)$ admits a nonzero monotone homogeneous quasi-morphism $\mu$ with the vanishing property, normalized by \eqref{eq:mu-normalization}. For every finite open cover $\mathcal U=\{U_1,\ldots,U_N\}$ of $Y$ by contact displaceable sets, with $N\ge2$, we have
\begin{equation} \label{eq:main-cover-bound}
pb_\alpha(\mathcal U)\ge\frac{1}{16(N-1)^2D(\mu)}>0.
\end{equation}
\end{mainthm}

\rev{The proof follows the method in \cite[Theorems 1.4 and 1.7]{EPZ}. First, we compare the flow of a sum of input Hamiltonian functions with the product of their flows. We then use the quasi-morphism defect} $D(\mu)$ \rev{to bound the failure of additivity. Finally, we sum these estimates over successive partial sums of a partition.} The flow comparison and its application to cover bounds are explained in \cite[Sections 4.6 and 9.2]{Polterovich-Rosen}. \rev{In contact geometry, however, an additional estimate is needed to control the Reeb derivatives. This is precisely where inequality~\eqref{eq:partition-reeb-control} enters: it bounds these derivatives by \keepmath{$\Palpha(\vec f)$}, allowing them to be absorbed into the left-hand side of (\ref{eq:main-cover-bound}).} The proof of Theorem \ref{prop-4} is given in Section~\ref{sec:proof-E}.

\subsection{A lower bound from contact spectral invariants}\label{subsec:spectral-cover}
\rev{As an alternative to quasi-morphisms, contact spectral invariants give a lower bound in terms of displacement energy.} \rev{Several constructions of these invariants are available \cite{Sandon,Albers-Merry,Dylan-eternal,DUZ}. For our purposes, we use the recent construction in \cite{CKZ}.} \rev{To state its hypotheses, recall that an aspherical} $G$-filling is a convex-at-infinity symplectic manifold $(W,\omega)$ with $\omega|_{\pi_2(W)}=0$ and a symplectic action of a finite group $G$, commuting with the Liouville flow near infinity, such that $G$ acts freely on the ideal contact boundary and $Y=\partial_\infty W/G$; see \cite[Definition 1.1]{CKZ}. If some point of $W$ has nontrivial stabilizer, \cite[Theorem 1.10]{CKZ} constructs, for every contact form $\alpha$ on $Y$, a real-valued spectral invariant
\[ c_\alpha:\widetilde{\Cont}_0(Y,\xi)\longrightarrow\R. \]
The construction uses equivariant contact Floer cohomology and does not require the ordinary symplectic cohomology of $W$ to vanish.

For a contactomorphism $\varphi$, its translated spectrum is
\begin{equation} \label{eq:translated-spectrum}
\operatorname{Spec}_\alpha(\varphi):=\{s\in\R \,| \, \varphi(x)=\mrev{\varphi_1^s}(x),\ \kappa(\varphi)(x)=0\text{ for some }x\in Y\}.
\end{equation}
The properties (R1)--(R5) in \cite[Section 1.1]{CKZ} give spectrality, monotonicity under nonnegative contact deformations, continuity from above, normalization, and subadditivity. In particular,
\begin{equation} \label{eq:spectral-properties}
c_\alpha(\widetilde\varphi)\in\operatorname{Spec}_\alpha(\varphi),\qquad c_\alpha(\widetilde{\mathrm{id}})=0,\qquad c_\alpha(ab)\le c_\alpha(a)+c_\alpha(b),
\end{equation}
and $c_\alpha(a)=\lim_{\rrev{\epsilon} \to 0^+}c_\alpha(\widetilde\varphi_1^{\,\rrev{\epsilon}}a)$. \rev{These properties determine the exact value on the class represented by the Reeb flow: \keepmath{$c_\alpha(\widetilde\varphi_1^{\,t})=t$}, as verified in Lemma~\ref{lem:spectral-calibration}.}

For an open set $U\subset Y$, its contact displacement energy is
\begin{equation} \label{eq:intro-displacement-energy}
e_\alpha(U):=\inf_{\psi_1(U)\cap U=\varnothing}\int_0^1\norm{h_s}_\infty\,ds,
\end{equation}
where $h_s$ generates the based contact isotopy $\psi_s$. \rev{For a displaceable open set, this number is finite, and it is positive when \keepmath{$U$} is nonempty by \cite{Shelukhin-cont}.}

\begin{mainthm}\label{thm:spectral-cover} Suppose that $(Y,\xi)$ admits an aspherical $G$-filling with a point having nontrivial stabilizer. For every contact form $\alpha$ and every finite open cover $\mathcal U=\{U_1,\ldots,U_N\}$ of $Y$ by contact displaceable sets, we have
\begin{equation} \label{eq:main-spectral-cover}
pb_\alpha(\mathcal U)\ge\frac{1}{16(N-1)\displaystyle\sum_{i=1}^N e_\alpha(U_i)}>0.
\end{equation}
\end{mainthm}

\begin{remark} \rev{The bound in Theorem \ref{thm:spectral-cover} involves the sum of the displacement energies of the cover members.} Since $\sum_i e_\alpha(U_i)\le N\max_i e_\alpha(U_i)$, we obtain the weaker bound $pb_\alpha(\mathcal U)\ge (16N(N-1)\max_i e_\alpha(U_i))^{-1}$. Note that this bound depends on $N$ (in order $2$), the number of cover members in $\mathcal U$. 
\end{remark}

The proof of Theorem~\ref{thm:spectral-cover} is given in Section \ref{sec:proof-F}. \rev{It uses only the axiomatic properties of the spectral invariant \keepmath{$c_\alpha$} and does not require a quasi-morphism or a periodic Reeb flow. The same conclusion therefore holds in any setting admitting an invariant \keepmath{$c_{\alpha}$} with these properties.} \rev{The following example illustrates the scope of Theorem \ref{thm:spectral-cover}: all the manifolds listed have infinite fundamental group} $\pi_1$. \rev{By comparison, the real projective and lens-space examples} in \cite{Givental,GKPS}, as well as the monotone toric prequantizations \cite{BZ}, all have finite fundamental group. \rev{Thus the cited quasi-morphism constructions do not directly establish the hypotheses of Theorem~\ref{prop-4} for these examples.} 

\begin{example}\label{ex:ckz-examples} (i) Let $\Sigma$ be a closed oriented surface of positive genus with an integral area form $\omega$, and let $L\to\Sigma$ be a complex line bundle with $c_1(L)=-[\omega]$. For any integer $m\ge2$, the unit circle bundle $Y=S(L^{\otimes m})$ satisfies the hypotheses of Theorem~\ref{thm:spectral-cover}. Indeed, equip the total space of $L$ with its usual convex-at-infinity symplectic structure. Since it retracts onto $\Sigma$, its second homotopy group vanishes. Fiberwise multiplication by $e^{2\pi\sqrt{-1}/m}$ fixes the zero section and acts freely on the ideal boundary. The quotient boundary is $S(L^{\otimes m})$, as in \cite[Section 1.2.2]{CKZ}. Therefore, \eqref{eq:main-spectral-cover} holds for every contact form defining this contact structure. Note that the projection $Y\to\Sigma$ induces a surjection on fundamental groups, so $\pi_1(Y)$ is infinite.

(ii) Let $W=T^*\T^2$, where $\T^2=\R^2/\Z^2$, with its canonical Liouville form. The involution $(q,p)\mapsto(-q,-p)$ preserves this form and has exactly four fixed points, all on the zero section. It acts freely on $S^*\T^2$, since a fixed covector would have to satisfy $p=-p$ and hence $p=0$. Therefore,
\[ Y=S^*\T^2/\langle(q,p)\mapsto(-q,-p)\rangle \]
admits an exact, and hence aspherical, $\Z/2\Z$-filling with nontrivial stabilizers. \rev{This satisfies the hypotheses of \cite[Theorem 1.10]{CKZ}}, and Theorem~\ref{thm:spectral-cover} applies for every contact form. Since $S^*\T^2\cong\T^3$ is a finite cover of $Y$, the group $\pi_1(Y)$ contains $\Z^3$ as a finite-index subgroup and is infinite.
\end{example}

\rev{A crucial step in the proof of Theorem \ref{thm:spectral-cover} is the following energy estimate, valid for every \keepmath{$t \geq 0$} (see Proposition \ref{prop:local-spectral-energy}):} 
\begin{equation} \label{contact-Usher-energy}
0\le c_\alpha(\widetilde\varphi_h^{\,t})\le e_\alpha(U)\bigl(1+e^{t\norm{R_\alpha h}_\infty}\bigr)
\end{equation}
for any contact Hamiltonian function $h\ge0$ compactly supported in a contact displaceable open set $U \subset (Y, \ker \alpha)$. Its proof adapts the spectral displacement argument of Frauenfelder--Ginzburg--Schlenk \cite[Proposition 2.1]{FGS} and Usher \cite[Proposition 3.1]{Usher-energy}. \rev{To adapt that argument to the contact setting, we use Reeb translation times as the spectral values, as in (\ref{eq:translated-spectrum}).} \rev{Compared with the symplectic analogue,} the extra term $e^{t\norm{R_\alpha h}_\infty}$ in (\ref{contact-Usher-energy}) comes from a cut-off trick (see the beginning of the proof of Proposition \ref{prop:local-spectral-energy}). \rev{More precisely, an appropriate cut-off of the original} $h$ \rev{ensures that its support and a prescribed Reeb translate remain in} $U$, \rev{which allows us to exclude one spectral value throughout the deformation.} 

\rev{To conclude this subsection, we compare the displacement-energy bound in Theorem~\ref{thm:spectral-cover} with the symplectic results in \cite[Theorem 4.8 and Sections 5.2--5.4]{Polterovich-noise} and those of} Seyfaddini \cite{Seyfaddini}, Tanny \cite{Tanny}, Buhovsky--Logunov--Tanny \cite{BLT}, Payette \cite{Payette}, and Shi--Lu \cite{Shi-Lu}. \rev{These results concern bounds of the form} 
\[ pb_\omega(\mathcal U)\ge\frac{C}{\max_i e_\omega(U_i)}, \]
\rev{where \keepmath{$e_\omega$} denotes Hofer displacement energy and the dependence of the constant \keepmath{$C>0$} varies with the hypotheses.} \rev{In particular, \cite[Corollary 1.9]{BLT} establishes the conjecture in \cite[Question 8.1]{Polterovich-noise} in dimension two, with a constant \keepmath{$C$} independent of both the number of sets and their overlap \cite[Remark 1.10 and Example 3.7]{BLT}.} \rev{The proof uses two-dimensional methods, namely intersections of level curves and area estimates.} Therefore, one may ask whether, in the setting of Theorem \ref{thm:spectral-cover}, \rev{the estimate (\ref{eq:main-spectral-cover}) can be improved to} 
\[ pb_\alpha(\mathcal U)\ge\frac{C}{\max_i e_\alpha(U_i)} \]
\rev{where the constant \keepmath{$C>0$} depends only on \keepmath{$(Y,\alpha)$},} independently of the number and overlap of the cover members. In dimension three, a possible approach is to adapt the method of \cite{BLT} to surfaces transverse to the Reeb flow. 

Also, imitating the proof of \cite[Theorem 1.15]{BLT}, one can adapt the main steps in the proof of Theorem~\ref{prop-3} (say \eqref{uig} and \eqref{Gronwall}) to prove that if $(Y,\ker\alpha)$ has dimension at least three and is equipped with a fixed Riemannian metric $\rho$, and every member of a given cover $\mathcal U$ has diameter at most $\ep$ with respect to $\rho$, where $0<\ep\le\ep_0$, then we have an estimate $pb_{\alpha}(\mathcal U)\geq\frac{c}{\ep^2}$. Both $\ep_0>0$ and $c>0$ depend only on $\alpha$ and $\rho$. We will provide more details and develop this cover-cardinality-independence phenomenon more systematically in future work.

\subsection{A quantitative contact big fiber theorem}\label{subsec:big-fiber}
\rev{The study of rigidity for an involutive system \keepmath{$(f_1, ..., \mrev{f_k}): Y \to \R^k$} goes back to the symplectic big fiber theorem in \cite[Theorem 2.1]{EP} and has recently appeared in contact geometry in \cite[Theorem 1.4]{DUZ}. For the contact result in \cite{DUZ}, involutivity requires a contact form for which all components are Reeb-invariant and all their pairwise contact Poisson brackets vanish.} \rev{The usual argument pulls back a partition of unity from the target, assuming that all the relevant inverse images are displaceable, and then derives a contradiction from the normalization and vanishing properties of an invariant. In our setting, the lower bound in either Theorem \ref{prop-4} or Theorem \ref{thm:spectral-cover} applies to the same partition of unity.} \rev{By combining this construction with quantitative control inspired by (\ref{Reeb-cond}) in Theorem \ref{prop-3}, we obtain} the following quantitative version of the contact big fiber result in \cite{DUZ}. \rev{Here, quantitative means that the conclusion concerns a thickened fiber and applies to systems whose component functions need not Poisson commute or be Reeb-invariant.} 

\rrev{In this subsection and Section~\ref{sec:proof-G}, we use the maximum norm $|v|=\max_i|v_i|$ on $\R^k$. The Euclidean-norm version follows by replacing $c_k$ with $k^{-(k+1)}c_k$.}

\begin{mainthm}\label{thm-qcbft} Suppose that $(Y,\xi=\ker\alpha)$ admits a nonzero monotone homogeneous quasi-morphism $\mu$ with the vanishing property, normalized by \eqref{eq:mu-normalization}, with defect $D(\mu)$. Then for each $k\ge1$ there is a constant $c_k>0$, depending only on $k$, with the following property. Let $F=(f_1,\ldots,f_k):Y\to\R^k$ be smooth, and write $M=\max_i\norm{f_i}_\infty>0$. If $0<\delta\le M$ and
\begin{equation} \label{eq:main-big-fiber}
D(\mu)\max\left\{M\max_i\norm{R_\alpha f_i}_\infty,\ \max_{i,j}\norm{\{f_i,f_j\}_\alpha}_\infty\right\}<c_k\frac{\delta^{2k+2}}{M^{2k}},
\end{equation}
then $\{x\in Y:|F(x)-a|\le\delta\}$ is contact non-displaceable for some $a\in\R^k$. \end{mainthm}

When all contact brackets and Reeb derivatives vanish, the left-hand side of \eqref{eq:main-big-fiber} is zero. Therefore, letting the radius $\delta$ tend to zero gives the following consequence.

\begin{cor}\label{cor:big-fiber} Under the hypotheses of Theorem~\ref{thm-qcbft}, let $F=(f_1,\ldots,f_k):Y\to\R^k$ be smooth. If $R_\alpha f_i=0$ and $\{f_i,f_j\}_\alpha=0$ for every $i,j$, then $F$ has a contact non-displaceable fiber. \end{cor}

\rev{For context, qualitative conclusions of the type in Corollary \ref{cor:big-fiber} have been obtained under different global hypotheses. For example, \cite[Theorem 1.17]{BZ} proves non-displaceability of a toric fiber in suitable prequantizations using a contact quasi-morphism.} \rev{Other approaches use symplectic cohomology with support and its Mayer--Vietoris property in \cite{SUV}, or partial contact quasi-states constructed from contact spectral invariants in \cite{DUZ}.} 

\begin{remark} \rev{We state Theorem \ref{thm-qcbft} and Corollary \ref{cor:big-fiber} under the hypothesis of Theorem \ref{prop-4}. Under the hypothesis of Theorem~\ref{thm:spectral-cover}, the same cover argument gives the qualitative conclusion of Corollary~\ref{cor:big-fiber}. A quantitative spectral version of this argument also requires control of the sum of the displacement energies of the sets in the refined cover.} \end{remark} 

\begin{remark} \rev{While this paper was being written, Adi Dickstein generously shared a draft \cite{Dickstein-thick-fiber}, establishing a quantitative refinement of the qualitative {\rm symplectic} big fiber theorem in \cite{EP}. His proof shares the main idea of the proof of Theorem \ref{thm-qcbft}: refining the displaceable preimages.} We emphasize that both proofs are influenced by \cite{EPZ}. \end{remark}  

The following example illustrates an application of Theorem~\ref{thm-qcbft}. 

\begin{example} \label{ex-ellipsoid-cont-non-dis} We recall a standard fact. Consider the boundary of an ellipsoid in $\C^2$, which is contactomorphic to the standard contact sphere. Then every proper closed subset is contact displaceable by \cite[Proposition 4.1]{MP16}. However, its antipodal quotient is contactomorphic to the standard contact $\R P^3$, which admits Givental's monotone homogeneous quasi-morphism with the vanishing property \cite[Proposition 1.2]{BZ}. This supplies the hypothesis of Theorem~\ref{thm-qcbft}.

Explicitly, consider the contact boundary of an irrational ellipsoid
\[
\Sigma:=\left\{(z_0,z_1)\in\C^2 \,\bigg| \, |z_0|^2+\frac{|z_1|^2}{\sqrt2}=1\right\} \quad \mbox{and} \quad Y:=\Sigma/\{z\sim-z\}
\]
where $z_\ell=x_\ell+\sqrt{-1}\,y_\ell$. On $Y$ take the form descended from $\sum_{\ell=0}^1(x_\ell\,dy_\ell-y_\ell\,dx_\ell)|_\Sigma$. Its Reeb field is $\partial_{\theta_0}+2^{-1/2}\partial_{\theta_1}$ away from the coordinate axes. Consider the function 
\[ u([z]):=\frac{|z_0|^2}{|z_0|^2+|z_1|^2}\in[0,1]. \]
Then $R_\alpha u=0$ and Corollary~\ref{cor:big-fiber} guarantees the existence of a point $a \in [0,1]$ such that $u^{-1}(a)$ is contact non-displaceable. In fact, it is easy to see that $a=1/2$ is the only such value. 

Now, let us consider a perturbation of $u$ as follows:
\begin{equation}\label{eq:function}
f:=u+\frac{\chi(u)}{16} \left(\cos\bigl(2(q\theta_0-p\theta_1)\bigr) +\cos\bigl(2(q'\theta_0-p'\theta_1)\bigr)\right),
\end{equation}
where $\theta_0, \theta_1$ denote the angular coordinates of $z_0, z_1$, respectively, and $q, p, q', p'$ are four fixed positive integers defined by
\[
p+q\sqrt2:=(1+\sqrt2)^{32}, \qquad p'+q'\sqrt2:=(1+\sqrt2)^{33},
\]
and $\chi$ is a smooth cutoff $0\le\chi\le1$ which is one on $[3/8,5/8]$ and zero outside $(1/4,3/4)$. Then one can check that $\|f-u\|_{\infty} = 1/8$, $\|f\|_{\infty} = 1$,  
\begin{equation} \label{ex:Reeb-bound}
f\left(u^{-1}\left(1/2\right)\right) = [3/8, 5/8] \quad \mbox{and} \quad \|R_{\alpha} f\|_{\infty} = \frac18(\sqrt2-1)^{32}.
\end{equation}
Since $R_\alpha f\not\equiv0$, Corollary~\ref{cor:big-fiber} does not apply. However, Theorem~\ref{thm-qcbft}, with $\delta=1/16$, yields a non-displaceable thickened fiber $\{|f-b|\le1/16\}$. Here, we can take $D(\mu)=64$ (by \cite[Theorem 1.2]{GKPS}) and $c_1=1/77760$ (following the proof of Theorem~\ref{thm-qcbft}). We emphasize that the parameters in this perturbation were carefully chosen so that the subset $\{|f-b|\le1/16\}$ does {\it not} contain the (non-displaceable) fiber $u^{-1}(\frac{1}{2})$. Otherwise, $f(u^{-1}(\frac{1}{2}))$ would be contained in an interval of length $1/8$, contradicting the first equality in (\ref{ex:Reeb-bound}). 
\end{example}

\subsection{Contact pb invariants}\label{subsec:pb34}

\rev{We now introduce contact invariants of triples and quadruples of compact subsets.} The definitions are motivated by the symplectic invariants $pb_3$ and $pb_4$ \rev{introduced in \cite[Section 1.2]{BEP}.} \rev{These symplectic invariants provide a powerful way to relate the rigidity of subsets to Hamiltonian dynamics.} 


\begin{definition}[Contact invariants of triples and quadruples]\label{def:contact-pb-subsets} (i) Let $X_0,X_1,X_2\subset Y$ be compact with $X_0\cap X_1\cap X_2=\varnothing$. An admissible triple is a smooth partition of unity $\vec h=(h_0,h_1,h_2)$ such that $h_i$ vanishes on a neighborhood of $X_i$. Define
\begin{equation} \label{eq:subset-triple-definition}
pb_{3,\alpha}(X_0,X_1,X_2):=\inf_{\vec h}\max_{0\le i<j\le2}\norm{\{h_i,h_j\}_\alpha}_\infty,
\end{equation}
where the infimum is taken over all admissible triples.

(ii) For compact sets $X_0,X_1,Y_0,Y_1\subset Y$ with $X_0\cap X_1=Y_0\cap Y_1=\varnothing$, an admissible pair consists of smooth functions $f,g:Y\to[0,1]$ satisfying
\[ f=i\text{ on a neighborhood of }X_i,\qquad g=j\text{ on a neighborhood of }Y_j,\qquad i,j\in\{0,1\}. \]
Define
\begin{equation} \label{eq:subset-square-definition}
pb_{4,\alpha}(X_0,X_1,Y_0,Y_1):=\inf_{(f,g)}\max_{\substack{u\in\{f,1-f\}\\v\in\{g,1-g\}}}\norm{\{u,v\}_\alpha}_\infty,
\end{equation}
where the infimum is taken over all admissible pairs. \end{definition}

The admissible classes in Definition~\ref{def:contact-pb-subsets} are non-empty under the stated intersection conditions, so both $pb_{3, \alpha}$ and $pb_{4, \alpha}$ are finite and non-negative in these cases. Note that the admissible triples take values in the simplex $\{(a_0,a_1,a_2) \in [0,1]^3\,| \, a_i\ge0,\ \sum_i a_i=1\}$, and the admissible pairs take values in the square $[0,1]^2$. \rev{Compared with the symplectic counterparts (see \cite[Proposition 1.3]{BEP}), Definition \ref{def:contact-pb-subsets} takes the maximum over additional brackets of complementary Hamiltonian functions. These additional brackets ensure the symmetries of the resulting Poisson bracket invariants (see Lemma \ref{prop:subset-basic} below).} \rev{Indeed, the following example shows that these symmetries can fail if \keepmath{$pb_{3,\alpha}$} and \keepmath{$pb_{4, \alpha}$} are defined by directly copying the symplectic definitions of \keepmath{$pb_3$} and \keepmath{$pb_4$}.} 

\begin{example} Let $(Y, \ker \alpha)$ be a contact manifold. Suppose we define 
\begin{equation} \label{symp-cont-pb-3} 
\widetilde{pb}_{3, \alpha}(X_0, X_1, X_2): = \inf_{\vec h}\norm{\{h_0,h_1\}_\alpha}_\infty
\end{equation}
where $\mrev{h_i\ge0}$, $h_0 + h_1 + h_2 = 1$ and $h_i = 0$ near $X_i$. Consider $(X_0, X_1, X_2) = (\{p\}, \{\varphi_1^{\rrev{\epsilon}}(p)\}, Y)$, where the Reeb flow at $t= \rrev{\epsilon}>0$ \rev{displaces the point \keepmath{$\{p\}$}}, that is, $\varphi_1^{\rrev{\epsilon}}(p) \neq p$. \rev{Since \keepmath{$X_2 = Y$}, we have \keepmath{$h_2 \equiv 0$} and hence \keepmath{$h_1 = 1-h_0$}.} \rev{By (\ref{eq:constant-bracket}), every admissible triple of contact Hamiltonian functions \keepmath{$(h_0, h_1, h_2)$} satisfies} $\{h_0, h_1\}_{\alpha} = \{h_0, 1-h_0\}_{\alpha} = - R_{\alpha} h_0.$
\rev{Moreover, the endpoint conditions} $h_0(p) =0$ and $h_0(\varphi_1^{\rrev{\epsilon}}(p)) =1$ \rev{imply} 
\[ 1 = h_0(\varphi_1^{\rrev{\epsilon}}(p)) - h_0(p) = \int_0^{\rrev{\epsilon}} R_{\alpha} h_0 (\varphi_1^{t}(p)) dt \leq \rrev{\epsilon} \rrev{\norm{R_{\alpha} h_0}_\infty} \]
which implies that $\widetilde{pb}_{3, \alpha}(X_0, X_1, X_2) \geq \frac{1}{\rrev{\epsilon}}>0$. On the other hand, permute $(X_0, X_1, X_2)$ to $(X_2, X_0, X_1)$. Then $h_0 \equiv 0$ and $h_2 = 1- h_1$. \rev{\rrev{Therefore}, a function} $h_1: Y \to [0,1]$, vanishing near $p$ and equal to $1$ near $\varphi_1^{\rrev{\epsilon}}(p)$ \rev{determines an admissible triple of contact Hamiltonian functions} for $\widetilde{pb}_{3, \alpha}$, \rev{and the bracket in the infimum vanishes. Hence} $\widetilde{pb}_{3, \alpha}(X_2, X_0, X_1) = 0$. In particular, $\widetilde{pb}_{3, \alpha}(X_0, X_1, X_2) \neq \widetilde{pb}_{3, \alpha}(X_2, X_0, X_1)$. 

Similarly, suppose we define 
\begin{equation} \label{symp-cont-pb-4}  
\widetilde{pb}_{4, \alpha}(X_0, X_1, Y_0, Y_1): = \inf_{(f,g)}\norm{\{f,g\}_\alpha}_\infty
\end{equation}
where $\mrev{f,g:Y\to[0,1]}$, $f=i$ near $X_i$, and $g=i$ near $Y_i$ for $i=0,1$. Consider $(X_0, X_1, Y_0, Y_1) = (\{z\}, \{\varphi_1^{t}(p)\}_{t \in [0, \rrev{\epsilon}]}, \mrev{\{p\}}, \mrev{\{\varphi_1^{\rrev{\epsilon}}(p)\}})$, \rev{where again \keepmath{$\varphi_1^{\rrev{\epsilon}}(p)$} differs from \keepmath{$p$} and \keepmath{$z$} lies outside the entire trajectory \keepmath{$\{\varphi_1^{t}(p)\}_{t \in [0, \rrev{\epsilon}]}$}.} \rrev{For an admissible pair $(f,g)$, we have $f=1$ and $df=0$ along the Reeb segment, so $\{f,g\}_\alpha=R_\alpha g$ there. Integrating along the segment, as above, gives} $\rrev{\norm{\{f, g\}_{\alpha}}_\infty} \geq \frac{1}{\rrev{\epsilon}}$, which implies that $\widetilde{pb}_{4, \alpha}(X_0, X_1, Y_0, Y_1) \geq \frac{1}{\rrev{\epsilon}}>0.$ \rev{On the other hand, after the permutation \keepmath{$(X_1, X_0, Y_0, Y_1)$}, we have} $\widetilde{pb}_{4, \alpha}(X_1, X_0, Y_0, Y_1) =0$. \rrev{Indeed, choose $f$ supported in a small neighborhood of $z$, equal to one near $z$, and choose $g$ supported in a disjoint small neighborhood of $\varphi_1^{\rrev{\epsilon}}(p)$, equal to one there. These neighborhoods can be chosen to avoid the Reeb segment and $p$, respectively. The pair is admissible for the permuted quadruple, and its bracket vanishes because its supports are disjoint.} 
\end{example}

\rev{Unlike (\ref{symp-cont-pb-3}) and (\ref{symp-cont-pb-4}), Definition \ref{def:contact-pb-subsets} incorporates these symmetries directly.} \rev{Besides preserving these symmetries, the complementary coordinates also provide Reeb-derivative control:}
\[
\begin{aligned}
\max_{0\le i\le2}\norm{R_\alpha h_i}_\infty&\le2\max_{0\le i<j\le2}\norm{\{h_i,h_j\}_\alpha}_\infty,\\
\max\{\norm{R_\alpha f}_\infty,\norm{R_\alpha g}_\infty\}&\le2\max_{\substack{u\in\{f,1-f\}\\v\in\{g,1-g\}}}\norm{\{u,v\}_\alpha}_\infty.
\end{aligned}
\]
Therefore, Theorem~\ref{prop-3} \rev{yields lower semicontinuity of the two bracket maxima} in Definition~\ref{def:contact-pb-subsets} \rev{under uniform convergence of the input contact Hamiltonian functions}, without an additional derivative hypothesis (see Proposition~\ref{prop:subset-lsc}). 

\begin{example} In this example, we relate contact $pb_{3, \alpha}, pb_{4, \alpha}$ to symplectic $pb_3, pb_4$ in the setting of prequantizations. Let $\pi:(Y,\alpha)\to(M,\omega)$ be a prequantization over a closed symplectic manifold, with $d\alpha=\pi^*\omega$ and Reeb flow equal to the period-one circle action on the fibers. Then
\[
\begin{aligned}
pb_{3,\alpha}(\pi^{-1}A_0,\pi^{-1}A_1,\pi^{-1}A_2)&\le pb_{3}(A_0,A_1,A_2),\\
pb_{4,\alpha}(\pi^{-1}A_0,\pi^{-1}A_1,\pi^{-1}B_0,\pi^{-1}B_1)&\le pb_{4}(A_0,A_1,B_0,B_1)
\end{aligned}
\]
for compact input sets satisfying $A_0\cap A_1\cap A_2=\varnothing$ in the first inequality and $A_0\cap A_1=B_0\cap B_1=\varnothing$ in the second. 

The proof is based on the following equality between two Poisson brackets: for $F,G\in C^\infty(M)$, we have $\{\pi^*F,\pi^*G\}_\alpha = \pi^*\{F,G\}_\omega$, which implies that 
\[ \norm{\{\pi^*F,\pi^*G\}_\alpha}_\infty=\norm{\{F,G\}_\omega}_\infty. \]
For $pb_{3, \alpha}$,  \cite[Proposition~1.3]{BEP} allows us to use normalized Hamiltonian functions $F,G\ge0$ with $F+G\le1$, where $F=0$ near $A_0$, $G=0$ near $A_1$, and $F+G=1$ near $A_2$. Set
\[ h_0:=\pi^*F,\qquad h_1:=\pi^*G,\qquad h_2:=1-h_0-h_1. \]
Then $(h_0,h_1,h_2)$ is an admissible (contact) partition for the lifted $\pi^{-1}A_0,\pi^{-1}A_1,\pi^{-1}A_2$. Since the functions $h_0, h_1, h_2$ are Reeb-invariant, the brackets with the constant function vanish. So, $\{h_0,h_2\}_\alpha=-\{h_0,h_1\}_\alpha$ and $\{h_1,h_2\}_\alpha=\{h_0,h_1\}_\alpha$, which implies that 
\[ \max_{0\le i<j\le2}\norm{\{h_i,h_j\}_\alpha}_\infty=\norm{\{F,G\}_\omega}_\infty. \]
Taking the infimum over normalized (symplectic) Hamiltonian functions proves the inequality between $pb_{3, \alpha}$ and $pb_3$. A similar argument works for $pb_{4, \alpha}$ and $pb_4$. 
\end{example}

\rev{There are again two approaches to obtaining lower bounds for \keepmath{$pb_{3, \alpha}$} and \keepmath{$pb_{4, \alpha}$}, depending on the hypotheses. To formulate the quasi-morphism bounds, recall the notions in \cite[Definitions 1.7 and 1.9]{BZ}.} For a nonzero monotone homogeneous quasi-morphism $\mu:\widetilde{\Cont}_0(Y,\xi)\to\R$, a compact set $S\subset Y$ is $\mu$-subheavy if $\mu(\widetilde\varphi_h)=0$ for every autonomous $h\in C^\infty(Y)$ with $h|_S=0$; it is $\mu$-superheavy if $\mu(\widetilde\varphi_h)>0$ for every autonomous $h\in C^\infty(Y)$ with $h|_S>0$. These properties do not depend on the choice of contact form, and every $\mu$-superheavy set is $\mu$-subheavy by \cite[Proposition 1.10 and Theorem 1.11]{BZ}. 

\begin{mainthm}\label{thm:subset-bounds} Suppose that $(Y,\xi=\ker\alpha)$ admits a nonzero monotone homogeneous quasi-morphism $\mu$, normalized by \eqref{eq:mu-normalization}. For compact input subsets, the following statements hold.
\begin{enumerate}[(i)]
\item If $X_0\cap X_1\cap X_2=\varnothing$ and each $X_i$ is $\mu$-subheavy, then
\[ pb_{3,\alpha}(X_0,X_1,X_2)\ge\frac{1}{256D(\mu)}. \]
\item If $X_0\cap X_1=Y_0\cap Y_1=\varnothing$ and each $X_i\cup Y_j$, $i,j\in\{0,1\}$, is $\mu$-subheavy, then
\[ pb_{4,\alpha}(X_0,X_1,Y_0,Y_1)\ge\frac{1}{8064D(\mu)}. \]
If, in addition, $Y_1$ is $\mu$-subheavy, then $pb_{4,\alpha}(X_0,X_1,Y_0,Y_1)\ge\frac{1}{768D(\mu)}.$
\end{enumerate}
\end{mainthm}

\begin{remark} \rev{Theorem \ref{thm:subset-bounds} does not assume that \keepmath{$\mu$} has the vanishing property.} \end{remark}

Theorem~\ref{thm:subset-bounds} is a contact counterpart of \cite[Theorem 1.15]{BEP}, which gives symplectic $pb_3$ and $pb_4$ lower bounds from superheavy sets and a Poisson bracket inequality for a symplectic quasi-state; see also \cite[Section 7.3]{Polterovich-Rosen}. In our proof, subheaviness makes the quasi-morphism vanish on the isotopy class generated by each test function in a triple and by each product $fg$, $f(1-g)$, $(1-f)g$, and $(1-f)(1-g)$ in a quadruple. Hence, the partition estimate used for Theorem~\ref{prop-4} applies directly. \rev{To obtain the strengthened bound in the last conclusion, we also use the relation}
\begin{equation} \label{eq:subset-four-to-three}
pb_{3,\alpha}(X_0\cup Y_0,X_1\cup Y_0,Y_1)\le3pb_{4,\alpha}(X_0,X_1,Y_0,Y_1),
\end{equation}
proved in Lemma~\ref{prop:subset-four-to-three}. \rev{The hypothesis for that strengthened bound can equivalently be stated by requiring only} $X_0\cup Y_0$, $X_1\cup Y_0$, and $Y_1$ to be $\mu$-subheavy, since every compact set containing a $\mu$-subheavy set is $\mu$-subheavy. This is the contact version of the three-set condition in \cite[Theorem 1.15(ii)]{BEP}. \rev{Alternatively, the spectral invariant gives a separate conclusion under displacement hypotheses, without requiring a quasi-morphism.}

\begin{cor}\label{cor:subset-energy-bounds} Suppose that $(Y,\xi=\ker\alpha)$ satisfies the filling hypothesis of Theorem~\ref{thm:spectral-cover}. Then for compact input subsets, the following statements hold.
\begin{enumerate}[(i)]
\item If $X_0\cap X_1\cap X_2=\varnothing$ and each $Y\setminus X_i$ is contact displaceable, then
\begin{equation} \label{eq:subset-triple-bound}
pb_{3,\alpha}(X_0,X_1,X_2)\ge\frac{1}{128\displaystyle\sum_{i=0}^2e_\alpha(Y\setminus X_i)}.
\end{equation}
\item If $X_0\cap X_1=Y_0\cap Y_1=\varnothing$ and each $Y\setminus(X_i\cup Y_j)$ is contact displaceable, then
\begin{equation} \label{eq:subset-square-bound}
pb_{4,\alpha}(X_0,X_1,Y_0,Y_1)\ge\frac{1}{2688\displaystyle\sum_{i,j=0}^1e_\alpha(Y\setminus(X_i\cup Y_j))}.
\end{equation}
\end{enumerate}
\end{cor}

The proofs of Theorem \ref{thm:subset-bounds} and Corollary \ref{cor:subset-energy-bounds} are given in Section~\ref{sec:proof-H}. \rev{The constants in these lower bounds reflect the complexity of the contact Poisson bracket formula (\ref{eq:subset-chain-rule}) (see also (\ref{eq:contact-chain-rule})).} 

\subsection{Dynamical applications}\label{subsec:subset-dynamics}
\rev{The first application concerns the existence of contact Hamiltonian chords, using \keepmath{$pb_{4, \alpha}$} as defined and estimated in Section \ref{subsec:pb34}.} In this application, the four inputs have different roles: $X_0,X_1$ are the endpoint sets that a trajectory must connect, and $Y_0,Y_1$ prescribe a separation between the values of the Hamiltonian. Thus a lower bound for $pb_{4,\alpha}(X_0,X_1,Y_0,Y_1)$ gives a time bound for Hamiltonian flows. \rev{To apply our earlier estimates,} we consider smooth autonomous Hamiltonians $h$ satisfying $R_\alpha h=0$. 

\begin{mainthm}\label{thm:subset-chords} Let $X_0,X_1,Y_0,Y_1\subset Y$ be non-empty compact sets with $X_0\cap X_1=Y_0\cap Y_1=\varnothing$, and suppose that
\[
p:=pb_{4,\alpha}(X_0,X_1,Y_0,Y_1)>0.
\]
Let $h\in C^\infty(Y)$ satisfy $R_\alpha h=0$ and assume that $\min_{Y_1}h-\max_{Y_0}h>0$. If $T,L>0$ satisfy
\begin{equation}\label{eq:main-condition}
\max\left\{\frac{\min_{Y_1}h-\max_{Y_0}h}{L},\ \frac1T+\frac{\max\{|\max_{Y_0}h|,|\min_{Y_1}h|\}}{L}\right\}\le p\left(\min_{Y_1}h-\max_{Y_0}h\right),
\end{equation}
then there exist $x\in X_0$ and (signed) times $t,s\in\R$ such that
\begin{equation}\label{eq:assisted-chord}
\varphi_1^s\varphi_h^t(x)\in X_1,
\qquad |t|\le T,\qquad |s|\le L,\qquad |t|+|s|>0.
\end{equation}
If, furthermore, $\left(\bigcup_{|s|\le L}\varphi_1^s(X_0)\right)\cap X_1=\varnothing$, then $|t|$ in \eqref{eq:assisted-chord} is nonzero.
\end{mainthm}

Since $\max\{|\max_{Y_0}h|,|\min_{Y_1}h|\} \leq \|h\|_{\infty}$, one can impose the following condition, which is simpler than (\ref{eq:main-condition}): 
\begin{equation}\label{eq:simple-condition}
\frac1T+\frac{\norm{h}_\infty+\min_{Y_1}h-\max_{Y_0}h}{L}<p\left(\min_{Y_1}h-\max_{Y_0}h\right).
\end{equation}
However, the condition \eqref{eq:main-condition} uses only values of $h$ on the input subsets $Y_0, Y_1$. Also, if $h=0$ near $Y_0$ and $h=1$ near $Y_1$, then the condition \eqref{eq:main-condition} reduces to $\frac{1}{T} + \frac{1}{L} \leq p$. Taking $T=L=2/p$ thus bounds both time parameters by $2/p$. The following corollary gives an ordinary Hamiltonian chord, with no Reeb translation, under an additional endpoint hypothesis. 

\begin{cor}\label{cor:one-endpoint}
In the setting of Theorem~\ref{thm:subset-chords}, assume that $p>0$, $R_\alpha h=0$, and $\min_{Y_1}h-\max_{Y_0}h>0$. If at least one of $X_0,X_1$ is invariant under the entire Reeb flow, then there exist $x\in X_0$ and a (signed) time $t$ such that
\begin{equation}\label{eq:ordinary-chord}
\varphi_h^t(x)\in X_1 \quad \mbox{where} \,\, 0<|t|\le\frac{1}{p(\min_{Y_1}h-\max_{Y_0}h)}.
\end{equation}
\end{cor}

\begin{proof} Fix $T>1/(p(\min_{Y_1}h-\max_{Y_0}h))$. Since $1/T<p(\min_{Y_1}h-\max_{Y_0}h)$, we can choose $L$ sufficiently large so that \eqref{eq:main-condition} holds. Then Theorem~\ref{thm:subset-chords} gives $x\in X_0$ and $t,s\in\R$ with $|t|\le T$ and $\varphi_1^s\varphi_h^t(x)\in X_1$. If $X_1$ is Reeb-invariant, then applying the inverse Reeb translation shows directly that $\varphi_h^t(x)\in X_1$. If instead $X_0$ is Reeb-invariant, put $x':=\varphi_1^s(x)\in X_0$. Then commutation of the two flows (since $h$ is Reeb-invariant) gives
\[
\varphi_h^t(x')=\varphi_1^s\varphi_h^t(x)\in X_1.
\]
Therefore, in either case there is an ordinary $h$-chord from $X_0$ to $X_1$ with $|t|\le T$. Letting $T$ decrease to $1/(p(\min_{Y_1}h-\max_{Y_0}h))$ and passing to a convergent subsequence of initial points and times gives the stated bound. The limit time is nonzero since $X_0\cap X_1=\varnothing$. \end{proof}

The proof of Theorem~\ref{thm:subset-chords}, given in Section~\ref{sec:proof-I}, adapts the averaging argument of \cite[Lemma 4.1 and Theorem 1.10]{BEP} to two commuting flows; see also the travel-time interpretation of $pb_4$ in \cite[Section 7.4]{Polterovich-Rosen}. Since $h$ is strict, its Hamiltonian flow commutes with the Reeb flow, and $X_hf=-\{f,h\}_\alpha$. Assuming that no connection with the prescribed Hamiltonian and Reeb time bounds exists, we first separate the two sets swept out from the endpoints by these flows. Averaging this separating function over both time variables preserves its prescribed values near the endpoint sets and gives simultaneous bounds on $X_hf$ and $R_\alpha f$. A smooth cutoff of $h$ then converts these derivative bounds into estimates for all four complementary brackets in the contact invariant. Therefore, positivity of the invariant forces a Hamiltonian trajectory followed by a quantitatively controlled Reeb translation. 

\begin{remark} We compare Theorem \ref{thm:subset-chords} and Corollary \ref{cor:one-endpoint} with \cite[Sections 1.4--1.5 and 4.5]{EP-Legendrian}. Entov and Polterovich use persistence modules from Legendrian contact homology to bound symplectic pb-invariants in the symplectization. Their interlinking results give directed chords for suitable Legendrian pairs and positive contact Hamiltonians without requiring the Hamiltonians to be Reeb-invariant. Our results use the contact pb-invariant of four compact sets and apply to strict autonomous Hamiltonians. These Hamiltonians may change sign, and the endpoint sets need not be Legendrian. Also, the conclusion allows either direction of travel. \end{remark}

The following examples illustrate the conclusions of Theorem \ref{thm:subset-chords} and Corollary \ref{cor:one-endpoint} in several concrete geometric settings.

\begin{example} \label{ex-s1-copshere}
Let $S^1=\R/\mathbb Z$ and equip $D^*S^1=(\R/\mathbb Z)\times[-1,1]$ with the Liouville form $\lambda=p\,dq$. Its contact boundary $S^*S^1=\partial D^*S^1$ has two connected components:
\[
Y=S^*S^1=(\R/\mathbb Z)\times\{+1,-1\}.
\]
For our input contact Hamiltonian $h$, the condition $R_\alpha h=0$ implies that it is constant on each component. Fix $0<d\le1/2$ and put
\[
X_0 =\{(0,+1)\}, \,\,\,\, X_1=\{(d,+1)\}, \,\,\,\, Y_0=(\R/\mathbb Z)\times\{-1\}, \,\,\,\, Y_1=(\R/\mathbb Z)\times\{+1\}.
\]
It is readily verified that \begin{equation}\label{eq:circle-pb}
pb_{4,\alpha}(X_0,X_1,Y_0,Y_1)=\frac1d.
\end{equation}
For a Reeb-invariant Hamiltonian $h: Y  \to \R$ such that $h|_{Y_0}=0$ and $h|_{Y_1}=1$, Theorem~\ref{thm:subset-chords} applies whenever the inequality $\frac1T+\frac1L\le\frac1d$ holds. For instance, take $T= L = 2d$. The conclusion of Theorem~\ref{thm:subset-chords} is therefore equivalent to the arithmetic relation $t+s\equiv d\pmod1$ for some $|t|\le T$ and $|s|\le L$, under the stated time bound. 
\end{example}

\begin{example}\label{ex-2-sphere} In this example, we will demonstrate how the contact Poisson bracket invariant \(pb_{4,\alpha}\) detects fast diffusion in the precession of a great circle on the round sphere \(S^2\). Here, ``fast diffusion'' (see \cite[Section~1.3]{EP-nonlinearity}) means that as its plane turns at angular speed \(\ep\), the first component of its unit normal can pass between two fixed levels in time \(O(\ep^{-1})\). We call such a passage a crossing. Such a crossing persists, on the same time scale \(O(\ep^{-1})\), if the perturbing function is sufficiently small relative to \(\ep\) on just two arcs of the sphere of normals (see the assumption in Corollary \ref{cor-fast-diffusion} below).

\medskip

\noindent{\bf (i) Precession of the geodesic flow on $S^2$}. Let $q_0$ be a point on the unit round sphere $S^2\subset\R^3$ and let $v_0\in T_{q_0}S^2$ be a unit vector. The geodesic with these initial data is $\gamma(t)=\cos t\,q_0+\sin t\,v_0$. Its image is the great circle in the plane perpendicular to $n(0):=q_0\times v_0$. In particular, the point moves, but the plane does not. Now fix $0<c<1/\sqrt{2}$, $q_0=(0,1,0)$, and $v_0=(-\sqrt{1-c^2},0,-c)$, so that the initial normal direction is $n(0)=(-c,0,\sqrt{1-c^2}).$ For $\ep>0$, consider the system 
\begin{equation} \label{eq:geo-precession}
\dot n(t)=(0, \ep, 0) \times n(t) 
\end{equation}
whose solution is
\begin{equation} \label{solution-nt}
n(t)=\bigl(\sin(\ep t-\arcsin c),\,0,\,\cos(\ep t-\arcsin c)\bigr).
\end{equation}
 Equivalently, if $A_\theta: \R^3 \to \R^3$ denotes rotation through angle $\theta$ about the second coordinate axis, in the direction $(0,1,0)$ of $\R^3$, then the solution above is $n(t)= A_{\ep t}n(0)$. In other words, the plane of the great circle rotates with its normal $n(t)$, a motion called precession. Moreover, under this precession, the perturbed flow with state $(q_{\ep}(t), v_{\ep}(t))$ on $S^*S^2$ satisfies the following system of equations: 
 \begin{equation} \label{geo-flow-precession}
 \dot{q}_{\ep}(t) = v_{\ep} + (0, \ep, 0) \times q_{\ep},\qquad  \dot{v}_{\ep}(t) = -q_{\ep} + (0, \ep, 0) \times v_{\ep}. 
 \end{equation}
 Figure~\ref{fig:round-geodesic-precession} shows the two motions at $t_\ast=5\pi/2$ for $c=\ep=1/10$, when the precessing plane has turned through angle $\ep t_\ast=\pi/4$. The supplementary movie shows one full period of precession~\cite{Zhang-precession-video}.

\begin{figure}[htbp]
\centering
{\color{black}
\includegraphics[width=.65\linewidth]{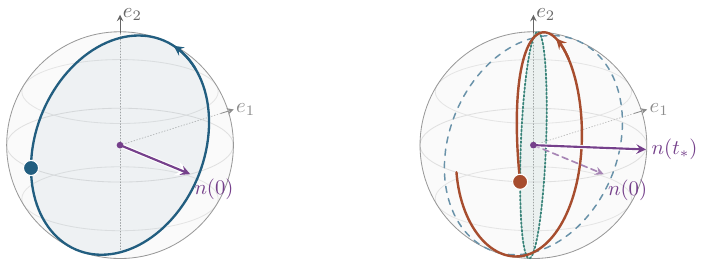}
\caption{Left: The round geodesic, with fixed normal $n(0)$. Right: The precessing geodesic plane, whose normal turns from $n(0)$ to $n(t_\ast)$.}\label{fig:round-geodesic-precession}
}
\end{figure}
\FloatBarrier

\noindent{\bf (ii) Fast diffusion within a finite time interval}. From (\ref{solution-nt}), consider the time interval $[0, \frac{2\arcsin(c)}{\ep}]$. Since $|\ep t-\arcsin c|\le\arcsin c$ in this interval, the velocity $\dot n_1(t)=\ep n_3(t)$ lies between $\ep\sqrt{1-c^2}$ and $\ep$. Then, integrating from $0$ to $t$ gives
\[
 \ep t\sqrt{1-c^2} + n_1(0) \ \le\ n_1(t) =\ep\int_0^t n_3(s)\,ds + n_1(0) \ \le\ \ep t  + n_1(0)
\]
for $0\le t\le \frac{2\arcsin(c)}{\ep}.$ In particular, $n_1(t)$ increases from $n_1(0) = -c$ to $n_1(2\arcsin(c)/\ep) = c$ in time $O(\ep^{-1})$. The left picture of Figure~\ref{fig:round-geodesic-n1} compares the three functions $\ep t\sqrt{1-c^2} + n_1(0)$, $n_1(t)$, and $\ep t  + n_1(0)$, for parameters $c=1/10$ and $\ep=1/100$. The graph is nearly linear on the indicated interval. The right picture of 
Figure~\ref{fig:round-geodesic-n1} shows a crossing from level $n_1 = -c$ to $n_1 = c$ in the {\bf sphere of unit normals} (i.e., $\{n \in \R^3 \,| \, |n| =1\}$), where a point $n$ can also represent the oriented great circle perpendicular to $n$.
\begin{figure}[htbp]
\centering
{\color{black}
\includegraphics[width=.78\linewidth]{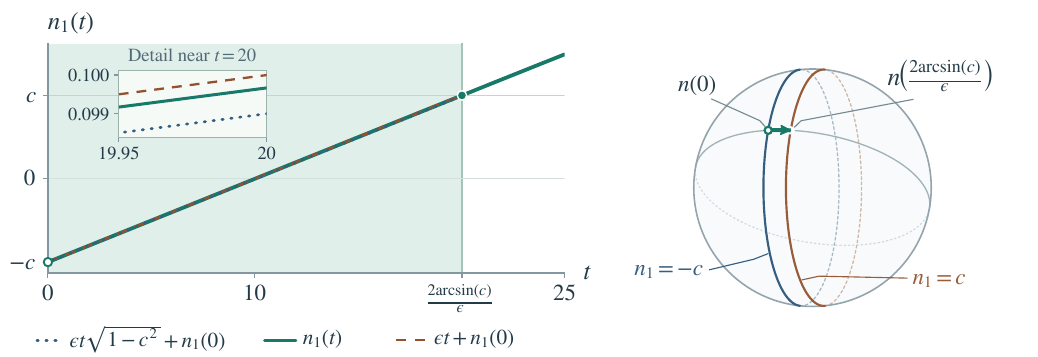}
\caption{Left: The graph of $n_1(t)$ on the time interval $[0,25]$, with the interval $[0,\frac{2\arcsin(c)}{\ep}]$ shaded. Right: A crossing between $n_1= -c$ and $n_1=c$.}\label{fig:round-geodesic-n1}
}
\end{figure}
\FloatBarrier

\noindent{\bf (iii) Robustness under $C^0$-perturbations}. In the upper hemisphere of the sphere of unit normals, consider the closed quadrilateral $Q:=\{n_3>0,\ |n_1|\le c,\ |n_2|\le c\}$, which can be identified with the square $[-c,c]^2$. Denote its left and right sides by $A_0,A_1$ and its bottom and top sides by $B_0,B_1$:
\begin{equation}\label{eq:S2-four sets}
 A_i:=Q\cap\{n_1=(2i-1)c\},\qquad B_i:=Q\cap\{n_2=(2i-1)c\},\qquad i=0,1.
\end{equation}
See Figure \ref{fig:round-geodesic-chord} (and ignore the orange chord temporarily). 
\begin{figure}[!htbp]
\centering
{\color{black}
\includegraphics[width=.78\linewidth]{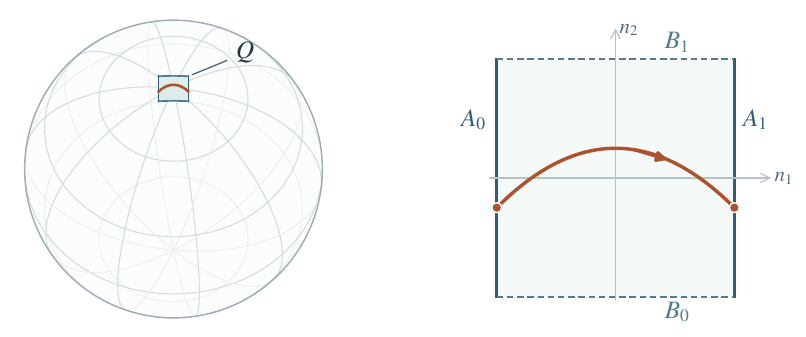}
\caption{Left: The region $Q$ on the sphere of unit normals. Right: A crossing from $A_0$ to $A_1$ for a specific perturbation $k(n)=\frac{\ep n_1^2}{2c}$.}\label{fig:round-geodesic-chord}
}
\end{figure}
To perturb \eqref{eq:geo-precession}, take a smooth function $k: S^2 \to \R$ on the sphere of unit normals and consider 
\begin{equation} \label{eq-presession-perturb}
\dot n=\left((0,\ep,0) +\nabla_{S^2}k(n)\right)\times n
\end{equation}
where the perturbation term $\nabla_{S^2}k(n)$ means the gradient of \(k\) at \(n\), taken with respect to its round metric. 

\begin{cor}\label{cor-fast-diffusion}
For any smooth $k:S^2\to\R$ satisfying $\norm{k}_{C^0(B_0\cup B_1)}\le \frac{c \ep}{2}$, the perturbed precession in (\ref{eq-presession-perturb}) admits a trajectory segment with one endpoint on each of $A_0$ and $A_1$ and duration at most $7\operatorname{area}(Q)/(c\ep)$.
\end{cor}

Figure \ref{fig:round-geodesic-chord} illustrates such a crossing, with the sign of time specified, and with a specific choice of function $k$ on $S^2$. Corollary~\ref{cor-fast-diffusion} requires only $C^0$-control on $B_0\cup B_1$ and places no bound on the derivatives of $k$. For a similar conclusion but in a different setting, see  \cite[Corollary~1.8]{EP-nonlinearity}. 

\medskip

\noindent{\bf (iv) Contact geometric interpretation}. The precession in (\ref{eq:geo-precession}) can be reformulated in terms of contact geometry. Identify the round unit cosphere bundle with $Y:=S^*S^2=\{(q,v) \in \R^3 \times \R^3 \,| \, |q|=|v|=1,\ q\cdot v=0\}$ and take $\alpha:=v\cdot dq$. The Reeb flow on $(Y,\alpha)$, generated by the constant contact Hamiltonian $1$, corresponds to the round geodesic flow on the base $S^2$. The map $\pi:Y\to S^2$, $(q,v)\mapsto q\times v$, identifies the quotient by the Reeb flow with the sphere of unit normals. Consider the contact Hamiltonian $h_{\ep}: Y \to \R$ defined by 
\begin{equation} \label{cont-ham-precession}
 h_\ep:=1+\ep(n_2\circ\pi) 
\end{equation}
where $n_2:S^2\to\R$ is the second coordinate function. The flow of $h_{\ep}$ satisfies (\ref{geo-flow-precession}), and its normals satisfy (\ref{eq:geo-precession}). In this way, we reinterpret the precession of the geodesic flow in terms of contact Hamiltonian dynamics. 

Moreover, the four subsets $A_0, A_1, B_0, B_1 \subset S^2$ in (\ref{eq:S2-four sets}) lift to four compact annuli in $Y$: 
\begin{equation}\label{eq:S2-four sets-2}
 X_i:=\pi^{-1}(A_i),\qquad Y_j:=\pi^{-1}(B_j),\qquad i,j\in\{0,1\},
\end{equation}
which are Reeb invariant, and opposite annuli are disjoint. The perturbation $k$ in Corollary~\ref{cor-fast-diffusion} adds the term $k\circ\pi$ to $h_\ep$, and the resulting contact Hamiltonian $h:=h_\ep+k\circ\pi$ satisfies $R_\alpha h=0$ because it is constant on the Reeb fibers. It remains to verify the Hamiltonian gap and the positive lower bound on $pb_{4,\alpha}$ required by Corollary~\ref{cor:one-endpoint}. 

\smallskip

\noindent{\bf Proof of Corollary \ref{cor-fast-diffusion}} Write $u(n):=\ep n_2+k(n)$, so $h=1+u\circ\pi$. The hypothesis gives $u\le-c\ep/2$ on $B_0$ and $u\ge c\ep/2$ on $B_1$. Therefore, 
\[ \min_{Y_1}h-\max_{Y_0}h\ge c\ep.\]
With the notation of (\ref{eq:S2-four sets-2}), it suffices to check that $pb_{4,\alpha}(X_0,X_1,Y_0,Y_1)$ admits a positive lower bound. Put $D:=\pi^{-1}(Q)$ and let $P$ be the maximum in \eqref{eq:subset-square-definition} for an admissible pair $(f,g)$ satisfying 
\[ f|_{X_0} = 0, \quad f|_{X_1}=1, \quad g|_{Y_0} = 0, \quad g|_{Y_1} = 1. \]
Consider the 2-form $f dg \wedge \alpha$. With the outward round area form $\omega$ on the sphere of normals, $d\alpha=-\pi^*\omega$. Only $X_1$ contributes to the boundary integral. On this annulus $d\alpha=0$, so $dg\wedge\alpha=d(g\alpha)$. Its two boundary fibers have $g=0,1$ and Reeb period $2\pi$, giving 
\[\left| \int_{\partial D} f dg \wedge \alpha \right| = 2\pi. \]
On the other hand, by a direct computation, we have 
\[ d(f dg \wedge \alpha) = \bigl(\{f,g\}_\alpha+gR_\alpha f-2fR_\alpha g\bigr)\,\alpha\wedge d\alpha. \]
Therefore, by Stokes' theorem,
\[ 2\pi=\left|\int_D\bigl(\{f,g\}_\alpha+gR_\alpha f-2fR_\alpha g\bigr)\,\alpha\wedge d\alpha\right|\le 7P\cdot 2\pi\,\operatorname{area}(Q),
\]
where the coefficient in the inequality is bounded by $P+2P+4P$ by \eqref{eq:subset-square-reeb}. Taking the infimum yields $pb_{4,\alpha}(X_0,X_1,Y_0,Y_1)\ge 1/(7\operatorname{area}(Q))$. Then Corollary~\ref{cor:one-endpoint} yields a crossing of duration at most $7\operatorname{area}(Q)/(c\ep)$. \qed

\begin{remark} The robustness of the fast diffusion does {\rm not} violate Nekhoroshev's exponential stability in \cite[1.4 Theorem]{Neh1977}. To make a comparison, we lift the contact geometric interpretation in part (iv) into a (symplectic) Hamiltonian dynamical system. Consider an identification 
\[ ((0,\infty)\times Y,d(r\alpha)) \simeq T^*S^2 \backslash o_{T^*S^2} \quad \mbox{by $(r,q,v)\mapsto(q,\eta=rv)$} \]
where $\eta \in T^*_q S^2$ denotes the covector corresponding to $rv$ under the round metric. The unperturbed Hamiltonian on $T^*S^2 \backslash o_{T^*S^2}$ that generates the geodesic flow on $S^2$ is 
\[ H_0(q, \eta)  = |\eta|. \]
This corresponds to $h_0$ in (\ref{cont-ham-precession}) when $\ep=0$. On a regular action chart one can take $I_1=|\eta|$ and $I_2=(q\times\eta^\sharp)_3$ (where $\eta^{\sharp}$ is the dual of $\eta$): their flows are the $2\pi$-periodic round geodesic flow and axial rotation, respectively. Thus $H_0=I_1$ is constant on every line $I_1=\mathrm{constant}$ in this chart, so it fails the steepness criterion in \cite[Theorem 1.5]{Nie06}. Therefore, Nekhoroshev's result does not apply. \end{remark}  \end{example}

\rev{The second application gives a geometric interpretation of a \keepmath{$C^0$}-bound on the contact Poisson bracket in terms of deformations of the contact form.} \rev{The symplectic model is the deformation considered in \cite[Section 1.6]{BEP}, \keepmath{$\omega_s=\omega-s\,df\wedge dg$}, for which the identity} $\omega_s^n=(1-s\{f,g\}_\omega)\omega^n$ \rev{gives an interval on which the deformed forms remain symplectic. Combined with persistence of holomorphic curves, this calculation yields lower bounds for Poisson brackets. In the contact setting, however, the deformation} $\alpha_s=\alpha-sf\,dg$ \rev{has the coefficient} $X_fg=\{f,g\}_\alpha+gR_\alpha f$. \rev{To control this additional Reeb term, we again use the brackets of complementary Hamiltonian functions: a direct computation shows that these brackets, as in (\ref{eq:subset-square-definition}), control this coefficient and hence determine an explicit interval on which the deformed forms remain contact.}

\begin{mainthm}\label{thm:contact-form-deformation} Let $(Y, \xi = \ker \alpha)$ be a closed contact manifold and $f,g\in C^\infty(Y)$. \rev{Consider the deformation \keepmath{$\alpha_s=\alpha-sf\,dg$},} $s\in\R$. Then
\begin{equation}\label{eq:contact-deformation-volume}
\alpha_s\wedge(d\alpha_s)^n=(1-sX_fg)\,\alpha\wedge(d\alpha)^n.
\end{equation}
In particular, for $|s|\norm{X_fg}_\infty<1$, the forms $\alpha_s$ remain contact forms, and the maps $\Psi_s(x)=\varphi_f^{-sg(x)}(x)$ form an isotopy of diffeomorphisms starting at the identity such that
\begin{equation}\label{eq:contact-deformation-pullback}
\Psi_s^*\alpha=e^{\kappa_s}\alpha_s,\qquad \norm{\kappa_s}_\infty\le |s|\norm{g}_\infty\norm{R_\alpha f}_\infty.
\end{equation}
\end{mainthm}

\begin{remark} If $R_\alpha f=0$, then (\ref{eq:contact-deformation-pullback}) in Theorem~\ref{thm:contact-form-deformation} gives $\Psi_s^*\alpha=\alpha_s$, so $\Psi_s$ conjugates the Reeb flow of $\alpha_s$ to the Reeb flow of $\alpha$ through this deformation. \end{remark}

\rev{The following observation shows how the contact Poisson bracket controls the range of \keepmath{$s$} in Theorem \ref{thm:contact-form-deformation} through the norm \keepmath{$\norm{X_fg}_\infty$}.} Since $\{f,1\}_\alpha=-R_\alpha f$, bilinearity gives
\[ X_fg=\{f,g\}_\alpha+gR_\alpha f=(1-g)\{f,g\}_\alpha-g\{f,1-g\}_\alpha. \]
\rev{\rrev{Therefore}, if} $0\le g\le1$, then
\[
\begin{aligned}
\norm{X_fg}_{\infty} \leq \max\bigl\{\norm{\{f,g\}_\alpha}_{\infty}, \norm{\{f,1-g\}_\alpha}_{\infty}\bigr\}.
\end{aligned}
\]
\rev{Therefore, in this case, whenever} $|s|\max\bigl\{\norm{\{f,g\}_\alpha}_{\infty}, \norm{\{f,1-g\}_\alpha}_{\infty}\bigr\}<1$, \rev{the deformed 1-form \keepmath{$\alpha_s$} remains contact.} \rev{For linear combinations} $f, g$ \rev{of a partition of unity, we obtain the following consequence.}

\begin{cor}\label{cor:cover-deformation} Let $\mathcal U=\{U_1,\ldots,U_N\}$ be a finite open cover of $Y$. For every $T>0$ with $T \cdot pb_\alpha(\mathcal U)<1$, there exists a partition of unity $\vec f=(f_1,\ldots,f_N)$ subordinate to $\mathcal U$ such that, for every $a,b\in[0,1]^N$ and $|s|\le T$, the deformed form
\postdisplaypenalty=10000
\begin{equation}\label{eq:cover-contact-deformation}
\alpha-s\left(\sum_{i=1}^Na_if_i\right)d\left(\sum_{i=1}^Nb_if_i\right)
\end{equation}
is contact and its kernel is isotopic to $\xi$.
\end{cor}

For an admissible quadruple, the same argument shows that $T \cdot pb_{4,\alpha}(X_0,X_1,Y_0,Y_1)<1$ guarantees an admissible pair $(f,g)$ for which $\alpha-sf\,dg$ remains contact for all $|s|\le T$. The proofs of Theorem \ref{thm:contact-form-deformation} and Corollary \ref{cor:cover-deformation} are given in Section~\ref{sec:proof-J}.

\medskip\noindent \textbf{Acknowledgements.}  \rev{This project began during the author's work on contact Poisson bracket rigidity for contact coisotropic submanifolds in \cite{Rosen-Zhang}. Substantial progress was made recently, partly inspired by \cite{DUZ} (see Remark \ref{rmk-conformal-factor} above). The paper was completed during a visit to the University of Notre Dame, and the author thanks Richard Hind for his encouragement and interest in this work. The author also thanks Lev Buhovsky and Leonid Polterovich for useful discussions, in particular, for suggestions on adding Example~\ref{ex-ellipsoid-cont-non-dis}, Example~\ref{ex-2-sphere} and Corollary~\ref{cor-fast-diffusion}; Dylan Cant for discussions concerning the proof of Lemma \ref{lem:spectral-barrier}; Adi Dickstein for discussions about the quantitative big fiber theorem, Theorem \ref{thm-qcbft}; Shira Tanny for Remark~\ref{rmk-change-contact-form} and the discussion immediately preceding Section~\ref{subsec:big-fiber}; and Michael Usher for discussions about a quantitative version of the characterization of contact coisotropic submanifolds by the contact Poisson bracket in \cite[Proposition 1.3]{Rosen-Zhang}, to be explored in future work.} The author is partially supported by the National Key R\&D Program of China (Grant No.~2023YFA1010500) and NSFC Grants No.~12301081, No.~12511540054, No.~12361141812, and No.~12631007.

\medskip

\noindent {\bf About AI use}: Large language models (LLMs) were used for proofreading and language polishing, as well as for the following aspects of this work: (i) calculating the contact Poisson bracket formulas in \eqref{eq:contact-chain-rule} and \eqref{eq:subset-chain-rule}; (ii) suggesting the perturbation (\ref{eq:function}) in Example \ref{ex-ellipsoid-cont-non-dis}, and determining an explicit value of the constant $c_1$ in the same example; (iii) constructing the counterexample in Example \ref{ex-counter}, detailed in Section \ref{sec:counterexample}. The author constructed a counterexample to an unconditional version of Theorem \ref{prop-3} by adapting Example 1.2.2 in \cite{Buhovsky-C0}; LLMs produced the example presented here, which disproves unconditional versions of both Theorem \ref{prop-2} and Theorem \ref{prop-3}; (iv) providing details for the proof of Lemma~\ref{lem:spectral-displacer}; (v) producing the figures in Example \ref{ex-2-sphere} and the illustrative movie \cite{Zhang-precession-video}. Finally, (vi) LLMs suggested collecting several useful identities and estimates in Lemma \ref{lem:contact-flow-identities} and simplified the proof of estimate \eqref{eq:transport-growth}, for which the author's original argument used Gronwall's inequality.

\section{Proof of Theorem \ref{prop-1}}\label{sec:proof-A}

For a based contact isotopy $\psi_s$, $0\le s\le t$, generated by $h_s$, write
\begin{equation}\label{eq:contact-path-length}
\ell_\alpha(\psi):=\int_0^t\norm{h_s}_\infty\,ds.
\end{equation}
Whenever such a path is regarded as an element of $\widetilde{\Cont}_0(Y,\xi)$, we reparametrize the path on $[0,1]$. We use this convention throughout the proofs.

For $\varphi\in\Cont_0(Y,\xi)$, the Shelukhin norm is
\begin{equation} \label{she-norm}
|\varphi|_\alpha =\inf_{\varphi_h^1=\varphi}\int_0^1\norm{h_t}_\infty\,dt.
\end{equation}
It is nondegenerate, symmetric, and satisfies the triangle inequality; see \cite{Shelukhin-cont}. However, it need not be conjugation invariant. \rev{Nevertheless, conjugation can be controlled by the following estimate:}
\begin{equation} \label{eq:conjugation-bound}
|\psi\varphi\psi^{-1}|_\alpha \le e^{\max_Y\kappa(\psi)}|\varphi|_\alpha.
\end{equation}
Indeed, conjugating an isotopy generated by $h_t$ gives the Hamiltonian $(e^{\kappa(\psi)}h_t)\circ\psi^{-1}$. Therefore, taking its uniform norm, integrating, and then taking the infimum proves \eqref{eq:conjugation-bound}.

\rev{To compare the contact flows that appear in the displacement, quasi-morphism, and spectral arguments, we first collect a few useful identities.}

\begin{lemma}\label{lem:contact-flow-identities} Let $h_s,k_s$ be smooth time-dependent contact Hamiltonians with based flows $\varphi_h^s,\varphi_k^s$, for $0\le s\le t$. Then
\begin{equation}\label{conformal}
\kappa(\varphi_h^s)(x)=\int_0^s (R_\alpha h_r)(\varphi_h^r(x))\,dr,\qquad \norm{\kappa(\varphi_h^s)}_\infty\le\int_0^s\norm{R_\alpha h_r}_\infty\,dr,
\end{equation}
and the relative path $(\varphi_h^s)^{-1}\varphi_k^s$ is generated by
\begin{equation}\label{eq:relative-flow-hamiltonian}
e^{-\kappa(\varphi_h^s)}(k_s-h_s)\circ\varphi_h^s.
\end{equation}
For every fixed smooth function $g$,
\begin{equation}\label{eq:weighted-transport}
\frac{d}{ds}\bigl(e^{-\kappa(\varphi_h^s)}g\circ\varphi_h^s\bigr)=e^{-\kappa(\varphi_h^s)}\{h_s,g\}_\alpha\circ\varphi_h^s.
\end{equation}
In particular, for autonomous $f,g$ and $t\ge0$,
\begin{equation}\label{eq:relative-flow-length}
\begin{split}
\ell_\alpha\bigl(\{\varphi_f^{-s}\varphi_g^s\}_{0\le s\le t}\bigr)&\le\norm{g-f}_\infty\int_0^t e^{s\norm{R_\alpha f}_\infty}\,ds \le t e^{t\norm{R_\alpha f}_\infty}\norm{g-f}_\infty,
\end{split}
\end{equation}
and, for every $x\in Y$,
\begin{equation}\label{eq:transport-growth}
|g(\varphi_f^t(x))|\le e^{t\norm{R_\alpha f}_\infty}\bigl(|g(x)|+t\norm{\{f,g\}_\alpha}_\infty\bigr).
\end{equation}
\end{lemma}

\begin{proof} Cartan's formula and the contact Hamiltonian equations yield $\mathcal L_{X_{h_s}}\alpha=(R_\alpha h_s)\alpha$. Hence, differentiating $(\varphi_h^s)^*\alpha=e^{\kappa(\varphi_h^s)}\alpha$ proves \eqref{conformal}. \rev{Next, to compare the two flows, note that the time-dependent vector field generating} $(\varphi_h^s)^{-1}\varphi_k^s$ is $(\varphi_h^s)^{-1}_*(X_{k_s}-X_{h_s})$. Evaluating $\alpha$ on this vector field gives \eqref{eq:relative-flow-hamiltonian}. The product rule now yields
\[ \frac{d}{ds}\bigl(e^{-\kappa(\varphi_h^s)}g\circ\varphi_h^s\bigr)=e^{-\kappa(\varphi_h^s)}\bigl(X_{h_s}g-gR_\alpha h_s\bigr)\circ\varphi_h^s, \]
which is \eqref{eq:weighted-transport}.

\rev{For the length estimate, specialize to autonomous} $f,g$\rev{.} \rev{The norm of} \eqref{eq:relative-flow-hamiltonian} is at most $e^{s\norm{R_\alpha f}_\infty}\norm{g-f}_\infty$. Therefore, integration gives \eqref{eq:relative-flow-length}. \rev{To obtain the growth estimate, integrate} \eqref{eq:weighted-transport} with $h_s=f$ \rev{and then multiply by} $e^{\kappa(\varphi_f^t)(x)}$\rev{.} \rev{This gives}
\[ g(\varphi_f^t(x))=e^{\kappa(\varphi_f^t)(x)}g(x)+\int_0^t e^{\kappa(\varphi_f^t)(x)-\kappa(\varphi_f^s)(x)}\{f,g\}_\alpha(\varphi_f^s(x))\,ds. \]
By \eqref{conformal}, the exponent in the integrand has absolute value at most $(t-s)\norm{R_\alpha f}_\infty$. Hence, taking absolute values proves \eqref{eq:transport-growth}. \end{proof}

\begin{proof}[Proof of Theorem \ref{prop-1}] We first compare each approximating flow with its fixed limiting flow. By \eqref{eq:relative-flow-length},
\[ |(\varphi_f^1)^{-1}\varphi_{f_j}^1|_\alpha\le\ell_\alpha\bigl(\{\varphi_f^{-s}\varphi_{f_j}^s\}_{0\le s\le1}\bigr)\le e^{\norm{R_\alpha f}_\infty}\norm{f_j-f}_\infty\longrightarrow0. \]
The same conclusion holds for $g_j$ and $g$. \rev{In both comparisons, the conformal estimate involves only the fixed limiting Hamiltonian. Thus no derivative bound on} $f_j$ \rev{or} $g_j$ \rev{is needed.}

\rev{Next, we use these comparisons to pass commutation to the limit. Write} $a=\varphi_f^1$, $b=\varphi_g^1$, $a_j=\varphi_{f_j}^1$, and $b_j=\varphi_{g_j}^1$. The bracket hypothesis gives $[X_{f_j},X_{g_j}]=0$, and hence
\begin{equation} \label{commute}
a_jb_j=b_ja_j.
\end{equation}
The factorization
\[ (ab)^{-1}a_jb_j =\bigl(b^{-1}(a^{-1}a_j)b\bigr)(b^{-1}b_j) \]
and \eqref{eq:conjugation-bound} imply
\[ |(ab)^{-1}a_jb_j|_\alpha \le e^{\max_Y\kappa(b^{-1})}|a^{-1}a_j|_\alpha +|b^{-1}b_j|_\alpha\longrightarrow0. \]
\rev{Interchanging} $a$ \rev{and} $b$ \rev{gives the same estimate in the opposite order:} $|(ba)^{-1}b_ja_j|_\alpha\to0$. By the triangle inequality, symmetry, and \eqref{commute},
\[ |(ab)^{-1}ba|_\alpha \le |(ab)^{-1}a_jb_j|_\alpha+|(ba)^{-1}b_ja_j|_\alpha \longrightarrow0. \]
Therefore, nondegeneracy yields $ab=ba$.

\rev{This proves commutation at time one. To obtain commutation at arbitrary times, apply the same argument to} $tf_j$ and $sg_j$, which converge to $tf$ and $sg$ and still have zero bracket, for arbitrary fixed $s,t\in\R$. Hence, $\varphi_f^t\varphi_g^s=\varphi_g^s\varphi_f^t$ for all $s,t$. Differentiating this identity gives $[X_f,X_g]=0$, and thus $\{f,g\}_\alpha=\alpha([X_f,X_g])=0$. \end{proof}

\section{Proof of Theorem \ref{prop-2}}\label{sec:proof-B}

\rev{The absence of common zeros allows us to rewrite the contact bracket using polar variables. We first establish the required identity. For a smooth pair with no common zero, put}
\[ \rho=\sqrt{f^2+g^2},\qquad u=\frac{f+\sqrt{-1}g}{\rho}=e^{\sqrt{-1}\theta}, \]
where $\theta$ is circle-valued and admits real-valued local lifts.

\begin{lemma} \label{lem-1} With the notation above, the contact vector field satisfies $\rrev{X_\rho \theta}=\{f,g\}_\alpha/\rho$ for every real-valued local lift of $\theta$. Therefore, $\rrev{X_\rho u}=\sqrt{-1}u\{f,g\}_\alpha/\rho$ on $Y$. \end{lemma}

\begin{proof} \rev{Using the polar variables \keepmath{$\rho, \theta$}, we have} $f = \rho \cos \theta$ and $g = \rho \sin \theta$. Hence, $f dg - g df = \rho^2 d\theta$, which implies that
\[  d\theta (X_\rho) = \frac{ f\rrev{X_\rho g} - g \rrev{X_{\rho} f}}{\rho^2}. \]
\rev{To evaluate the numerator, we first express the new contact vector field in terms of the original ones:}
\begin{equation} \label{X-rho}
X_\rho = \frac{f}{\rho} X_f + \frac{g}{\rho} X_g.
\end{equation}
\rev{Indeed, we verify the defining equations of the contact Hamiltonian vector field:}
\[ \alpha\left(\frac{f}{\rho} X_f + \frac{g}{\rho} X_g\right) = \frac{f^2 + g^2}{\rho} = \rho \]
and
\begin{align*}
d\alpha\left(\frac{f}{\rho} X_f + \frac{g}{\rho} X_g, \cdot\right) & = d\alpha\left(\frac{f}{\rho} X_f, \cdot \right) + d\alpha \left(\frac{g}{\rho} X_g, \cdot \right) \\
 & = \frac{f}{\rho}\left(df(R_{\alpha}) \alpha - df\right)(\cdot) + \frac{g}{\rho}\left(dg(R_{\alpha}) \alpha - dg\right)(\cdot)\\
 & = \left(\frac{fdf+gdg}{\rho}\right)(R_{\alpha})\alpha(\cdot) - \left(\frac{fdf+gdg}{\rho}\right)(\cdot).
\end{align*}
Since $d\rho = \frac{fdf+gdg}{\rho}$, we get
\[ d\alpha\left(\frac{f}{\rho} X_f + \frac{g}{\rho} X_g, \cdot\right) = d\rho(R_{\alpha}) \alpha(\cdot) - d\rho (\cdot).  \]
This verifies (\ref{X-rho}). \rev{Applying this identity to the first component gives}
\begin{align*}
\rrev{X_\rho f} &= \frac{f}{\rho}\rrev{X_f f} + \frac{g}{\rho}\rrev{X_g f} \\
 & = \frac{f}{\rho} (f\rrev{R_\alpha f}) +  \frac{g}{\rho} \left(f\rrev{R_\alpha g} - \{f, g\}_{\alpha} \right)\\
 & = \frac{f^2\rrev{R_\alpha f} + gf\rrev{R_\alpha g}}{\rho} -  \frac{g}{\rho}\{f, g\}_{\alpha}\\
 & = f \left(\frac{f\rrev{R_\alpha f} + g\rrev{R_\alpha g}}{\rho} \right) - \frac{g}{\rho}\{f, g\}_{\alpha} = f \rrev{R_\alpha \rho} - \frac{g}{\rho}\{f, g\}_{\alpha}.
\end{align*}
Here, the second equality uses $\{f,f\}_{\alpha} = \rrev{X_f f} - f \rrev{R_\alpha f}$ and $-\{f,g\}_{\alpha} = \rrev{X_g f} - f\rrev{R_\alpha g}$ \rev{(both follow from the definition in (\ref{dfn-CBP}))}. Similarly, we have $\rrev{X_{\rho} g} = g \rrev{R_\alpha \rho} + \frac{f}{\rho}\{f,g\}_{\alpha}$. \rev{Substituting both component identities into the expression above, we obtain}
\begin{align*}
d\theta (X_\rho) & = \frac{ f\rrev{X_\rho g} - g \rrev{X_{\rho} f}}{\rho^2} \\
 & = \frac{f\left(g \rrev{R_\alpha \rho} + \frac{f}{\rho}\{f,g\}_{\alpha} \right) - g\left(f \rrev{R_\alpha \rho} - \frac{g}{\rho}\{f, g\}_{\alpha}\right)}{\rho^2}\\
 & = \frac{\rho \{f, g\}_{\alpha}}{\rho^2} = \frac{\{f, g\}_{\alpha}}{\rho}.
\end{align*}
\rev{Finally, for \keepmath{$u = e^{\sqrt{-1} \theta}$}, we have}
\[ \rrev{X_{\rho} u} = du(X_{\rho}) = \sqrt{-1} e^{\sqrt{-1} \theta} d\theta(X_{\rho}) =\sqrt{-1} u \cdot \frac{\{f, g\}_{\alpha}}{\rho}. \]
This completes the proof. \end{proof}

To apply the polar identity in the proof of Theorem \ref{prop-2}, \rev{we also need a lower-semicontinuity property for the directional derivative \keepmath{$\rrev{X_f g}$}.} \rev{The key point is that a decrease in the directional derivative would force displacement, whereas displacement energy is positive for the norm} $|\cdot|_{\alpha}$. Explicitly, \rev{for any nonempty open subset \keepmath{$U \subset Y$}, we have}
\begin{equation} \label{disp-energy}
e_{\alpha}(U)=\inf_{\varphi(U)\cap U=\emptyset}|\varphi|_\alpha>0\text{\rev{.}}
\end{equation}
This is the displacement energy already defined in \eqref{eq:intro-displacement-energy}, since taking the infimum first over paths with a fixed endpoint and then over displacing endpoints gives the same infimum. The infimum of the empty set is $+\infty$. Positivity for a nonempty open set follows from the energy-capacity inequality in \cite{Shelukhin-cont}.

\begin{lemma} \label{lem-displace} Let $(Y,\xi=\ker\alpha)$ be closed, and let $f_i,g_i\in C^\infty(Y)$ converge uniformly to $f,g\in C^\infty(Y)$. Then
\[ \norm{\rrev{X_f g}}_\infty\le\liminf_{i\to\infty}\norm{\rrev{X_{f_i} g_i}}_\infty. \]
\end{lemma}

\begin{proof} We argue by contradiction. If the conclusion fails, there is a subsequence, still indexed by $i$, such that
\[ \max_Y |\rrev{X_{f_i} g_i}| \leq A < \max_Y |\rrev{X_f g}| \,\,\,\,\mbox{for all $i$} \]
\rev{for a constant \keepmath{$A$} independent of \keepmath{$i$}.} Assume $\rrev{\max_Y}|\rrev{X_f g}| = (dg(X_f))(x_0)$ for some point $x_0 \in Y$ \rev{(replacing \keepmath{$g$} and all \keepmath{$g_i$} by their negatives if necessary)}. By definition of the directional derivative, there exists a sufficiently small $t>0$ such that
\[ \Delta:= g(\varphi_f^t(x_0)) - g(x_0) - At >0. \]
\rev{To extend this strict inequality to nearby points, continuity of} \(g\) at \(x_0\) gives an open neighborhood \(U_1\) of \(x_0\) such that $|g(z)-g(x_0)|<\frac{\Delta}{4}$ for any $z\in U_1$; continuity of \(g\circ\varphi_f^t\) at \(x_0\) gives an open neighborhood \(U_2\) of \(x_0\) such that $|g(\varphi_f^t(z))-g(\varphi_f^t(x_0))| <\frac{\Delta}{4}$ for any $z\in U_2$. \rev{To ensure that both estimates hold, set} \(U=U_1\cap U_2\), which is an open neighborhood of \(x_0\). For arbitrary \(x,y\in U\), we have
\begin{align*}
g(\varphi_f^t(y))-g(x) &=g(\varphi_f^t(x_0))-g(x_0) +\bigl[g(\varphi_f^t(y))-g(\varphi_f^t(x_0))\bigr] -\bigl[g(x)-g(x_0)\bigr]\\
 &>At+\Delta-\frac{\Delta}{4}-\frac{\Delta}{4}\\
 &=At+\frac{\Delta}{2}.
\end{align*}
\rev{By contrast, along each approximating flow, the bound \keepmath{$\max_Y |\rrev{X_{f_i} g_i}| \leq A$} implies that, for any \keepmath{$x \in Y$}, we have} $g_i(\varphi_{f_i}^t(x)) - g_i(x) \leq At$ for all $i$. \rev{These different growth bounds force the two images of the neighborhood to separate. \rrev{Therefore},} $(\varphi^t_f)^{-1} \varphi_{f_i}^t$ must displace $U$ for every sufficiently large $i$. Indeed, otherwise $\varphi_{f_i}^t(U) \cap \varphi^t_f(U) \neq \emptyset$ \rev{for some sufficiently large \keepmath{$i$}, as specified below.} \rev{Such an intersection would give points} $x, y \in U$ such that $\varphi_{f_i}^t(x) = \varphi_{f}^t(y)$, which implies that
\begin{align*}
g_i(\varphi_{f_i}^t(x)) - g_i(x) &= g_i(\varphi_{f}^t(y)) - g_i(x) \\
 & \geq g(\varphi_{f}^t(y)) - g(x) - 2 \max_Y |g_i - g| \\
 & \geq At + \frac{\Delta}{2} - 2 \max_Y |g_i - g|.
\end{align*}
When $i$ is sufficiently large (so that $\max_Y |g_i - g| \leq \frac{\Delta}{8}$), we get
\[ At \geq g_i(\varphi_{f_i}^t(x)) - g_i(x) \geq At + \frac{\Delta}{2} - \frac{\Delta}{4} = At + \frac{\Delta}{4}, \]
which is a contradiction. \rev{Thus displacement of \keepmath{$U$} is unavoidable, and positivity of its energy implies that} $|(\varphi^t_f)^{-1} \varphi_{f_i}^t|_{\alpha} \geq e_{\alpha}(U)$ for all sufficiently large $i$. \rev{However, uniform convergence makes the relative flow arbitrarily short: by (\ref{eq:relative-flow-length}), we have} $|(\varphi^t_f)^{-1} \varphi_{f_i}^t|_{\alpha}\to 0$ as $i \to \infty$ \rev{by the hypothesis \keepmath{$f_i \xrightarrow{C^0} f$}.} This contradiction completes the proof.\end{proof}

\begin{remark} By definition (\ref{dfn-CBP}), the contact Poisson bracket $\{f, g\}_{\alpha}$ \rev{is the difference of the two terms \keepmath{$\rrev{X_f g}$} and \keepmath{$g\rrev{R_\alpha f}$}, so control of the terms separately does not directly imply lower semicontinuity of the \keepmath{$C^0$}-norm of \keepmath{$\{f,g\}_{\alpha}$}.} Therefore, we use the transformation in Lemma \ref{lem-1}, which requires the absence of common zeros. \end{remark}

\rev{We now apply lower semi-continuity to the directional derivative in the polar identity.}

\begin{proof}[Proof of Theorem \ref{prop-2}] Since $\min_Y\rho>0$ and $\rho_i:=\sqrt{f_i^2+g_i^2}$ converges uniformly to $\rho$, we have $\rho_i>0$ for all sufficiently large $i$. \rev{After discarding the finitely many exceptional indices, we may therefore use polar variables for every approximating pair. For each pair} $(f_i,g_i)$, \rev{write \keepmath{$(\rho_i, u_i)$} for the polar variables from Lemma \ref{lem-1}.} Recall that the map $u: Y \to \C \simeq \R^2$, like each $u_i$, is globally defined. \rev{Hence, the hypotheses \keepmath{$f_i \xrightarrow{C^0} f$} and \keepmath{$g_i \xrightarrow{C^0} g$} imply that} $\rho_i \xrightarrow{C^0} \rho$ (valued in $\R$) and $u_i \xrightarrow{C^0} u$ (valued in $\R^2$).

\rev{However, Lemma \ref{lem-displace} is stated for real-valued functions. To apply it to \keepmath{$u$} and \keepmath{$u_i$}, we reduce to \keepmath{$\R$}-valued functions using a fixed linear projection.} Suppose $\max_Y|\rrev{X_\rho u}| = |\rrev{(X_\rho u)}(x_0)|$ where $|\cdot| = |\cdot|_{\R^2}$. \rev{If \keepmath{$\rrev{(X_\rho u)}(x_0) \neq 0 \in \R^2$}, choose a linear function \keepmath{$\ell: \R^2 \to \R$} given by projection onto the direction \keepmath{$\rrev{(X_\rho u)}(x_0)$}, that is,}
\[ \ell(v) = \left<v, \frac{\rrev{(X_\rho u)}(x_0)}{|\rrev{(X_\rho u)}(x_0)|} \right>. \]
\rev{By construction,} $\ell$ has operator norm one and $\ell(\rrev{(X_\rho u)}(x_0)) = |\rrev{(X_\rho u)}(x_0)|$. Then for the function $\ell \circ u: Y  \to \R$, we have $X_\rho(\ell \circ u) = \ell(\rrev{X_\rho u})$ since $\ell$ is linear. Moreover, for any $y \in Y$,
\[ |X_\rho(\ell \circ u)(y)| \leq |\rrev{(X_\rho u)}(y)| \leq \max_Y |\rrev{X_\rho u}| \]
\rev{whereas equality is attained at the point} $x_0 \in Y$\rev{:}
\[ |X_\rho(\ell \circ u)(x_0)| = \ell(\rrev{(X_\rho u)}(x_0)) = |\rrev{(X_\rho u)}(x_0)|. \]
\rrev{Therefore}, we have $\max_Y |X_\rho(\ell \circ u)|_{\R} = \max_Y|\rrev{X_\rho u}|_{\R^2}$. Since $u_i \xrightarrow{C^0} u$ implies $\ell \circ u_i \xrightarrow{C^0} \ell \circ u$, Lemma \ref{lem-displace} applies and yields
\begin{align*}
\max_Y |\rrev{X_\rho u}| & =  \max_Y |X_\rho(\ell \circ u)| \\
 & \leq \liminf_i \max_Y |X_{\rho_i}(\ell \circ u_i)| \leq \liminf_i \max_Y |\rrev{X_{\rho_i} u_i}|.
\end{align*}

\rev{The same inequality is immediate if \keepmath{$\rrev{(X_\rho u)}(x_0)=0$}. It remains to show that its right-hand side vanishes. Since \keepmath{$\rho>0$} on \keepmath{$Y$} and \keepmath{$Y$} is closed,} there exists some $\delta>0$ such that $\min_Y \rho \geq \delta$. Since $\rho_i \xrightarrow{C^0} \rho$, for all sufficiently large $i$, we have $\min_Y \rho_i \geq \frac{\delta}{2} >0$. Hence, by the second conclusion in Lemma \ref{lem-1},
\[ |\rrev{X_{\rho_i} u_i}| = \left|\sqrt{-1} u_i \cdot \frac{\{f_i, g_i\}_{\alpha}}{\rho_i} \right| \leq |u_i| \cdot \frac{|\{f_i, g_i\}_{\alpha}|}{\delta/2}.\]
In particular,
\[ \max_Y |\rrev{X_{\rho_i} u_i}| \leq \max_Y \frac{|\{f_i, g_i\}_{\alpha}|}{\delta/2} \to 0 \,\,\,\,\mbox{as $i \to \infty$} \]
\rev{because \keepmath{$|u_i| = 1$} and the brackets \keepmath{$\{f_i, g_i\}_{\alpha}$} converge uniformly to zero.} \rev{The lower-semicontinuity inequality above therefore gives} $\max_Y |\rrev{X_\rho u}|=0$ and so $\rrev{X_{\rho} u}=0$. Hence, Lemma \ref{lem-1} gives $\{f, g\}_{\alpha} =0$, which completes the proof.  \end{proof}


\section{Proof of Theorem \ref{prop-3}}\label{sec:proof-C}

Fix a Riemannian metric $\mathfrak g$ on $Y$, and let $B_{\mathfrak g}(p,\delta)$ be the open metric ball with center $p$ and radius $\delta$. \rev{We first need a lower bound for displacement energy that keeps track of the radius.}

\begin{lemma} \label{lem-refine-dis} For each $p\in Y$, there is a constant $c_p>0$ such that $e_\alpha(B_{\mathfrak g}(p,\delta))\ge c_p\delta^2$ for all sufficiently small $\delta>0$. \end{lemma}

\begin{proof} \rev{Recall the energy-capacity inequality from \cite{Shelukhin-cont}:}
\[ e_{\alpha}(B_{\mathfrak g}(p, \delta)) \geq \frac{1}{4} \sup_{A \subset \R_{>0} \times B_{\mathfrak g}(p, \delta)} \frac{w_G(A)}{h_{\alpha}(A)} \]
\rrev{The supremum runs over compact subsets. Here} $\R_{>0} \times B_{\mathfrak g}(p, \delta)$ can be viewed as the preimage of $B_{\mathfrak g}(p, \delta)$ under the projection from the symplectization $\R_{>0} \times Y \to Y$, \rev{while \keepmath{$w_G(A)$} denotes the Gromov width of \keepmath{$A$} and \keepmath{$h_{\alpha}(A)$} is the supremum of its \keepmath{$\R_{>0}$}-coordinate.} \rev{Thus it suffices, for each sufficiently small \keepmath{$\delta>0$}, to construct a compact subset} $A_\delta \subset \R_{>0} \times B_{\mathfrak g}(p, \delta)$ such that
\begin{equation} \label{delta-ball}
w_G(A_\delta) \geq C\delta^2 \,\,\,\,\mbox{and} \,\,\,\, h_{\alpha}(A_{\delta}) \leq C'
\end{equation}
for some $C, C'$ independent of $\delta$.

\rev{We construct these subsets by rescaling in a fixed Darboux ball} $Q$ near $p$, identified with a neighborhood of the origin in $\R^{2n+1}$ in strict Darboux coordinates $(x_1,y_1,\ldots,x_n,y_n,z)$ for the fixed form, so that $\alpha=dz-\sum_{j=1}^n y_j\,dx_j$ and $p$ corresponds to the origin. Shrink $Q$ to a bounded Euclidean ball centered at the origin. \rev{On \keepmath{$Q$}, we have \keepmath{$d_{\mathfrak g}(p, q) \leq C |q|$}, where \keepmath{$|q|$} denotes the Euclidean norm of \keepmath{$q$} and \keepmath{$C$} is a constant depending on \keepmath{$p$} and the metric \keepmath{$\mathfrak g$}.} Consider the rescaling $\sigma_s(x_1, ..., y_n, z) = (s x_1, ..., s y_n, s^2z)$, which is a contactomorphism with respect to the standard contact 1-form $dz - \sum y_i dx_i$. Then for any $q \in Q$ and $0<s<1$, we have
\[ d_{\mathfrak g}(p, \sigma_s(q)) \leq C|\sigma_s(q)| \leq C \cdot s|q| \leq C \cdot s \sup_{q \in Q} |q|. \]
The second inequality comes from the following estimate: if $q = (x_1, ..., y_n, z)$, then
\[ |\sigma_s(q)|^2 = s^2(x_1^2 + ... + y_n^2) + s^4 z^2 \leq s^2(x_1^2 + ...+ z^2) = s^2 |q|^2. \]
Choose a constant $a>0$ such that $a \cdot C \sup_{q \in Q} |q|<1$. Then
\[ d_{\mathfrak g}(p, \sigma_{a\delta}(q)) < \delta \,\,\,\,\mbox{for any $q \in Q$}. \]
In particular, $\sigma_{a\delta}(Q) \subset B_{\mathfrak g}(p, \delta)$.

\rev{We now lift this rescaling to the symplectization to control the width and height simultaneously. Fix a compact subset \keepmath{$A \subset (1/2, 2) \times Q$} containing a symplectic ball of radius \keepmath{$R$}.} Since $\sigma_s^*\alpha=s^2\alpha$, the map $\Phi_s(r,q)=(r,\sigma_s(q))$ satisfies $\Phi_s^*d(r\alpha)=s^2d(r\alpha)$. Thus, for the conformally symplectic map on the symplectization $\Phi_{a\delta}: (r,q) \mapsto (r, \sigma_{a\delta}(q))$, $\Phi_{a\delta}(A)$ contains a symplectic ball of radius $a\delta R$. Indeed, compose a symplectic embedding of the radius-$R$ ball into $A$ with dilation by $1/(a\delta)$ on its domain and then with $\Phi_{a\delta}$; the two conformal factors cancel. \rev{\rrev{Therefore},} $\Phi_{a\delta}(A)$ \rev{is the desired} $A_{\delta}$. \rev{The inclusion established above gives} $\Phi_{a\delta}(A) \subset \R_{>0} \times B_{\mathfrak g}(p, \delta)$\rev{, and the width and height satisfy}
\[ w_G(\Phi_{a\delta}(A)) \geq \pi (a\delta R)^2 \,\,\,\,\mbox{and}\,\,\,\, h_{\alpha}(\Phi_{a\delta}(A)) = h_{\alpha}(A). \]
This verifies (\ref{delta-ball}) with $C = \pi a^2 R^2$ and $C' =h_{\alpha}(A)$. Note that both $C$ and $C'$ depend on $Q$ \rev{(and hence on \keepmath{$p$} and the metric \keepmath{$\mathfrak g$})}. Therefore, $c_p=C/(4C')$ gives the required bound for all sufficiently small $\delta>0$. \end{proof}

\rev{We now compare two flows on shrinking neighborhoods. The comparison will force displacement, and the local energy bound will then yield a contradiction.}

\begin{proof}[Proof of Theorem \ref{prop-3}] For this proof, write $\rrev{\epsilon}_i=\max\{\norm{f_i-f}_\infty,\norm{g_i-g}_\infty\}$. We argue by contradiction. If the conclusion fails, there is a subsequence, still indexed by $i$, such that
\[ \max_Y |\{f_i, g_i\}_{\alpha} | \leq A < \max_Y |\{f, g\}_{\alpha}| \,\,\,\,\mbox{for all $i$} \]
\rev{for a constant \keepmath{$A\geq 0$} independent of \keepmath{$i$}.} By replacing $f$ and all $f_i$ with their negatives if necessary, one can assume that there exists some $p \in Y$ such that $\{f, g\}_{\alpha}(p)>A$. In particular, the set $\{\{f, g\}_{\alpha} >A\} (\subset Y)$ \rev{is open and nonempty.} \rev{Moreover, \keepmath{$g$} cannot vanish identically on this open set: otherwise, \keepmath{$\{f, g\}_{\alpha} = \rrev{X_f g} - g \rrev{R_\alpha f} =0$} there, a contradiction.} \rev{Thus we can choose a point in this open set, still denoted by} $p \in Y$, \rev{such that}
\[ \{f, g\}_{\alpha}(p)>A \,\,\,\,\mbox{and}\,\,\,\, g(p) \neq 0. \]

\rev{To turn the bracket inequality into a directional-derivative inequality at this point, set} $c:= \frac{f(p)}{g(p)}$ and consider the function $u:= cg-f$. This linear combination preserves the relevant bracket while making its Reeb term vanish at $p$. Explicitly,
\[ \{g, u\}_{\alpha} = c\{g, g\}_{\alpha} - \{g, f\}_{\alpha} = \{f, g\}_{\alpha}. \]
Hence, $\{g, u\}_{\alpha}(p) > A$. Since $u(p) = 0$ by definition, we have $\{g, u\}_{\alpha}(p) = \rrev{(X_g u)}(p) (>A)$. Similarly, for any approximating pair $(f_i, g_i)$, one can define $u_i := c g_i - f_i$ (with the same constant $c$). Thus, we still have $\{g_i, u_i\}_{\alpha} = \{f_i, g_i\}_{\alpha}$ \rev{(with \keepmath{$C^0$}-norm at most \keepmath{$A$} by our contradiction assumption)}. However, $u_i(p)$ is not necessarily zero. \rev{Nevertheless, its size is controlled by the approximation error. Indeed, since}
\[ \max_Y |u_i - u| \leq |c| \left(\max_Y |g_i - g|\right) + \left(\max_Y |f_i -f|\right) \leq (1+|c|) \rrev{\epsilon}_i \]
we have $|u_i(p)| \leq (1+|c|) \rrev{\epsilon}_i$, which tends to zero.

\rev{We use these facts to compare the values of \keepmath{$u_i$} along the flows \keepmath{$\varphi_g^t$} and \keepmath{$\varphi_{g_i}^t$} near \keepmath{$p$}. The fixed flow gives a lower bound, whereas the bracket bound controls growth along the approximating flow.} Consider $B_{\mathfrak g}(p, \delta)$ for sufficiently small $\delta$. For any $x \in B_{\mathfrak g}(p, \delta)$, we have the local expansion (for sufficiently small $t$, subject to further conditions specified below):
\[ u_i\left(\varphi_g^t(x)\right) = u(x) + t \rrev{(X_g u)}(x) + O(t^2 + \rrev{\epsilon}_i) =t \rrev{(X_g u)}(p) + O(\delta + t^2 + \rrev{\epsilon}_i). \]
Here, we use the estimate that $\rrev{(X_g u)}(x) = \rrev{(X_g u)}(p) + O(\delta)$. \rrev{Since $g$ and $u$ are fixed and $\norm{u_i-u}_\infty\le(1+|c|)\rrev{\epsilon}_i$, the constants in these expansions can be chosen independently of $i$ for all sufficiently small $t,\delta>0$. Moreover, $u(p)=0$ gives $|u_i(x)|\le C_0(\delta+\rrev{\epsilon}_i)$ on the ball, for a fixed constant $C_0>0$.} Hence,
\begin{equation} \label{uig}
\inf_{\varphi_g^t(B_{\mathfrak g}(p, \delta))} u_i \geq t \rrev{(X_g u)}(p) - C(\delta + t^2 + \rrev{\epsilon}_i)
\end{equation}
for some constant $C>0$. On the other hand, \eqref{eq:transport-growth}, applied to the Hamiltonian $g_i$ and the function $u_i$, gives
\begin{equation} \label{Gronwall}
|u_i(\varphi_{g_i}^t(x))|\le e^{t\norm{R_\alpha g_i}_\infty}\bigl(|u_i(x)|+At\bigr)
\end{equation}
for every $x\in B_{\mathfrak g}(p,\delta)$, since $\{g_i,u_i\}_\alpha=\{f_i,g_i\}_\alpha$ has norm at most $A$.
\rev{To compare this exponential bound with the preceding linear expansion, suppose} $t$ \rev{is so small that} $\max_Y |\rrev{R_\alpha g_i}| t \leq 1$. Then $e^{\max_Y |\rrev{R_\alpha g_i}| t} \leq 1+ t e \max_Y |\rrev{R_\alpha g_i}|$. \rrev{Using the initial-value bound above and $e^{t\norm{R_\alpha g_i}_\infty}\le e$, we obtain}
\begin{equation} \label{uigi}
\sup_{\varphi_{g_i}^t(B_{\mathfrak g}(p, \delta))} u_i  \leq At + C\left(\delta + \rrev{\epsilon}_i + \max_Y |\rrev{R_\alpha g_i}| t^2\right)
\end{equation}
for some constant $C>0$ (enlarged if necessary so that it also works in (\ref{uig})). \rev{Subtracting the upper bound from the lower bound now gives}
\begin{align*}
\inf_{\varphi_g^t(B_{\mathfrak g}(p, \delta))} u_i - \sup_{\varphi_{g_i}^t(B_{\mathfrak g}(p, \delta))} u_i  & \geq t (\rrev{(X_g u)}(p) - A) - C\left(2\delta + 2\rrev{\epsilon}_i + \left(1+ \max_Y |\rrev{R_\alpha g_i}|\right) t^2\right)  \\
 & \geq t (\rrev{(X_g u)}(p) - A) - C' \left(\delta + \rrev{\epsilon}_i + \left(1+ \max_Y |\rrev{R_\alpha g_i}|\right) t^2\right)
\end{align*}
where $C' = 2C$. In particular, if this lower bound is strictly positive, then $\varphi_g^t(B_{\mathfrak g}(p, \delta))$ is disjoint from $\varphi_{g_i}^t(B_{\mathfrak g}(p, \delta))$.

\rev{Whenever this positivity condition holds,} $(\varphi_g^t)^{-1} \varphi_{g_i}^t$ displaces $B_{\mathfrak g}(p, \delta)$. By \eqref{eq:relative-flow-length} and Lemma~\ref{lem-refine-dis}, we have
\begin{equation} \label{two-side}
c_p \delta^2 \leq |(\varphi_g^t)^{-1} \varphi_{g_i}^t|_{\alpha} \leq te^{t \max_Y |\rrev{R_\alpha g}|} \max_Y |g_i - g| \leq C' t \rrev{\epsilon}_i
\end{equation}
with $C'>0$ enlarged if necessary and $c_p>0$ as in the lemma. Consider
\[ t_i = \sqrt{\frac{\rrev{\epsilon}_i}{1 + \max_Y |\rrev{R_\alpha g_i}|}} \,\,\,\,\mbox{and}\,\,\,\, \delta_i = \frac{\rrev{(X_g u)}(p) - A}{4C'} t_i. \]
Since $\rrev{\epsilon}_i \to 0$ and $\rrev{\epsilon}_i \max_Y |\rrev{R_\alpha g_i}| \to 0$ by our assumption, we have
\[ 0<t_i\le\sqrt{\rrev{\epsilon}_i}\longrightarrow0, \qquad t_i\max_Y|\rrev{R_\alpha g_i}| \le\sqrt{\rrev{\epsilon}_i\max_Y|\rrev{R_\alpha g_i}|}\longrightarrow0. \]
Here $\rrev{\epsilon}_i>0$, since equality $\rrev{\epsilon}_i=0$ would already contradict the strict bracket bound defining the subsequence. \rev{Thus, for this choice of times, the exponential estimate used in \eqref{uigi} holds} for all sufficiently large $i$. \rev{It remains to verify separation of the two images. With these choices, the lower bound for the difference of the values of \keepmath{$u_i$} becomes}
\begin{align*}
t_i (\rrev{(X_g u)}(p) - A) & - C' \left(\delta_i + \rrev{\epsilon}_i  +  \left(1+ \max_Y |\rrev{R_\alpha g_i}|\right) t_i^2\right)  \\
 &= 4C'\delta_i - C'(\delta_i + \rrev{\epsilon}_i + \rrev{\epsilon}_i) \\
 & = 3 C' \delta_i - 2\rrev{\epsilon}_i C'\\
 & = \frac{3}{4} (\rrev{(X_g u)}(p) - A) \cdot t_i - 2\rrev{\epsilon}_i C' = t_i \left(\frac{3}{4} (\rrev{(X_g u)}(p) - A) - 2 C' \frac{\rrev{\epsilon}_i}{t_i}\right).
\end{align*}
Since $\frac{\rrev{\epsilon}_i}{t_i} = \sqrt{\rrev{\epsilon}_i (1+ \max_Y |\rrev{R_\alpha g_i}|)} \to 0$, \rev{the lower bound above is strictly positive for sufficiently large \keepmath{$i$}, so the displacement argument applies.} \rev{The resulting displacement, however, must satisfy both bounds in (\ref{two-side}). Hence we have} $c_p \delta_i^2 \leq C' t_i \rrev{\epsilon}_i$, which is
\[ c_p \left(\frac{\rrev{(X_g u)}(p) - A}{4C'} \right)^2 t_i^2 \leq C' t_i \rrev{\epsilon}_i. \]
This is impossible since, after dividing by $t^2_i$, we get
\[ 0< c_p \left(\frac{\rrev{(X_g u)}(p) - A}{4C'} \right)^2 \leq C' \frac{\rrev{\epsilon}_i}{t_i} \to 0 \,\,\,\,\mbox{as $i \to \infty$}. \]
This completes the proof. \end{proof}

\section{\rev{Counterexample to} \texorpdfstring{$C^0$}{C0}-rigidity}\label{sec:counterexample}

On $\T^3=(\R/2\pi\Z)^3$, consider $\alpha=\cos x\,dz+\sin x\,dy$. Write $c:=\cos x$ and $s:=\sin x$. Since $\alpha\wedge d\alpha=dx\wedge dy\wedge dz$, this is a contact form. Its Reeb field is $R=s\partial_y+c\partial_z$. With $V=c\partial_y-s\partial_z$, the contact Hamiltonian equations give
\begin{equation} \label{ce:coordinate-bracket}
X_h=hR-(Vh)\partial_x+h_xV,\qquad \{h,k\}_\alpha=hRk-kRh+h_xVk-(Vh)k_x.
\end{equation}
For the limiting functions
\begin{equation} \label{ce:limit-torus-recalled}
f=c\cos z,\qquad g=c\sin z,
\end{equation}
we obtain
\begin{equation} \label{ce:base-fields}
X_f=\cos z\,\partial_z-sc\sin z\,\partial_x,\qquad X_g=\sin z\,\partial_z+sc\cos z\,\partial_x.
\end{equation}
In particular, both fields have zero $y$-component. Also, $Rf=-c^2\sin z$ and $Rg=c^2\cos z$, so \eqref{ce:coordinate-bracket} gives $\{f,g\}_\alpha=c$ and then $\norm{\{f,g\}_\alpha}_\infty=1$.

\rev{To approximate these functions while canceling their nonzero bracket, introduce oscillations in the remaining coordinate. For each integer \keepmath{$j\ge1$}, set \keepmath{$A_j(x) =\sqrt{4jc^2-2sc+j^{-1}}$} and}
\begin{equation} \label{ce:seq-torus}
\begin{aligned}
u_j&=j^{-2}A_j(x)\cos(j^4y),\qquad v_j=j^{-2}A_j(x)\sin(j^4y), \qquad f_j&=f+u_j,\qquad g_j=g+v_j.
\end{aligned}
\end{equation}
The inequality $2|sc|\le2jc^2+s^2/(2j)$ implies
\begin{equation} \label{ce:Aj-bounds}
2jc^2+\frac1{2j}\le A_j^2\le4j+1+j^{-1}\le6j.
\end{equation}
Thus $A_j$ is smooth and periodic, and the integer frequency $j^4$ makes $f_j,g_j$ well-defined on $\T^3$. \rev{The upper bound also gives uniform convergence:}
\begin{equation} \label{ce:C0-est}
\max\{\norm{f_j-f}_\infty,\norm{g_j-g}_\infty\}\le\sqrt6\,j^{-3/2}\longrightarrow0.
\end{equation}

\rev{It remains to estimate the bracket of the perturbed pair. We use \eqref{ce:coordinate-bracket} to separate the oscillatory term from the mixed terms.} Since $A_jA_j'=-4jcs-c^2+s^2$,
\[ \{u_j,v_j\}_\alpha=cA_jA_j'+sA_j^2=-c+\frac{s}{j}. \]
\rev{The leading term therefore cancels the bracket of the limiting pair. Meanwhile, since the vector fields in \eqref{ce:base-fields} have no} $\partial_y$ component, the mixed terms contain no factor $j^4$. More precisely,
\begin{align*}
\{f,v_j\}_\alpha&=j^{-2}(c^2A_j-scA_j')\sin z\sin(j^4y),\\
\{u_j,g\}_\alpha&=j^{-2}(c^2A_j-scA_j')\cos z\cos(j^4y).
\end{align*}
Since $A_j(c^2A_j-scA_j')=(4j+j^{-1})c^2-sc$, adding these expressions gives
\begin{equation} \label{ce:exact-bracket}
\{f_j,g_j\}_\alpha=\frac{s}{j}+j^{-2}B_j(x)\cos(z-j^4y),\qquad B_j(x)=\frac{(4j+j^{-1})c^2-sc}{A_j(x)}.
\end{equation}
For $c\ne0$, the lower bound in \eqref{ce:Aj-bounds} gives
\[ |B_j|\le\frac{(4j+j^{-1})|c|+|s|}{\sqrt{2j}}\le3\sqrt2\,\sqrt j. \]
At $c=0$ the same estimate holds because $B_j=0$. \rev{For a lower bound, evaluate at} $x=\pi/2$, \rev{where} the bracket in \eqref{ce:exact-bracket} equals $1/j$. Therefore,
\begin{equation} \label{ce:bracket-est}
\frac1j\le\norm{\{f_j,g_j\}_\alpha}_\infty\le\frac1j+3\sqrt2\,j^{-3/2}\longrightarrow0.
\end{equation}
Together with \eqref{ce:C0-est} and $\norm{\{f,g\}_\alpha}_\infty=1$, this proves the failure of unrestricted weak rigidity and \rev{lower semicontinuity of the bracket norm.} 

\rev{Finally, we check which additional hypotheses in the rigidity theorems fail for this construction. The limiting functions have common zeros precisely where} $c=0$. \rev{The Reeb-derivative condition in Theorem~\ref{prop-3} fails as well, despite the uniform convergence.} Indeed, evaluating $f_j-f$ at $x=y=0$ gives an approximation error at least $j^{-2}\sqrt{4j+j^{-1}}$. Also,
\[ Rg_j=c^2\cos z+j^2sA_j\cos(j^4y), \]
whose value at $(x,y,z)=(\pi/4,0,\pi/2)$ is $j^2\sqrt{2j-1+j^{-1}}/\sqrt2$. Hence
\begin{align*}
\max\{\norm{f_j-f}_\infty,\norm{g_j-g}_\infty\} \cdot \norm{R_\alpha g_j}_\infty&\ge\frac1{\sqrt2}\sqrt{(4j+j^{-1})(2j-1+j^{-1})}\\
&\ge\sqrt2\,j\longrightarrow+\infty.
\end{align*}
\rev{Thus this example satisfies neither the additional hypothesis of Theorem~\ref{prop-2} nor that of Theorem~\ref{prop-3}.}

\section{Proof of Theorem \ref{cor-1}}\label{sec:proof-D}

\begin{proof} If $\liminf_j\Palpha(\vec f_j)=+\infty$, there is nothing to prove. \rev{Otherwise, pass to a subsequence realizing this liminf. The bracket functional is then bounded: after discarding finitely many terms,} there is a constant $C$ such that $\Palpha(\vec f_j)\le C$ along this subsequence.

\rev{To apply the rigidity theorem for pairs of functions, choose coefficients that attain the limiting bracket maximum. The expression defining} $\Palpha(\vec f)$ \rev{is continuous} in the coefficient vectors and the point of $Y$. Therefore, compactness gives $a,b\in[-1,1]^N$ such that, with
\[ h=\sum_{i=1}^N a_if_i,\qquad k=\sum_{i=1}^N b_if_i, \]
we have $\Palpha(\vec f)=\norm{\{h,k\}_\alpha}_\infty$. Use these same coefficients for every approximating partition:
\[ h_j=\sum_{i=1}^N a_if_{j,i},\qquad k_j=\sum_{i=1}^N b_if_{j,i}. \]
Since the coefficient vectors lie in $[-1,1]^N$,
\[ \max\{\norm{h_j-h}_\infty,\norm{k_j-k}_\infty\}\le N\max_i\norm{f_{j,i}-f_i}_\infty\longrightarrow0. \]
\rev{Uniform convergence alone is not enough for that theorem; we also need Reeb control. This follows from the partition identity: since} $\sum_i f_{j,i}=1$,
\[ \norm{R_\alpha k_j}_\infty =\rrev{\norm{\left\{\sum_{i=1}^N f_{j,i},\sum_{i=1}^N b_if_{j,i} \right\}_\alpha}_\infty} \le\Palpha(\vec f_j)\le C. \]
Therefore, the uniform convergence and this Reeb bound imply the condition in Theorem~\ref{prop-3}, which gives
\[ \Palpha(\vec f)=\norm{\{h,k\}_\alpha}_\infty \le\liminf_j\norm{\{h_j,k_j\}_\alpha}_\infty \le\liminf_j\Palpha(\vec f_j). \]
Since the subsequence was chosen to realize the original liminf, this is the desired inequality for the full sequence. \end{proof}

\section{Proof of Theorem \ref{prop-4}}\label{sec:proof-E}

\rev{We begin by comparing the flow of a sum with the product of the two flows. The following geometric estimate holds on every closed contact manifold and will also be used in the proof of Theorem~\ref{thm:spectral-cover}.}

\begingroup\samepage
\begin{lemma}\label{lem:flow-comparison} For autonomous $f,g\in C^\infty(Y)$ and $t\ge0$, the based path
\[ \eta_s=\varphi_g^{-s}\varphi_f^{-s}\varphi_{f+g}^s,\qquad 0\le s\le t, \]
represents $(\widetilde\varphi_f^{\,t}\widetilde\varphi_g^{\,t})^{-1}\widetilde\varphi_{f+g}^{\,t}$ and satisfies
\begin{equation}\label{eq:sum-flow-length}
\ell_\alpha(\eta)\le\frac{t^2}{2}e^{t(\norm{R_\alpha f}_\infty+\norm{R_\alpha g}_\infty)}\norm{\{f,g\}_\alpha}_\infty.
\end{equation}
\end{lemma}
\endgroup

\begin{proof} Apply \eqref{eq:relative-flow-hamiltonian} first to $\varphi_f^s,\varphi_{f+g}^s$, and then to $\varphi_g^s,\varphi_f^{-s}\varphi_{f+g}^s$. The resulting Hamiltonian of $\eta_s$ is
\begin{equation}\label{modify-cont-ham}
h_s=e^{-\kappa(\varphi_g^s)}\bigl(e^{-\kappa(\varphi_f^s)}g\circ\varphi_f^s-g\bigr)\circ\varphi_g^s.
\end{equation}
\rev{The weighted transport identity \eqref{eq:weighted-transport} expresses the difference in parentheses in terms of the bracket: it equals}
\[ \int_0^s e^{-\kappa(\varphi_f^r)}\{f,g\}_\alpha\circ\varphi_f^r\,dr. \]
Combining this identity with \eqref{conformal} gives
\[ \norm{h_s}_\infty\le s e^{s(\norm{R_\alpha f}_\infty+\norm{R_\alpha g}_\infty)}\norm{\{f,g\}_\alpha}_\infty. \]
Since the exponent is nonnegative, integration from $0$ to $t$ proves \eqref{eq:sum-flow-length}. The pointwise identity $\varphi_{f+g}^s=\varphi_f^s\varphi_g^s\eta_s$ gives the stated identity of lifted classes. \end{proof}

\rev{To convert this geometric comparison into an estimate for the bracket invariant, we now fix a nonzero monotone homogeneous quasi-morphism} $\mu$, normalized by \eqref{eq:mu-normalization}, and derive several further estimates based on Lemma \ref{lem:flow-comparison}. The vanishing property will be needed only in the final application to a displaceable cover.

\begin{remark} \label{rmk-normalization} A homogeneous quasi-morphism is additive on commuting elements. Indeed, if $ab=ba$, then $(ab)^m=a^mb^m$, so its defect inequality gives
\[ m\,|\mu(ab)-\mu(a)-\mu(b)|\le D(\mu) \]
for every positive integer $m$. Therefore, dividing by $m$ and passing to the limit proves the assertion. In particular, $t\mapsto\mu(\widetilde\varphi_1^{\,t})$ is additive. It is also monotone, because each forward segment of the Reeb flow is generated by a positive Hamiltonian. Since an additive monotone function on $\R$ is linear, we have
\[ \mu(\widetilde\varphi_1^{\,t})=t\mu(\widetilde\varphi_1). \]
Monotonicity shows that $\mu(\widetilde\varphi_1)\ge0$. Moreover, this number cannot be zero for a nonzero $\mu$: every smooth contact Hamiltonian $h$ on $[0,1]\times Y$ satisfies $-m \leq h(t,x) \leq m$ for some integer $m\in\Z_{\geq0}$ and all $(t,x) \in [0,1] \times Y$, which implies that $-m\mu(\widetilde\varphi_1) \leq \mu(\widetilde\varphi_h) \leq m\mu(\widetilde\varphi_1)$; then $\mu(\widetilde\varphi_1)=0$ would imply that $\mu$ is identically zero. Therefore, the positive normalization in \eqref{eq:mu-normalization} can be achieved. \end{remark} 

\rev{With this normalization, monotonicity gives control by path length. Let} $h_s$ \rev{generate a path} for $0\le s\le t$. Formula~\eqref{eq:relative-flow-hamiltonian} shows that pointwise comparison $a_s\le b_s$ gives a nonnegative relative Hamiltonian. Comparing $h_s$ with constant Hamiltonians $\min_Yh_s$ and $\max_Yh_s$, monotonicity and the Reeb normalization give
\[ \int_0^t\min_Yh_s\,ds\le\mu([\{\varphi_h^s\}_{0\le s\le t}])\le\int_0^t\max_Yh_s\,ds. \]
\rrev{If the extrema are not smooth functions of time, use smooth scalar approximations from below and above, respectively, and pass to the limit in their integrals.}
In particular,
\begin{equation} \label{stability}
\left|\mu([\{\varphi_h^s\}_{0\le s\le t}])\right|\le\ell_\alpha(\{\varphi_h^s\}_{0\le s\le t}).
\end{equation}

\rev{We also need homogeneity for arbitrary real times. For each autonomous} $f$, \rev{the function} $t\mapsto\mu(\widetilde\varphi_f^{\,t})$ is additive, by the commuting property above, and is locally bounded by $|t|\norm{f}_\infty$, by \eqref{stability} and homogeneity under inversion. Hence it is linear:
\begin{equation} \label{eq:real-homogeneity}
\mu(\widetilde\varphi_f^{\,t})=t\mu(\widetilde\varphi_f) \qquad(t\in\R).
\end{equation}
\rev{This justifies the use of arbitrary terminal times for contact Hamiltonian flows below.}

\begin{lemma}\label{lem-nu-est} For autonomous $f,g\in C^\infty(Y)$ and every $t>0$,
\begin{equation} \label{eq:mu-almost-additive}
\begin{split} &\left|\mu(\widetilde\varphi_{f+g}) -\mu(\widetilde\varphi_f)-\mu(\widetilde\varphi_g)\right| \le\frac{2D(\mu)}{t} +\frac t2 e^{t(\norm{R_\alpha f}_\infty+\norm{R_\alpha g}_\infty)} \norm{\{f,g\}_\alpha}_\infty. \end{split}
\end{equation}
\end{lemma}

\begin{proof} Put $A_t=\widetilde\varphi_f^{\,t}$, $B_t=\widetilde\varphi_g^{\,t}$, and $C_t=\widetilde\varphi_{f+g}^{\,t}$. Lemma~\ref{lem:flow-comparison} and \eqref{stability} imply
\begin{equation}\label{eq:comparison-mu}
|\mu((A_tB_t)^{-1}C_t)|\le\frac{t^2}{2}e^{t(\norm{R_\alpha f}_\infty+\norm{R_\alpha g}_\infty)}\norm{\{f,g\}_\alpha}_\infty.
\end{equation}
Applying the defect inequality first to $C_t=(A_tB_t)((A_tB_t)^{-1}C_t)$ and then to $A_tB_t$ gives
\[ |\mu(C_t)-\mu(A_t)-\mu(B_t)|\le2D(\mu)+|\mu((A_tB_t)^{-1}C_t)|. \]
Hence, combining this with \eqref{eq:comparison-mu} and \eqref{eq:real-homogeneity}, and dividing by $t$, proves \eqref{eq:mu-almost-additive}. \end{proof}

The following lemma will be used repeatedly later in the paper, so we state it separately.

\begin{lemma}\label{lem:vanishing-partition-bound} Let $\mu$ be a nonzero monotone homogeneous quasi-morphism normalized by \eqref{eq:mu-normalization}. If a smooth partition of unity $\vec f=(f_1,\ldots,f_N)$, $N\ge2$, satisfies $\mu(\widetilde\varphi_{f_i})=0$ for every $i$, then
\[ \Palpha(\vec f)\ge\frac{1}{16(N-1)^2D(\mu)}. \]
\end{lemma}

\begin{proof} Write $P=\Palpha(\vec f)$ and $h_k=\sum_{i=1}^k f_i$. By definition of $P$ and $\sum_i f_i=1$, for $2\le k\le N$ we have
\begin{equation}\label{eq:partition-controls}
\norm{R_\alpha h_{k-1}}_\infty\le P,\qquad \norm{R_\alpha f_k}_\infty\le P,\qquad \norm{\{h_{k-1},f_k\}_\alpha}_\infty\le P.
\end{equation}
For example, $R_\alpha h_{k-1}=\{1,h_{k-1}\}_\alpha$ and both inputs are combinations of the partition with coefficients in $[-1,1]$. The other inequalities follow in the same way.

By hypothesis, $\mu(\widetilde\varphi_{f_i})=0$ for every $i$. In particular, $\mu(\widetilde\varphi_{h_1})=0$. \rev{At the last partial sum,} $h_N=1$, so $\mu(\widetilde\varphi_{h_N})=1$ by normalization. \rev{To compare these endpoint values, apply Lemma~\ref{lem-nu-est} to} $h_{k-1}$ \rev{and} $f_k$ \rev{at each step to obtain}
\[ |\mu(\widetilde\varphi_{h_k}) -\mu(\widetilde\varphi_{h_{k-1}})| \le\frac{2D(\mu)}t+\frac t2 e^{2tP}P. \]
\rev{Summing from} $k=2$ \rev{to} $N$ \rev{and using the triangle inequality to compare the two endpoint values gives}
\begin{equation} \label{eq:telescoping-cover}
1\le(N-1)\left(\frac{2D(\mu)}t+\frac t2e^{2tP}P\right) \qquad(t>0).
\end{equation}
If $P=0$, letting $t\to\infty$ contradicts this inequality. Thus $P>0$, and we may choose $t=(2(N-1)P)^{-1}$. Substitution gives
\[ 1\le4(N-1)^2D(\mu)P+\frac14e^{1/(N-1)} \le4(N-1)^2D(\mu)P+\frac34, \]
since $N\ge2$ and $e<3$. Hence, rearrangement gives the required bound. \end{proof}

\begin{proof}[Proof of Theorem \ref{prop-4}] Let $\vec f$ be any smooth partition of unity subordinate to $\mathcal U$. Since each $f_i$ is supported in the displaceable set $U_i$, the vanishing property gives $\mu(\widetilde\varphi_{f_i})=0$. Apply Lemma~\ref{lem:vanishing-partition-bound} and take the infimum over all subordinate partitions. \end{proof}

\section{Proof of Theorem \ref{thm:spectral-cover}}\label{sec:proof-F}

Fix the contact spectral invariant $c_\alpha$ supplied by \cite[Theorem 1.10]{CKZ}. \rev{To obtain the cover bound, we first derive control by contact Hamiltonians and then prove a local spectral-energy inequality. We will combine these two estimates with the flow comparison in Lemma~\ref{lem:flow-comparison}.} \rev{The path length \keepmath{$\ell_{\alpha}$}, defined in \eqref{eq:contact-path-length}, also appears in the definition of \keepmath{$e_\alpha(U)$}.}

\begin{lemma}\label{lem:spectral-calibration} For every $t\in\R$, the class represented by the Reeb flow satisfies
\begin{equation} \label{eq:spectral-reeb-calibration}
c_\alpha(\widetilde\varphi_1^{\,t})=t.
\end{equation}
If $t\ge0$ and a based contact isotopy $\psi_s$, $0\le s\le t$, is generated by $h_s$ and represents $a\in\widetilde{\Cont}_0(Y,\xi)$, then
\begin{equation} \label{eq:spectral-length-bound}
\int_0^t\min_Y h_s\,ds\le c_\alpha(a)\le\int_0^t\max_Y h_s\,ds.
\end{equation}
In particular, both $|c_\alpha(a)|$ and $|c_\alpha(a^{-1})|$ are at most $\ell_\alpha(\psi)$. Moreover, the contact spectral invariant $c_\alpha$ is continuous along smooth families of contact isotopies.
\end{lemma}

\begin{proof}  We first establish a uniform lower bound for non-zero Reeb periods. Since $R_\alpha$ is nowhere zero and $Y$ is compact, there is $\tau>0$ such that $\varphi_1^s$ has no fixed point for $0<|s|<\tau$. Indeed, around each point choose a coordinate in which the component of $R_\alpha$ is strictly positive, and then choose a smaller neighborhood and a sufficiently short time interval on which every trajectory starting there remains in that coordinate neighborhood. The chosen coordinate is strictly increasing along each such trajectory. Hence, a finite covering of $Y$ gives a uniform time interval.

\rev{We now use this bound to calibrate the spectral invariant near time zero. Monotonicity and normalization give} $c_\alpha(\widetilde\varphi_1^{\,s})\ge0$ for $s\ge0$ and $c_\alpha(\widetilde\varphi_1^{\,s})\le0$ for $s\le0$. Continuity from above gives $c_\alpha(\widetilde\varphi_1^{\,s})\to0$ as $s \to 0^+$. \rev{To obtain the corresponding limit from below, subadditivity gives}
\[ 0=c_\alpha(\widetilde{\mathrm{id}})\le c_\alpha(\widetilde\varphi_1^{\,s})+c_\alpha(\widetilde\varphi_1^{\,-s}). \]
\rrev{Therefore}, the same convergence to zero holds as $s \to 0^-$. By spectrality, $c_\alpha(\widetilde\varphi_1^{\,s})-s$ is either zero or a nonzero signed period of the Reeb flow. Since its absolute value is less than $\tau$ for sufficiently small $|s|$, it must be zero. Thus \eqref{eq:spectral-reeb-calibration} holds near zero.

\rev{To extend this identity to arbitrary} $t\in\R$, \rev{choose a positive integer} $m$ for which both $t/m$ and $-t/m$ lie in this interval. Subadditivity and the group law for the lifted Reeb flow give
\[ c_\alpha(\widetilde\varphi_1^{\,t})\le m c_\alpha(\widetilde\varphi_1^{\,t/m})=t,\qquad c_\alpha(\widetilde\varphi_1^{\,-t})\le-t. \]
Applying subadditivity once more to their product, which is the identity, gives $c_\alpha(\widetilde\varphi_1^{\,t})\ge-c_\alpha(\widetilde\varphi_1^{\,-t})\ge t$. Therefore, (\ref{eq:spectral-reeb-calibration}) holds for all real times. \rev{With this calibration established, the argument following Remark \ref{rmk-normalization} gives the inequalities in (\ref{eq:spectral-length-bound}) directly from (\ref{eq:spectral-reeb-calibration}).}

\rev{The absolute bound for} $c_\alpha(a)$ \rev{now follows immediately. To obtain the same bound for the inverse, use the based path} $\widehat\psi_s=\psi_{t-s}\psi_t^{-1}$, $0\le s\le t$, which represents $a^{-1}$. Its Hamiltonian is $-h_{t-s}$, so it has exactly the same length as $\psi$. Hence, applying the preceding bound to this path proves $|c_\alpha(a^{-1})|\le\ell_\alpha(\psi)$. This choice of path avoids introducing the conformal factor of the pointwise inverse isotopy.

\rrev{It remains to prove continuity. If a path $\eta$ represents $a^{-1}b$, subadditivity gives}
{
\[ c_\alpha(b)-c_\alpha(a)\le c_\alpha(a^{-1}b)\le\ell_\alpha(\eta). \]
}
\rrev{The reversed, right-translated path represents $b^{-1}a$ and has the same length, so the reverse inequality also holds. Thus $|c_\alpha(a)-c_\alpha(b)|\le\ell_\alpha(\eta)$. For a smooth family of based isotopies $\psi_s^u$, generated by $h_s^u$ for $0\le s\le t$, fix $u_0$ and take $\eta_s^u=(\psi_s^{u_0})^{-1}\psi_s^u$. The time-dependent identity \eqref{eq:relative-flow-hamiltonian} yields}
{
\[ \ell_\alpha(\eta^u)\le\int_0^t e^{\norm{\kappa(\psi_s^{u_0})}_\infty}\norm{h_s^u-h_s^{u_0}}_\infty\,ds\longrightarrow0\qquad(u\to u_0). \]
}
\rrev{This proves the asserted continuity.} \end{proof}

Recall that a contactomorphism $\varphi$ has a discriminant point if $\varphi(x)=x$ and $\kappa(\varphi)(x)=0$ for some $x$. \rev{For the local spectral-energy estimate, we need a displacing isotopy whose translated spectrum omits a prescribed number. The next lemma constructs such an isotopy while preserving displacement of a fixed compact set.}

\begin{lemma}\label{lem:spectral-displacer} Let $K\subset U$ be compact, where $U$ is open, and let $\psi_s^0$, $0\le s\le1$, satisfy $\psi_1^0(U)\cap U=\varnothing$. If $L>\ell_\alpha(\psi^0)$, there is an arbitrarily small smooth perturbation of $\psi_s^0$ through based contact isotopies such that
\[ \ell_\alpha(\psi)<L,\qquad \psi_1(K)\cap U=\varnothing,\qquad L\notin\operatorname{Spec}_\alpha(\psi_1). \]
\end{lemma}

\begin{proof} Contactomorphisms without discriminant points are dense in the smooth topology. \rev{This density follows from a transversality argument, which is also contained in the parametric statement \cite[Proposition 4.15]{CKZ}.} If $d=\dim Y$ and $\varphi$ is a contactomorphism, consider
\[ x\longmapsto(x,\varphi(x),\kappa(\varphi)(x))\in Y\times Y\times\R. \]
A discriminant point corresponds to an intersection of this map with $\{(y,y,0) \,| \, y\in Y\}$, which has codimension $d+1$. Postcomposing $\varphi$ with the short contact flow of a function $k$ changes the last two entries, to first order at $y=\varphi(x)$, by $(X_k(y),R_\alpha k(y))$. \rev{Here, the value} $k(y)$, \rev{the restriction of} $dk(y)$ \rev{to} $\xi_y$, \rev{and} $R_\alpha k(y)$ \rev{can be prescribed independently. \rrev{Therefore}, the contact Hamiltonian equations show that these data prescribe an arbitrary element of} $T_yY\times\R$. Using cutoffs and a finite covering of the compact manifold, one obtains a finite-dimensional family of perturbations for which the total evaluation map is transverse to the indicated submanifold. Hence, parametric transversality gives arbitrarily small perturbations whose evaluation maps are transverse. Since the domain has dimension $d$ and the required codimension is $d+1$, its image is disjoint from the indicated submanifold. Therefore, such contactomorphisms are dense.

\rev{We apply this density statement to} $\varphi_1^{-L}\psi_1^0$. Since the Reeb flow preserves $\alpha$, the condition that $\varphi_1^{-L}\psi_1$ have no discriminant point is exactly $L\notin\operatorname{Spec}_\alpha(\psi_1)$. \rev{To realize the perturbation at the level of paths, postcompose the original based isotopy with an arbitrarily small based contact isotopy. The Hamiltonian depends smoothly on this perturbation, so the length remains less than} $L$. Moreover, $\psi_1^0(K)$ is a compact subset of the open set $\psi_1^0(U)$, which is disjoint from $U$. Thus the perturbed endpoint still sends $K$ outside $U$. \end{proof}

\begin{lemma}\label{lem:spectral-barrier} Let $g\ge0$ have compact support in a contact displaceable open set $U$. If $L>e_\alpha(U)$ and
\[ \varphi_1^L(\supp g)\subset U, \]
then $0\le c_\alpha(\widetilde\varphi_g^{\,t})<2L$ for every $t\ge0$.
\end{lemma}

\begin{proof} Put $K:=\supp g$. Choose an isotopy displacing $U$ with length less than $L$, and apply Lemma~\ref{lem:spectral-displacer} to obtain $\psi_s$ such that $\ell_\alpha(\psi)<L$, $\psi_1(K)\cap U=\varnothing$, and $L\notin\operatorname{Spec}_\alpha(\psi_1)$. \rev{We first use these properties to exclude the prescribed translation time. The flow of} $g$ \rev{preserves} $K$: it is the identity outside $K$, and each of its time maps is a diffeomorphism. For $x\in K$ and $t\ge0$, we therefore have
\[ \psi_1\varphi_g^t(x)\notin U,\qquad \varphi_1^L(x)\in U. \]
\rev{Thus these points cannot coincide. On the other hand, if} $x\notin K$, the flow $\varphi_g^t$ is the identity on a neighborhood of $x$ and has conformal factor zero there. A translated point of $\psi_1\varphi_g^t$ with translation time $L$ would then be a translated point of $\psi_1$ with the same translation time, which was excluded. \rev{The two cases together give}
\begin{equation} \label{eq:spectral-barrier-gap}
L\notin\operatorname{Spec}_\alpha(\psi_1\varphi_g^t)\qquad(t\ge0).
\end{equation}

\rev{We now turn this spectral exclusion into an upper bound. Let} $\widetilde\psi$ \rev{denote the lift represented by} $\psi_s$. \rev{Lemma~\ref{lem:spectral-calibration} shows that} $t\mapsto c_\alpha(\widetilde\psi\,\widetilde\varphi_g^{\,t})$ is continuous. At $t=0$, its value is at most $\ell_\alpha(\psi)<L$. Since spectrality and \eqref{eq:spectral-barrier-gap} prevent this value from crossing $L$, it remains strictly less than $L$ for every $t\ge0$. By subadditivity and the inverse bound in Lemma~\ref{lem:spectral-calibration},
\[ c_\alpha(\widetilde\varphi_g^{\,t})\le c_\alpha(\widetilde\psi^{-1})+c_\alpha(\widetilde\psi\,\widetilde\varphi_g^{\,t})<\ell_\alpha(\psi)+L<2L. \]
The lower bound follows from $g\ge0$, monotonicity, and normalization. \end{proof}

\begin{remark}[\rev{Communication with Dylan Cant}] \rev{The proof of Lemma \ref{lem:spectral-barrier} uses the same spectral argument as the proof of Theorem 8 in \cite[Section 3.5.4]{Dylan-eternal}. This argument goes back to \cite{Sch00,Ost-comparion}. In fact, Lemma \ref{lem:spectral-barrier} above also follows from Theorem 8 in \cite{Dylan-eternal}.}  
 \end{remark}

\begin{prop}[Local spectral-energy inequality]\label{prop:local-spectral-energy} Let $U$ be a contact displaceable open set, and let $h\ge0$ be smooth with $\supp h\subset U$. Then, for every $t\ge0$,
\begin{equation} \label{eq:local-spectral-energy}
0\le c_\alpha(\widetilde\varphi_h^{\,t})\le e_\alpha(U)\bigl(1+e^{t\norm{R_\alpha h}_\infty}\bigr).
\end{equation}
\end{prop}

\begin{proof} Write $q=\norm{R_\alpha h}_\infty$ and fix $L>e_\alpha(U)$ and $\epsilon>0$. \rev{To apply the preceding lemma, we cut off the Hamiltonian so that the required Reeb translate of its support remains in the displaceable set. Choose a smooth function} $\chi:\R\to\R$ \rev{such that}
\[ \chi(r)=0\ \text{for }r\le Lq+\epsilon/2,\qquad 0\le\chi'(r)\le1,\qquad \chi'(r)=1\ \text{for }r\ge Lq+\epsilon. \]
For example, integrate a smooth nondecreasing transition from zero to one on the indicated interval. For $r\ge0$, these conditions imply $0\le\chi(r)\le r$ and $0\le r-\chi(r)\le Lq+\epsilon$. Set $g=\chi\circ h$. Then
\begin{equation} \label{eq:spectral-cutoff-error}
0\le g\le h,\qquad \norm{h-g}_\infty\le Lq+\epsilon,\qquad \norm{R_\alpha g}_\infty\le q,
\end{equation}
and $\supp g\subset\{h\ge Lq+\mrev{\epsilon}/2\}\subset U$. For every $x\in\supp g$, integration along the Reeb flow gives
\[ h(\varphi_1^L(x))=h(x)+\int_0^L R_\alpha h(\varphi_1^s(x))\,ds\ge h(x)-Lq\ge\epsilon/2>0. \]
Thus $\varphi_1^L(\supp g)\subset\{h>0\}\subset U$. Therefore, Lemma~\ref{lem:spectral-barrier} applies and gives $c_\alpha(\widetilde\varphi_g^{\,t})<2L$ for all $t\ge0$.

\rev{The cutoff has the required spectral bound. It remains to estimate the error in replacing} $g$ \rev{by} $h$. \rev{For this purpose, put} $\eta_s=\varphi_g^{-s}\varphi_h^s$ for $0\le s\le t$. By \eqref{eq:relative-flow-length} and \eqref{eq:spectral-cutoff-error},
\[ \ell_\alpha(\eta)\le\norm{h-g}_\infty\int_0^t e^{s\norm{R_\alpha g}_\infty}\,ds\le(Lq+\epsilon)\int_0^t e^{sq}\,ds. \]
The identity $\varphi_h^s=\varphi_g^s\eta_s$ gives the corresponding factorization in the universal cover. Subadditivity and Lemma~\ref{lem:spectral-calibration} therefore give
\begin{equation} \label{eq:spectral-cutoff-comparison}
c_\alpha(\widetilde\varphi_h^{\,t})\le c_\alpha(\widetilde\varphi_g^{\,t})+c_\alpha([\eta])\le2L+(Lq+\epsilon)\int_0^t e^{sq}\,ds.
\end{equation}
If $q>0$, the last expression equals
\[ L(1+e^{tq})+\epsilon\frac{e^{tq}-1}{q}. \]
If $q=0$, it equals $2L+\epsilon t$. \rev{Thus letting} $\epsilon \to 0$ \rev{and then} $L \to e_\alpha(U)$ \rev{proves the upper bound in both cases. Although the auxiliary cutoff and displacer may depend on} $L$ \rev{and} $\mrev{\epsilon}$, \rev{the left-hand side of \eqref{eq:spectral-cutoff-comparison} is independent of those choices.} Finally, $h\ge0$ gives the lower bound by monotonicity. \end{proof}

\begin{cor}\label{cor:spectral-flow-comparison} For $f,g\in C^\infty(Y)$ and $t\ge0$,
\begin{equation} \label{eq:spectral-flow-comparison}
\begin{split}
c_\alpha(\widetilde\varphi_{f+g}^{\,t})&\le c_\alpha(\widetilde\varphi_f^{\,t})+c_\alpha(\widetilde\varphi_g^{\,t}) +\frac{t^2}{2}e^{t(\norm{R_\alpha f}_\infty+\norm{R_\alpha g}_\infty)}\norm{\{f,g\}_\alpha}_\infty.
\end{split}
\end{equation}
\end{cor}

\begin{proof} Let $\eta$ be the path in Lemma~\ref{lem:flow-comparison}. Its lifted factorization, subadditivity of $c_\alpha$, and \eqref{eq:spectral-length-bound} give
\[ c_\alpha(\widetilde\varphi_{f+g}^{\,t})\le c_\alpha(\widetilde\varphi_f^{\,t})+c_\alpha(\widetilde\varphi_g^{\,t})+c_\alpha([\eta])\le c_\alpha(\widetilde\varphi_f^{\,t})+c_\alpha(\widetilde\varphi_g^{\,t})+\ell_\alpha(\eta). \]
Therefore, substituting \eqref{eq:sum-flow-length} proves \eqref{eq:spectral-flow-comparison}. \end{proof}

\rev{We now combine the local spectral-energy estimate with the flow comparison to prove Theorem \ref{thm:spectral-cover}. The argument follows the strategy of successive partial sums used for Theorem \ref{prop-4}, specifically in Lemma \ref{lem:vanishing-partition-bound}.} 

\begin{proof}[Proof of Theorem \ref{thm:spectral-cover}] Let $\vec f=(f_1,\ldots,f_N)$ be any smooth partition of unity subordinate to $\mathcal U$. For this proof, write
\[ P=\Palpha(\vec f),\qquad E=\sum_{i=1}^N e_\alpha(U_i),\qquad h_k=\sum_{i=1}^k f_i. \]
The number $E$ is finite because each $U_i$ is displaceable, and it is \rev{positive because the cover has a nonempty member.} The partial sums and partition members satisfy \eqref{eq:partition-controls}, whose proof used only the definition of $P$. Moreover, \eqref{eq:partition-reeb-control} gives $\norm{R_\alpha f_i}_\infty\le P$ for every $i$, including $i=1$.

\rev{To compare the spectral invariant of the total sum with those of its summands, apply Corollary~\ref{cor:spectral-flow-comparison} to} $h_{k-1}$ \rev{and} $f_k$. \rev{Using \eqref{eq:partition-controls}, we obtain}
\[ c_\alpha(\widetilde\varphi_{h_k}^{\,t})\le c_\alpha(\widetilde\varphi_{h_{k-1}}^{\,t})+c_\alpha(\widetilde\varphi_{f_k}^{\,t})+\frac{t^2}{2}e^{2tP}P. \]
Sum these inequalities for $k=2,\ldots,N$. The intermediate terms cancel, $h_1=f_1$, and there are $N-1$ error terms. \rrev{Therefore},
\begin{equation} \label{eq:spectral-telescoping}
c_\alpha(\widetilde\varphi_{h_N}^{\,t})\le\sum_{i=1}^N c_\alpha(\widetilde\varphi_{f_i}^{\,t})+\frac{N-1}{2}t^2e^{2tP}P.
\end{equation}
\rev{The two sides of this inequality can now be controlled separately. Since} $h_N=1$, \rev{the left-hand side is} $t$ \rev{by \eqref{eq:spectral-reeb-calibration}. On the right, each \keepmath{$f_i$} is nonnegative and has compact support in \keepmath{$U_i$}.} Proposition~\ref{prop:local-spectral-energy} therefore gives
\[ c_\alpha(\widetilde\varphi_{f_i}^{\,t})\le e_\alpha(U_i)\bigl(1+e^{t\norm{R_\alpha f_i}_\infty}\bigr)\le e_\alpha(U_i)(1+e^{tP}). \]
Hence, substituting these bounds into \eqref{eq:spectral-telescoping} proves
\begin{equation} \label{eq:spectral-master}
t\le E(1+e^{tP})+\frac{N-1}{2}t^2e^{2tP}P\qquad(t\ge0).
\end{equation}

\rev{This estimate forces the desired lower bound. Indeed, suppose for contradiction that} $P<1/(16(N-1)E)$. \rev{Taking} $t=4E$ \rev{will make the right-hand side too small.} Then $tP<1/(4(N-1))\le1/4$, so $e^{tP}<2$ and $e^{2tP}<2$. Dividing \eqref{eq:spectral-master} by $t$ yields
\[ 1\le\frac{1+e^{tP}}4+\frac{N-1}{2}t e^{2tP}P<\frac34+4(N-1)EP<1. \]
This is impossible. Therefore, $P\ge1/(16(N-1)E)$ for every subordinate partition. Taking the infimum over all such partitions proves \eqref{eq:main-spectral-cover}. \end{proof}

\section{Proof of Theorem \ref{thm-qcbft}}\label{sec:proof-G}

\rev{We first prove the result when the image lies in the unit cube; the general statement will then follow by rescaling.} 

 \begin{prop} \label{prop-5}
Let $(Y, \xi = \ker \alpha)$ be a closed contact manifold that admits a non-zero homogeneous quasi-morphism $\mu$ satisfying monotonicity and the vanishing property, \rrev{normalized by \eqref{eq:mu-normalization}.} Let $F = (f_1, ..., f_k): Y \to [-1,1]^k$ be smooth and non-zero. \rev{Then there exists a constant \keepmath{$c_k>0$}, depending only on \keepmath{$k$}, such that, for any \keepmath{$0< \delta \leq 1$}, if} 
 \begin{equation} \label{quantify-poisson}
D(\mu) \max\left\{\max_{i=1}^k \rrev{\norm{R_\alpha f_i}_\infty}, \, \max_{i,j=1}^k \rrev{\norm{\{f_i, f_j\}_{\alpha}}_\infty} \right\}  < c_k \delta^{2k+2}  
 \end{equation}
 then there is a point $\vec{a} \in [-1,1]^k$ such that 
 \postdisplaypenalty=10000
 \[ K_{\vec{a}}(\delta): = \left\{ x \in Y \,| \, |F(x) - \vec{a}| \leq \delta\right\} \]
 \rev{is contact non-displaceable in \keepmath{$(Y, \xi)$}.} 
 \end{prop}
 
 \rev{The restriction of the target of \keepmath{$F$} from \keepmath{$\R^k$} to \keepmath{$[-1,1]^k$} is motivated by the definition of \keepmath{$\rrev{\Palpha}(\vec{f})$} in (\ref{dfn-Pf}), where the coefficient vectors \keepmath{$a,b$} lie in \keepmath{$[-1,1]^N$}, and by the estimate for \keepmath{$pb_{\alpha}(\mathcal U)$} in Theorem~\ref{prop-4}.} 

\begin{proof}[Proof of Proposition \ref{prop-5}] \rrev{We argue by contradiction. Suppose that \eqref{quantify-poisson} holds for some $0<\delta_*\leq1$, but $K_{\vec a}(\delta_*)$ is contact displaceable for every $\vec a\in[-1,1]^k$. We will specify $c_k$ near the end of the proof.} 

\rev{To apply the cover estimate, we first construct a finite open cover \keepmath{$\mathcal U$} of \keepmath{$Y$} by contact displaceable sets. For each lattice point \keepmath{$\nu \in \Z^k$}, consider the open cube} 
\[ Q_\nu: =  \delta_* \nu + \left(-\frac{3}{4}\delta_*,  \frac{3}{4}\delta_*\right)^k \subset \R^k. \]
\rev{The constant \keepmath{$\frac{3}{4}$} ensures that the cubes \keepmath{$Q_{\nu}$} cover the target \keepmath{$[-1,1]^k$}. We claim that the nonempty preimages \keepmath{$U_{\nu} : = F^{-1}(Q_{\nu})$} form a finite open cover of \keepmath{$Y$}. Indeed, for any \keepmath{$x \in Y$}, there exists a lattice point \keepmath{$\nu \in \Z^k$} such that} 
\[ \left| \frac{F(x)}{\delta_*} - \nu \right| \leq \frac{1}{2} \,\,\,\,\mbox{which implies} \,\,\,\, |F(x) - \delta_*\nu| \leq \frac{\delta_*}{2} < \frac{3}{4} \delta_*.\]
 \rev{Therefore, \keepmath{$F(x) \in Q_\nu$}, which proves the covering property. To show that only finitely many preimages are nonempty, observe that if \keepmath{$U_{\nu} \neq \emptyset$}, then the bounds \keepmath{$\rrev{\norm{f_i}_\infty} \leq 1$} imply \keepmath{$|\delta_*\nu| < 1+ \frac{3\delta_*}{4}$}, and hence} 
 \[ |\nu| < \frac{1}{\mrev{\delta_*}} + \frac{3}{4}. \]
 \rev{A \keepmath{$k$}-dimensional cube with side length} $\frac{2}{\mrev{\delta_*}} + \frac{3}{2}$ \rev{contains at most} $\left(\frac{2}{\mrev{\delta_*}} + \mrev{\frac{5}{2}}\right)^k$ \rev{lattice points. Moreover,} 
 \begin{equation} \label{est-open-cover}
 \left(\frac{2}{\mrev{\delta_*}} + \frac{5}{2}\right)^k = \frac{1}{\delta_*^k} \left(2 + \frac{5\mrev{\delta_*}}{2}\right)^k \leq \left(\frac{9}{2}\right)^k \frac{1}{\delta^k_*}.
 \end{equation}
\rev{Therefore, only finitely many sets \keepmath{$U_{\nu}$} are nonempty.} \rrev{Their lattice centers need not lie in $[-1,1]^k$. For each such index, project the center coordinatewise onto this cube by setting}
\[ \rrev{a_{\nu,i}:=\max\{-1,\min\{1,\delta_*\nu_i\}\},\qquad \vec a_\nu=(a_{\nu,1},\ldots,a_{\nu,k}).} \]
\rrev{Since $F(x)\in[-1,1]^k$, this projection cannot increase the center's distance from $F(x)$. Hence, for $x\in U_\nu$,}
\[ \rrev{|F(x)-\vec a_\nu|\le |F(x)-\delta_*\nu|<\tfrac34\delta_*<\delta_*.} \]
\rrev{Thus $U_\nu\subset K_{\vec a_\nu}(\delta_*)$, and $U_\nu$ is contact displaceable by assumption.} 

\rev{Having constructed the cover, we next construct a partition of unity \keepmath{$\vec{h} = \{h_\nu\}_{\nu \in \Z^k}$} subordinate to the open cover \keepmath{$\mathcal U = \{U_\nu\}_{\nu \in \Z^k}$}. We begin on the target cubes \keepmath{$Q_{\nu}$}.} Choose (and fix) a function $\eta: \R \to [0,1]$, compactly supported in $(-\frac{3}{4}, \frac{3}{4})$, such that $\sum_{n\in \Z} \eta(t - n) = 1$. Then for any lattice point $\nu = (\nu_1, ..., \nu_k) \in \Z^k$, consider the function 
\[ \chi_{\nu}: \R^k \to [0,1] \,\,\,\,\mbox{defined by} \,\,\,\,\chi_{\nu}(z) := \prod_{j=1}^k\eta\left(\frac{z_j}{\delta_*} - \nu_j\right). \]
 Each such function is compactly supported in $Q_{\nu}$, and together they satisfy $\sum_{\nu} \chi_{\nu} = 1$. \rev{We pull this partition back by defining} $h_{\nu}: = F^*\chi_\nu = \chi_\nu \circ F$. \rev{By construction,} $h_{\nu}$ \rev{is compactly supported in} $U_{\nu}$, \rev{and, for any} $x \in Y$, \rev{we have} 
 \[ \sum_{\nu \in \Z^k} h_\nu(x) = \sum_{\nu \in \Z^k} \chi_{\nu}(F(x)) = 1. \]
\rev{To compare this partition with the cover lower bound, we now estimate \keepmath{$\rrev{\Palpha}(\vec{h})$}.} 

\rev{By definition (\ref{dfn-Pf}), the brackets \keepmath{$\left\{\sum_{\nu} a_\nu h_\nu, \sum_\nu b_\nu h_\nu\right\}_{\alpha}$} appearing in \keepmath{$\rrev{\Palpha}(\vec{h})$} involve derivatives of \keepmath{$\chi_{\nu}$}, and hence of \keepmath{$\eta$}. Write \keepmath{$L: = \sup_{t \in \R} \sum_{n \in \Z} |\eta'(t-n)|$} and express the two linear combinations as} 
\[ \sum_{\nu} a_\nu h_\nu = A \circ F \,\,\,\,\mbox{and}\,\,\,\, \sum_{\nu} b_\nu h_\nu = B \circ F\]
where $A = \sum_{\nu} a_{\mrev{\nu}} \chi_{\mrev{\nu}}$ and $B = \sum_{\nu} b_{\mrev{\nu}} \chi_{\mrev{\nu}}$. \rev{The contact chain rule gives}
{
\begin{equation}\label{eq:contact-chain-rule}
\keepmath{\begin{aligned}
\{A\circ F,B\circ F\}_{\alpha}
&=
\sum_{i,j=1}^{k} A_i(F)B_j(F)
\Bigl(
\{f_i,f_j\}_{\alpha}
-f_i\rrev{R_\alpha f_j}
+f_j\rrev{R_\alpha f_i}
\Bigr)
\\
&\quad
+A(F)\sum_{j=1}^{k}B_j(F)\rrev{R_\alpha f_j} -B(F)\sum_{i=1}^{k}A_i(F)\rrev{R_\alpha f_i}
\end{aligned}}
\end{equation}
}
where $A_i$ denotes $\partial_i A$ and $B_j$ denotes $\partial_j B$. \rev{The partition construction also gives}
\[ |A(z)| \leq 1\,\,\,\,\mbox{and}\,\,\,\, \rrev{\norm{A_i}_\infty} = \rrev{\norm{\partial_i A}_\infty} \leq \frac{L}{\delta_*}.\]
\rev{Combining these bounds with the chain rule, we obtain}
\begin{align*}
\rrev{\norm{\left\{\sum_{\nu} a_\nu h_\nu, \sum_\nu b_\nu h_\nu\right\}_{\alpha}}_\infty} & = \rrev{\norm{\{A\circ F,B\circ F\}_{\alpha}}_\infty} \\
& \leq \left(\frac{3k^2L^2}{\delta_*^2} + \frac{2kL}{\delta_*} \right)\max\left\{\max_{i=1}^k \rrev{\norm{R_\alpha f_i}_\infty}, \, \max_{i,j=1}^k \rrev{\norm{\{f_i, f_j\}_{\alpha}}_\infty} \right\} \\
& \leq \frac{C_k}{\delta_*^2} \max\left\{\max_{i=1}^k \rrev{\norm{R_\alpha f_i}_\infty}, \, \max_{i,j=1}^k \rrev{\norm{\{f_i, f_j\}_{\alpha}}_\infty} \right\} 
\end{align*}
\rev{where \keepmath{$C_k : = 3k^2 L^2 + 2kL$} depends only on \keepmath{$k$}.} The last inequality holds since $\delta_* \leq 1$ (so $\frac{2kL}{\delta_*} \leq \frac{2kL}{\delta_*^2}$). \rev{The estimate is uniform in the coefficient vectors, so}
\begin{equation} \label{ph-big max}
\rrev{\Palpha}(\vec{h}) \leq \frac{C_k}{\delta_*^2} \max\left\{\max_{i=1}^k \rrev{\norm{R_\alpha f_i}_\infty}, \, \max_{i,j=1}^k \rrev{\norm{\{f_i, f_j\}_{\alpha}}_\infty} \right\}.
\end{equation}

\rev{Finally, applying Theorem~\ref{prop-4} to} $\mathcal U = \{U_{\nu}\}_{\nu \in \Z^k}$ and $\vec{h} = \{h_{\nu}\}_{\nu \in \Z^k}$, we have 
\[ \frac{1}{16(N-1)^2 D(\mu)} \leq pb_{\alpha}(\mathcal U) \leq \rrev{\Palpha}(\vec{h}) \]
\rev{where \keepmath{$N$} is the number of indices \keepmath{$\nu$} for which \keepmath{$U_\nu \neq \emptyset$}.} \rrev{Here $N\ge2$, since a single set covering $Y$ would be $Y$ itself and could not be displaced.} \rev{Combining the cardinality bound (\ref{est-open-cover}) with the bracket estimate (\ref{ph-big max}), we obtain} 
\begin{equation} \label{ck-inequ}
 \frac{\delta_*^{2k}}{16 \left(\frac{9}{2}\right)^{2k} D(\mu)} \leq \frac{C_k}{\delta_*^2} \max\left\{\max_{i=1}^k \rrev{\norm{R_\alpha f_i}_\infty}, \, \max_{i,j=1}^k \rrev{\norm{\{f_i, f_j\}_{\alpha}}_\infty} \right\}. 
 \end{equation}
\rev{Set \keepmath{$c_k: = \left(32 \left(\frac{9}{2}\right)^{2k} C_k \right)^{-1}$}. Then (\ref{ck-inequ}) and the quantitative condition (\ref{quantify-poisson}) imply} 
\[\frac{1}{16 \left(\frac{9}{2}\right)^{2k} C_k} \leq c_k = \frac{1}{32 \left(\frac{9}{2}\right)^{2k} C_k} \]
\rev{This contradicts the choice of the constant and completes the proof.} \end{proof} 

\rev{We now prove Theorem~\ref{thm-qcbft} by rescaling.} 

\begin{proof}[Proof of Theorem \ref{thm-qcbft}] Without loss of generality, assume $\max_{i=1}^k \rrev{\norm{f_i}_\infty}>0$. Then fix $0<\delta \mrev{\leq} \max_{i=1}^k \rrev{\norm{f_i}_\infty}$. Consider the rescaling 
\[ \tilde{F} = \left( \frac{f_1}{\max_{i=1}^k \rrev{\norm{f_i}_\infty}}, \mrev{\ldots,} \frac{f_k}{\max_{i=1}^k \rrev{\norm{f_i}_\infty}} \right) = :(\tilde{f}_1, ..., \tilde{f}_k) \]
\rev{so that} $\tilde{F}: Y \to [-1,1]^k$. \rev{To check the quantitative hypothesis after rescaling, we compute} 
\[ \rrev{R_\alpha \tilde{f}_i} = \frac{\rrev{R_\alpha f_i}}{\max_{i=1}^k \rrev{\norm{f_i}_\infty}} \,\,\,\,\Longrightarrow \,\,\,\, \rrev{\norm{R_\alpha \tilde{f}_i}_\infty} = \frac{\rrev{\norm{R_\alpha f_i}_\infty}}{\max_{i=1}^k \rrev{\norm{f_i}_\infty}}\]
 and 
 \[ \{\tilde{f}_i, \tilde{f}_j\}_{\alpha} = \frac{\{f_i, f_j\}_{\alpha}}{\left(\max_{i=1}^k \rrev{\norm{f_i}_\infty}\right)^2}\,\,\,\,\Longrightarrow \,\,\,\, \rrev{\norm{\{\tilde{f}_i, \tilde{f}_j\}_{\alpha}}_\infty} = \frac{\rrev{\norm{\{f_i, f_j\}_{\alpha}}_\infty}}{\left(\max_{i=1}^k \rrev{\norm{f_i}_\infty}\right)^2}.\]
\rev{The quantitative condition in Theorem~\ref{thm-qcbft} therefore takes the following form in terms of the functions \keepmath{$\tilde{f}_i$}:} 
\begin{align*}
D(\mu) \max\left\{\max_{i=1}^k \rrev{\norm{R_\alpha \tilde{f}_i}_\infty}, \, \max_{i,j=1}^k \rrev{\norm{\{\tilde{f}_i, \tilde{f}_j\}_{\alpha}}_\infty} \right\}  < c_k \left(\frac{\delta}{\max_{i=1}^k \rrev{\norm{f_i}_\infty}} \right)^{2k+2}.  
 \end{align*}
 \rev{Setting} $\tilde{\delta}: = \frac{\delta}{\max_{i=1}^k \rrev{\norm{f_i}_\infty}}$ \rev{gives} $\tilde{\delta} \in \mrev{(0,1]}$.  \rev{Thus the rescaled system satisfies the hypotheses of Proposition \ref{prop-5}.} \rrev{If $\tilde a$ is the center supplied by that proposition, multiply it by $M=\max_i\norm{f_i}_\infty$. The resulting set is exactly $\{x:|F(x)-M\tilde a|\le\delta\}$, which proves the theorem.} 
\end{proof}

\begin{proof}[Proof of Corollary \ref{cor:big-fiber}] \rrev{If $F\equiv0$, its zero fiber is $Y$ and is non-displaceable. We may therefore assume $\max_i\norm{f_i}_\infty>0$.} For the given map $F = (f_1, ..., f_k): Y \to \R^k$, \rev{the hypotheses give} 
\[ \max_{i=1}^k \rrev{\norm{R_\alpha f_i}_\infty} = \max_{i, j=1}^k\rrev{\norm{\{f_i, f_j\}_{\alpha}}_\infty} =0. \]
\rev{Thus the quantitative condition in Theorem~\ref{thm-qcbft} holds for any sequence \keepmath{$\delta_n \to 0$} of admissible positive radii. For each \keepmath{$\delta_n$}, there exists \keepmath{$\vec{a}_n \in \R^k$} such that} 
\begin{equation} \label{k-an}
K_{\vec{a}_n}(\delta_n) = \{x \in Y \,| \, |F(x) - \vec{a}_n| \leq \delta_n\} 
\end{equation}
\rev{is contact non-displaceable and, in particular, nonempty. To pass to a limiting fiber, we first observe that the points \keepmath{$\mrev{\vec{a}_n}$} lie in a bounded region of \keepmath{$\R^k$}.} Indeed, for any $x_n \in K_{\vec{a}_n}(\delta_n)$, by the condition in (\ref{k-an}), we have 
\[ |\vec{a}_n| \leq |F(x_n)| + \delta_n \leq \max_Y|F| + \delta_n \]
\rev{where \keepmath{$\max_Y|F|$} is finite by compactness of \keepmath{$Y$}.} In particular, $|\vec{a}_n| < 2 \max_Y|F|$ when $n$ is sufficiently large. \rev{After passing to a subsequence, we may assume that \keepmath{$\vec{a}_n$} converges to some \keepmath{$\vec{a} \in \R^k$}. By compactness of \keepmath{$Y$}, a further subsequence satisfies \keepmath{$\mrev{x_n \to x}$} for some \keepmath{$x \in Y$}. These convergences give} 
\begin{equation} \label{Fx-a}
|F(x) -\vec{a}| \leq |F(x) - F(x_n)| + |F(x_n) - \vec{a}_n| + |\vec{a}_n - \vec{a}|  \to 0 
\end{equation}
which implies that $F(x) =\vec{a}$. Therefore, $\vec{a}$ lies in the image of $F$. \rev{Finally, we show that the fiber \keepmath{$F^{-1}(\vec{a})$} is contact non-displaceable. Otherwise, compactness of \keepmath{$F^{-1}(\vec{a})$} would give a displaceable neighborhood \keepmath{$U$} of this fiber.} \rrev{For all sufficiently large $n$, we have $K_{\vec a_n}(\delta_n)\subset U$. Otherwise, choose a subsequence $y_n\in K_{\vec a_n}(\delta_n)\setminus U$ converging to $y\in Y\setminus U$. The estimate}
\[ \rrev{|F(y_n)-\vec a|\le\delta_n+|\vec a_n-\vec a|\longrightarrow0} \]
\rrev{would put $y$ in $F^{-1}(\vec a)\subset U$, a contradiction.} \rev{Then \keepmath{$K_{\vec{a}_n}(\delta_n)$} would also be contact displaceable, a contradiction.} 
\end{proof}

\section{Proof of Theorem \ref{thm:subset-bounds}}\label{sec:proof-H}

\rev{We first establish the algebraic properties of \keepmath{$pb_{3, \alpha}$} and \keepmath{$pb_{4, \alpha}$} and their relation to the rigidity results. These properties will allow us to compare the bracket maxima in their definitions with \keepmath{$\Palpha(\cdot)$} and apply the lower bounds established earlier. Spectral invariants enter only in the proof of Corollary~\ref{cor:subset-energy-bounds}.}

\begin{lemma}\label{lem:subset-reeb-control} For smooth $f,g:Y\to[0,1]$, the complementary brackets in \eqref{eq:subset-square-definition} satisfy
\begin{equation} \label{eq:subset-four-brackets}
\begin{aligned}
\{1-f,g\}_\alpha&=R_\alpha g-\{f,g\}_\alpha,\\
\{f,1-g\}_\alpha&=-R_\alpha f-\{f,g\}_\alpha,\\
\{1-f,1-g\}_\alpha&=\{f,g\}_\alpha+R_\alpha f-R_\alpha g.
\end{aligned}
\end{equation}
\rrev{Therefore},
\begin{equation} \label{eq:subset-square-reeb}
\begin{aligned}
\norm{\{f,g\}_\alpha}_\infty&\le\max_{\substack{u\in\{f,1-f\}\\v\in\{g,1-g\}}}\norm{\{u,v\}_\alpha}_\infty,\\
\max\{\norm{R_\alpha f}_\infty,\norm{R_\alpha g}_\infty\}&\le2\max_{\substack{u\in\{f,1-f\}\\v\in\{g,1-g\}}}\norm{\{u,v\}_\alpha}_\infty.
\end{aligned}
\end{equation}
For a three-term partition $\vec h=(h_0,h_1,h_2)$,
\begin{equation} \label{eq:subset-triple-reeb}
R_\alpha h_i=\sum_{j=0}^2\{h_j,h_i\}_\alpha,\qquad \norm{R_\alpha h_i}_\infty\le2\max_{0\le j<k\le2}\norm{\{h_j,h_k\}_\alpha}_\infty.
\end{equation}
\end{lemma}

\begin{proof} Bilinearity, skew-symmetry, and $\{1,u\}_\alpha=R_\alpha u$ give \eqref{eq:subset-four-brackets}. In particular,
\[ R_\alpha g=\{f,g\}_\alpha+\{1-f,g\}_\alpha,\qquad -R_\alpha f=\{f,g\}_\alpha+\{f,1-g\}_\alpha. \]
Hence, the triangle inequality proves \eqref{eq:subset-square-reeb}. \rev{For the triple, the same mechanism applies: use} $\sum_j h_j=1$ to obtain the first identity in \eqref{eq:subset-triple-reeb}. Since its summand with $j=i$ vanishes, only two terms remain, which proves the desired estimate. \end{proof}

\begin{lemma}\label{prop:subset-basic} \rev{The contact Poisson bracket invariants} $pb_{3, \alpha}$ and $\mrev{pb_{4, \alpha}}$ defined in Definition \ref{def:contact-pb-subsets} satisfy the following symmetry identities:
\[ pb_{3, \alpha}(X_0, X_1, X_2) = pb_{3, \alpha}(X_{\sigma(0)}, X_{\sigma(1)}, X_{\sigma(2)}) \]
where $\mrev{(\sigma(0), \sigma(1), \sigma(2))}$ is any permutation of $\mrev{(0,1,2)}$, and 
\[ pb_{4, \alpha}(X_0, X_1, Y_0, Y_1) = pb_{4, \alpha}(X_1, X_0, Y_0, Y_1) = pb_{4, \alpha}(X_0, X_1, Y_1, Y_0) = pb_{4, \alpha}(Y_0, Y_1, X_0, X_1). \]
Both $pb_{3, \alpha}$ and $\mrev{pb_{4, \alpha}}$ are monotone under enlargement of the input sets and invariant under simultaneous application of a strict contactomorphism. Moreover, for any constant $a>0$,
\postdisplaypenalty=10000
\[ pb_{3,a\alpha}=a^{-1}pb_{3,\alpha},\qquad pb_{4,a\alpha}=a^{-1}pb_{4,\alpha}, \]
with the same input sets on both sides.
\end{lemma}

\begin{proof} Permuting the three sets permutes the three partition functions, and skew-symmetry preserves the maximum of their bracket norms. For the quadruple, exchanging $X_0,X_1$ corresponds to $f\mapsto1-f$, and exchanging $Y_0,Y_1$ to $g\mapsto1-g$. Each operation permutes the four brackets; exchanging the two ordered pairs of sets corresponds to $(f,g)\mapsto(g,f)$ \rev{and again preserves the bracket norm by skew-symmetry.} This proves the desired symmetry. \rev{Enlarging an input set restricts the class of admissible functions, which proves monotonicity.}

\rev{To prove invariance, suppose} $\psi^*\alpha=\alpha$. \rev{The contact Hamiltonian equations then imply}
\[ \{f\circ\psi^{-1},g\circ\psi^{-1}\}_\alpha=\{f,g\}_\alpha\circ\psi^{-1}. \]
Therefore, composing all input Hamiltonian functions with $\psi^{-1}$ preserves their ranges, their values near the corresponding sets, \rev{and the norms of their brackets.} Applying the same argument to $\psi^{-1}$ proves invariance. \rev{For the scaling relation, observe that} $R_{a\alpha}=a^{-1}R_\alpha$, and in general the contact Hamiltonian vector field of $f$ under $a \alpha$ is $a^{-1}X_f$, so $\{f,g\}_{a\alpha}=a^{-1}\{f,g\}_\alpha$. \rev{The admissibility of the input Hamiltonian functions does not depend on this constant rescaling, which proves the final assertion.} \end{proof}

\begin{prop}\label{prop:subset-lsc} Let $(Y, \ker \alpha)$ be a closed contact manifold. \rev{For three-term partitions of unity} $\vec h=(h_0,h_1,h_2)$ and smooth pairs $f,g:Y\to[0,1]$, the maps
\[
\begin{aligned}
\vec h&\longmapsto\max_{0\le i<j\le2}\norm{\{h_i,h_j\}_\alpha}_\infty \qquad 
(f,g)&\longmapsto\max_{\substack{u\in\{f,1-f\}\\v\in\{g,1-g\}}}\norm{\{u,v\}_\alpha}_\infty,
\end{aligned}
\]
\rev{are lower semicontinuous with respect to componentwise uniform convergence within the respective classes of admissible triples and pairs.} 
\end{prop}

\begin{proof} \rrev{For either functional, let $M_j$ denote its value at the $j$th term of a uniformly convergent sequence. If $\liminf_jM_j=+\infty$, there is nothing to prove. Otherwise, pass to a subsequence on which $M_j$ converges to its finite liminf. Lemma~\ref{lem:subset-reeb-control} bounds the Reeb derivatives of every input function uniformly on this subsequence; the same bound applies to the complementary functions in the four-set case. Thus condition~\eqref{Reeb-cond} holds for every bracket in the corresponding maximum. Theorem~\ref{prop-3} then bounds each limiting bracket norm by $\liminf_jM_j$. Taking the maximum proves the assertion.} \end{proof}

Recall that for a system of contact Hamiltonian functions $\vec{f} = (f_1, ..., f_N)$, we have defined $\Palpha(\vec{f})$ in (\ref{dfn-Pf}), and for an open cover $\mathcal U$ of $Y$, we have defined $pb_{\alpha}(\mathcal U)$ in (\ref{eq:contact-pb-cover}). 

\begin{prop}\label{prop:subset-cover-comparison} Let $(Y, \ker \alpha)$ be a closed contact manifold. Then for an admissible triple of compact sets $X_0, X_1, X_2$, we have 
\begin{equation} \label{eq:subset-triple-cover}
\frac14pb_\alpha(\{Y\setminus X_i\}_{i=0}^2)\le pb_{3,\alpha}(X_0,X_1,X_2)\le pb_\alpha(\{Y\setminus X_i\}_{i=0}^2).
\end{equation}
For an admissible quadruple, the sets $U_{ij}:=Y\setminus(X_i\cup Y_j)$ form an open cover, and
\begin{equation} \label{eq:subset-square-cover}
pb_{4,\alpha}(X_0,X_1,Y_0,Y_1)\ge\frac1{56}pb_\alpha(\{U_{ij}\}_{i,j=0}^1).
\end{equation}
\end{prop}

\begin{proof} \rev{We first compare the triple invariant with the cover invariant. A function on the compact manifold} $Y$ \rev{vanishes} on a neighborhood of $X_i$ precisely when its support is a compact subset of $Y\setminus X_i$. Thus the triples in \eqref{eq:subset-triple-definition} are exactly the partitions subordinate to the open cover $\{Y \backslash X_0, Y \backslash X_1, Y \backslash X_2\}$. Each individual bracket occurs among the combinations defining $\Palpha(\cdot)$, so
\[ \max_{0\le i<j\le2}\norm{\{h_i,h_j\}_\alpha}_\infty\le\Palpha(\vec h). \]
\rev{Taking infima gives} $pb_{3,\alpha}(X_0,X_1,X_2)\le pb_\alpha(\{Y\setminus X_i\}_{i=0}^2)$. 
\rev{For the reverse inequality, it suffices to bound arbitrary linear combinations in terms of the three pairwise brackets. Fix a point} $p \in Y$ and fix all coefficients in the two expressions $\sum_{i=0}^2 a_i h_i$ and $\sum_{i=0}^2 b_i h_i$ except $a_0 = t$. Since the bracket is bilinear, \rev{its value at \keepmath{$p$}} has the form $At + B$ where $A, B$ are real numbers independent of $t$. Then for $-1 \leq t \leq 1$, we have 
\[ At+B = \frac{1-t}{2}(B-A) + \frac{1+t}{2}(B+A) \] 
which implies that 
\[ |At+B| \leq \frac{1-t}{2}|B-A| + \frac{1+t}{2} |B+A| \leq \max\{|B-A|, |B+A|\}. \]
\rev{Thus at least one of the endpoints \keepmath{$t=-1$} or \keepmath{$t=1$} gives an absolute value no smaller than that at \keepmath{$t$}. Applying this argument successively to all six coefficients shows that the maximum is attained at corner vectors \keepmath{$a,b\in\{-1,1\}^3$}.} \rev{It therefore remains to estimate the bracket at these vectors, where it has the form}
\[ \sum_{0\le i<j\le2}(a_i b_j-a_j b_i)\{h_i,h_j\}_\alpha. \]
Since $a_i^2 = 1$, we have $a_i b_j - a_j b_i = a_i a_j (a_jb_j - a_ib_i)$. \rev{Since every product} $a_i b_i$ \rev{is equal to either} $1$ \rev{or} $-1$, \rev{this identity shows that} 
\[ |a_ib_j-a_jb_i|=\begin{cases}0,&a_i b_i=a_j b_j,\\2,&a_i b_i\ne a_j b_j.\end{cases} \]
Now, consider the three products $a_0b_0, a_1b_1, a_2b_2$. \rev{If they are all equal, all the coefficients vanish. Otherwise, two have the same sign and the third has the opposite sign. Thus only two coefficients can be nonzero, and we obtain} 
\begin{align*}
\left|\left\{ \sum_{i=0}^2 a_i h_i, \sum_{\mrev{j}=0}^2 b_j h_j\right\}_{\alpha}(p) \right| & \leq 2 \max_{0\le i<j\le2}\norm{\{h_i,h_j\}_\alpha}_\infty + 2\max_{0\le i<j\le2}\norm{\{h_i,h_j\}_\alpha}_\infty \\
& = 4\max_{0\le i<j\le2}\norm{\{h_i,h_j\}_\alpha}_\infty. 
\end{align*}
Since this holds for every point $p \in Y$, taking the uniform norm gives
\[ \Palpha(\vec h)\le4\max_{0\le i<j\le2}\norm{\{h_i,h_j\}_\alpha}_\infty. \]
Therefore, taking infima proves \eqref{eq:subset-triple-cover}.

\rev{For the quadruple, we use a product partition. First, every point lies outside at least one of} $X_0,X_1$ and outside at least one of $Y_0,Y_1$, so it belongs to some $U_{ij}$. Fix an admissible pair $(f,g)$ and put
\[ f_0=f,\quad f_1=1-f,\qquad g_0=g,\quad g_1=1-g,\qquad h_{ij}=f_ig_j. \]
The four functions $h_{ij}$ are non-negative and sum to $1$. Since $f_i$ vanishes near $X_i$ and $g_j$ vanishes near $Y_j$, their product has support compactly contained in $U_{ij}$. Thus they give a subordinate partition of unity. \rev{To estimate the bracket functional of this partition, fix coefficient vectors} $a,b\in[-1,1]^4$ \rev{and consider the bilinear polynomials}
\[ A(f,g)=\sum_{i,j=0}^1a_{ij}f_ig_j,\qquad B(f,g)=\sum_{i,j=0}^1b_{ij}f_ig_j\mrev{.} \]
\rev{Thus \keepmath{$A(f,g)$} and \keepmath{$B(f,g)$} are linear combinations of \keepmath{$fg, f(1-g), (1-f)g, (1-f)(1-g)$}, and \keepmath{$|A|,|B|\le1$}.} Moreover, the derivative of $A$ with respect to $f$, denoted by $A_f$, is a convex combination, with weights $g$ and $1-g$, of two differences of coefficients in $[-1,1]$. This gives
\begin{equation} \label{eq:subset-bilinear-bounds}
|A_f|,|A_g|,|B_f|,|B_g|\le2.
\end{equation}
For brevity, write $q:=\{f,g\}_\alpha$, $r:=R_\alpha f$, $s:=R_\alpha g$, and $J:=A_fB_g-A_gB_f$. The formula \eqref{eq:contact-chain-rule} specializes to
\begin{equation} \label{eq:subset-chain-rule}
\begin{split}
\{A(f,g),B(f,g)\}_\alpha&=Jq+(AB_f-BA_f+gJ)r +(AB_g-BA_g-fJ)s.
\end{split}
\end{equation}
 By \eqref{eq:subset-bilinear-bounds}, $|J|\le8$ and each coefficient of $r$ and $s$ in \eqref{eq:subset-chain-rule} has absolute value at most $12$. \rev{The bracket is therefore controlled by the original bracket and the two Reeb derivatives. Applying Lemma~\ref{lem:subset-reeb-control} bounds its norm by}
\[ 8\norm q_\infty+12\norm r_\infty+12\norm s_\infty\le56\max_{\substack{u\in\{f,1-f\}\\v\in\{g,1-g\}}}\norm{\{u,v\}_\alpha}_\infty. \]
The bound is independent of $a,b$, so the same right-hand side bounds $\Palpha((h_{ij}))$, \rev{and hence also bounds \keepmath{$pb_\alpha(\{U_{ij}\})$}. Taking the infimum over all admissible pairs proves} \eqref{eq:subset-square-cover}. \end{proof}

\begin{lemma}\label{prop:subset-four-to-three} Every admissible quadruple $X_0, X_1, Y_0, Y_1$ satisfies 
\[ pb_{3,\alpha}(X_0\cup Y_0,X_1\cup Y_0,Y_1)\le3pb_{4,\alpha}(X_0,X_1,Y_0,Y_1). \]
\end{lemma}

\begin{proof} The three sets $X_0\cup Y_0$, $X_1\cup Y_0$, and $Y_1$ have empty common intersection, since $(X_0\cup Y_0)\cap(X_1\cup Y_0)=Y_0$. For an admissible pair $(f,g)$, set
\[ u=fg,\qquad v=(1-f)g,\qquad w=1-g. \]
These functions form an admissible triple for the three sets above. For brevity, write $q:=\{f,g\}_\alpha$, $r:=R_\alpha f$, \rev{and \keepmath{$s:=R_\alpha g$}. Formula~\eqref{eq:subset-chain-rule} gives}
\begin{align*}
\{u,v\}_\alpha&=gq-fgs,\\
\{u,w\}_\alpha&=-g(q+r)-f(1-g)s,\\
\{v,w\}_\alpha&=g(q+r-s)-(1-f)(1-g)s.
\end{align*}
By \eqref{eq:subset-four-brackets}, each of $q$, $q+r$, and $q+r-s$ is, up to sign, one of the four complementary brackets. \rev{To control the remaining terms, we use \eqref{eq:subset-square-reeb} to bound} $\norm{s}_\infty$ by twice the maximum of their norms. Therefore, we have 
\[ \max_{a,b\in\{u,v,w\}}\norm{\{a,b\}_\alpha}_\infty\le3\max_{\substack{a\in\{f,1-f\}\\b\in\{g,1-g\}}}\norm{\{a,b\}_\alpha}_\infty. \]
Hence, taking infima proves the conclusion. \end{proof}

\begin{lemma}\label{lem:displaceable-complement-subheavy} Suppose that $\mu$ is a nonzero monotone homogeneous quasi-morphism with the vanishing property. If $S\subset Y$ is compact and $Y\setminus S$ is contact displaceable, then $S$ is $\mu$-subheavy.
\end{lemma}

\begin{proof} Without loss of generality, we may normalize $\mu$ by \eqref{eq:mu-normalization}, since multiplying it by a positive constant does not change subheaviness. Let $h\in C^\infty(Y)$ satisfy $h|_S=0$. For each $\epsilon>0$, choose a smooth cutoff $\chi_\epsilon:Y\to[0,1]$ which vanishes near $S$ and equals $1$ outside $\{|h|<\epsilon\}$. Then $h_\epsilon=\chi_\epsilon h$ has support compactly contained in $Y\setminus S$, and $\norm{h-h_\epsilon}_\infty\le\epsilon$. \rev{The cutoff makes the vanishing property applicable, giving} $\mu(\widetilde\varphi_{h_\epsilon})=0$.

\rev{We next compare the cutoff Hamiltonian with the original one. Fix} $t>0$. \rev{By \eqref{eq:relative-flow-hamiltonian}, the path} $\varphi_h^{-s}\varphi_{h_\epsilon}^{s}$, $0\le s\le t$, is generated by $e^{-\kappa(\varphi_h^s)}(h_\epsilon-h)\circ\varphi_h^s$, \rev{and has \keepmath{$\alpha$}-length at most} $\epsilon\int_0^t e^{s\norm{R_\alpha h}_\infty}\,ds$. \rev{Combining this length estimate with \eqref{stability} and \eqref{eq:real-homogeneity}, we obtain}
\[ t\,|\mu(\widetilde\varphi_h)|\le D(\mu)+\epsilon\int_0^t e^{s\norm{R_\alpha h}_\infty}\,ds. \]
\rev{Letting \keepmath{$\epsilon\to0$} with \keepmath{$t$} fixed, and then letting \keepmath{$t\to\infty$}, we obtain} $\mu(\widetilde\varphi_h)=0$, which proves subheaviness. \end{proof}

\rev{We now prove Theorem~\ref{thm:subset-bounds}.} 

\begin{proof}[Proof of Theorem \ref{thm:subset-bounds}] Let $\vec h=(h_0,h_1,h_2)$ be admissible for the triple $(X_0,X_1,X_2)$. Since $h_i$ vanishes on $X_i$ and $X_i$ is $\mu$-subheavy, we have $\mu(\widetilde\varphi_{h_i})=0$ for every $i$. Lemma~\ref{lem:vanishing-partition-bound}, applied with $N=3$, gives
\[ \Palpha(\vec h)\ge\frac{1}{16(3-1)^2D(\mu)}=\frac{1}{64D(\mu)}. \]
\rev{To turn this bound for the bracket functional into a bound for the triple invariant, use the pointwise estimate from the proof of Proposition~\ref{prop:subset-cover-comparison}:}  $\Palpha(\vec h)\le4\max_{0\le i<j\le2}\norm{\{h_i,h_j\}_\alpha}_\infty.$ Therefore, every admissible triple has bracket maximum at least $1/(256D(\mu))$. Taking the infimum proves (i).

For (ii), let $(f,g)$ be admissible for the quadruple and put
\[ f_0=f,\qquad f_1=1-f,\qquad g_0=g,\qquad g_1=1-g,\qquad h_{ij}=f_i g_j. \]
These four nonnegative functions sum to $1$. For each $i,j$, the product $h_{ij}$ vanishes on $X_i\cup Y_j$, since $f_i$ vanishes on $X_i$ and $g_j$ vanishes on $Y_j$. \rev{\rrev{Therefore}, subheaviness of the union gives} $\mu(\widetilde\varphi_{h_{ij}})=0$. \rev{We can therefore apply Lemma~\ref{lem:vanishing-partition-bound}, now with} $N=4$. \rev{Combining its lower bound with the product-partition estimate in Proposition~\ref{prop:subset-cover-comparison} yields}
\[ \frac{1}{16(4-1)^2D(\mu)}\le\Palpha((h_{ij}))\le56\max_{\substack{u\in\{f,1-f\}\\v\in\{g,1-g\}}}\norm{\{u,v\}_\alpha}_\infty. \]
Since $56\cdot16\cdot9=8064$, taking the infimum over $(f,g)$ proves the first bound in (ii).

\rev{The additional assumption gives a sharper bound through the triple invariant. Namely, suppose that} $Y_1$ \rev{is} $\mu$-subheavy. The three sets $X_0\cup Y_0$, $X_1\cup Y_0$, and $Y_1$ are then $\mu$-subheavy, and their common intersection is empty because
\[ (X_0\cup Y_0)\cap(X_1\cup Y_0)\cap Y_1=Y_0\cap Y_1=\varnothing. \]
Therefore, part (i) and \rev{Lemma}~\ref{prop:subset-four-to-three} give
\postdisplaypenalty=10000
\[ \frac{1}{256D(\mu)}\le pb_{3,\alpha}(X_0\cup Y_0,X_1\cup Y_0,Y_1)\le3pb_{4,\alpha}(X_0,X_1,Y_0,Y_1). \]
This proves the last conclusion. \end{proof}

\begin{proof}[Proof of Corollary \ref{cor:subset-energy-bounds}] For triples, the three complements form a displaceable open cover. Apply \eqref{eq:subset-triple-cover} and Theorem~\ref{thm:spectral-cover} to obtain
\[ pb_{3,\alpha}(X_0,X_1,X_2)\ge\frac14\frac{1}{16(3-1)\sum_{i=0}^2e_\alpha(Y\setminus X_i)}=\frac{1}{128\sum_{i=0}^2e_\alpha(Y\setminus X_i)}. \]
\rev{The same comparison applies to quadruples: the four sets} $U_{ij}=Y\setminus(X_i\cup Y_j)$ \rev{form a displaceable open cover.} By \eqref{eq:subset-square-cover} and Theorem~\ref{thm:spectral-cover},
\[ pb_{4,\alpha}(X_0,X_1,Y_0,Y_1)\ge\frac1{56}\frac{1}{16(4-1)\sum_{i,j=0}^1e_\alpha(U_{ij})}=\frac{1}{2688\sum_{i,j=0}^1e_\alpha(U_{ij})}. \]
This completes the proof. \end{proof}

\section{Proof of Theorem \ref{thm:subset-chords}}\label{sec:proof-I}

\begin{lemma}\label{lem:averaging}
Assume $R_\alpha h=0$. Suppose that $T,L>0$ and
\[
\varphi_1^s\varphi_h^t(X_0)\cap X_1=\varnothing
\qquad\text{for }|t|\le T\text{ and }|s|\le L.
\]
Then there is a smooth function $f:Y\to[0,1]$, equal to $i$ near $X_i$ for $i=0,1$, such that
\begin{equation}\label{eq:averaging-bounds}
\norm{X_hf}_\infty\le\frac1T,
\qquad \norm{R_\alpha f}_\infty\le\frac1L.
\end{equation}
\end{lemma}

\begin{proof}
Consider the two compact sets
\[
K_i:=\{\varphi_1^s\varphi_h^t(x):x\in X_i,\ 0\le t\le T,\ 0\le s\le L\},
\qquad i=0,1.
\]
They are disjoint. Indeed, if $\varphi_1^{s_0}\varphi_h^{t_0}(x_0)=\varphi_1^{s_1}\varphi_h^{t_1}(x_1)$ with $x_i\in X_i$ and the indicated parameters, then commutation gives
\[
x_1=\varphi_1^{s_0-s_1}\varphi_h^{t_0-t_1}(x_0).
\]
Since $|t_0-t_1|\le T$ and $|s_0-s_1|\le L$, this contradicts the hypothesis. Choose a smooth function $f_0:Y\to[0,1]$ equal to $i$ on a neighborhood of $K_i$, and define
\begin{equation}\label{eq:double-average}
f(x):=\frac1{TL}\int_0^T\int_0^L
f_0(\varphi_1^s\varphi_h^t(x))\,ds\,dt.
\end{equation}
Then $f$ is smooth and $0\le f\le1$. Moreover, compactness of $[0,T]\times[0,L]\times X_i$ gives a neighborhood of $X_i$ whose image under every $\varphi_1^s\varphi_h^t$ in this rectangle lies where $f_0=i$. Therefore, $f=i$ on that neighborhood.

To estimate the derivative along the Hamiltonian flow, substitute $\varphi_h^\tau(x)$ for $x$ in \eqref{eq:double-average}. By the flow property and commutation,
\[
f(\varphi_h^\tau(x))
=\frac1{TL}\int_0^L\int_\tau^{T+\tau}
f_0(\varphi_1^s\varphi_h^t(x))\,dt\,ds.
\]
Differentiating at $\tau=0$ gives
\[
X_hf(x)=\frac1{TL}\int_0^L
\bigl[f_0(\varphi_1^s\varphi_h^T(x))-f_0(\varphi_1^s(x))\bigr]\,ds.
\]
Since each difference has absolute value at most one, $\norm{X_hf}_\infty\le1/T$. Similarly, substituting $\varphi_1^\sigma(x)$ in \eqref{eq:double-average} shifts the inner integration interval to $[\sigma,L+\sigma]$. Hence
\[
R_\alpha f(x)=\left.\frac{d}{d\sigma}\right|_{\sigma=0}f(\varphi_1^\sigma(x))
=\frac1{TL}\int_0^T
\bigl[f_0(\varphi_1^L\varphi_h^t(x))-f_0(\varphi_h^t(x))\bigr]\,dt.
\]
The same estimate proves $\norm{R_\alpha f}_\infty\le1/L$.
\end{proof}

We next convert these derivative bounds into a bound for $pb_{4, \alpha}$. In the setting of Theorem~\ref{thm:subset-chords}, write 
\[
a:=\max_{Y_0}h,\qquad b:=\min_{Y_1}h,\qquad \Delta:=b-a>0.
\]

\begin{lemma}\label{lem:cutoff}
In the setting of Theorem~\ref{thm:subset-chords}, let $f:Y\to[0,1]$ be smooth with $f=i$ near $X_i$. Then
\begin{equation}\label{eq:cutoff-bound}
p\le\max\left\{\norm{R_\alpha f}_\infty,
\frac{\norm{X_hf}_\infty+\max\{|a|,|b|\}\norm{R_\alpha f}_\infty}{\Delta}\right\}.
\end{equation}
\end{lemma}

\begin{proof}
In general, for a smooth scalar function $\chi: \R \to \R$, the defining equations of a contact Hamiltonian vector field give
\begin{equation}\label{eq:cutoff-field}
X_{\chi(h)}=\chi'(h)X_h+\bigl(\chi(h)-h\chi'(h)\bigr)R_\alpha.
\end{equation}
Indeed, the right-hand side has $\alpha$-value $\chi(h)$, and its contraction with $d\alpha$ is $-\chi'(h)dh=-d(\chi(h))$. Also, $R_\alpha(\chi(h))=0$, so these are precisely the defining equations for $X_{\chi(h)}$. Now, put $g:=\chi(h)$. Since $R_\alpha g=0$ and $\{1,f\}_\alpha=R_\alpha f$, the four brackets are
\[
\begin{aligned}
\{f,g\}_\alpha&=-X_gf,&
\{1-f,g\}_\alpha&=X_gf,\\
\{f,1-g\}_\alpha&=X_gf-R_\alpha f,&
\{1-f,1-g\}_\alpha&=R_\alpha f-X_gf.
\end{aligned}
\]
Therefore, their maximum norm is
\begin{equation}\label{eq:two-derivatives}
\max\{\norm{X_gf}_\infty,\ \norm{(R_\alpha-X_g)f}_\infty\}.
\end{equation}
In particular, if $\chi(r)=(r-a)/\Delta$ on $(a,b)$, then at any $x$ with $a<h(x)<b$ the two vector fields are $(X_h-aR_\alpha)/\Delta$ and $(bR_\alpha-X_h)/\Delta$. 

For the smoothing step, fix $0<\ep<\Delta/4$ and define
\[
q_\ep(r):=
\begin{cases}
0,&r\le a+2\ep,\\[2mm]
\dfrac{r-a-2\ep}{\Delta-4\ep},&a+2\ep<r<b-2\ep,\\[3mm]
1,&r\ge b-2\ep.
\end{cases}
\]
Set $g_\ep(x):=\int_{\mathbb R}q_\ep(h(x)-v)\rho_\ep(v)\,dv$, where
\[
\rho_\ep\in C_c^\infty(\mathbb R),\quad\rho_\ep\ge0,\qquad
\int_{\mathbb R}\rho_\ep(v)\,dv=1,\qquad
\operatorname{supp}\rho_\ep\subset[-\ep,\ep].
\]
Here $v$ is a scalar shift of the cutoff, not a time parameter for a flow. Fix a point $x$. We first examine the shifted function $q_\ep(h-v)$, before taking its average. The useful identity is
\[
X_{Ah+B}=A X_h+B R_\alpha
\]
for constants $A,B$. Then there are three possibilities.
\begin{itemize}
\item[(i)] If $h(x)-v<a+2\ep$, then $q_\ep(h-v)$ is locally the constant $0$. Its vector field at $x$ is $0$.
\item[(ii)] If $a+2\ep<h(x)-v<b-2\ep$, then locally $q_\ep(h-v)=\frac{h-(a+2\ep+v)}{\Delta-4\ep}$. Hence its vector field at $x$ is
\[
\frac{X_h-(a+2\ep+v)R_\alpha}{\Delta-4\ep}.
\]
\item[(iii)] If $h(x)-v>b-2\ep$, then $q_\ep(h-v)$ is locally the constant $1$. Its vector field at $x$ is $R_\alpha$.
\end{itemize}

Write $V_{\ep,v}$ for the vector field given by these three cases. Differentiating the defining integral for $g_{\ep}$ gives
\[
X_{g_\ep}(x)=\int_{\mathbb R}V_{\ep,v}(x)\rho_\ep(v)\,dv.
\]
All vectors in this integral belong to the same tangent space $T_xY$. The three formulas describe the integrand; the smoothed vector field is their average. Also, we have
\[
R_\alpha(x)-X_{g_\ep}(x)
=
\int_{\mathbb R}\bigl(R_\alpha(x)-V_{\ep,v}(x)\bigr)\rho_\ep(v)\,dv.
\]
In particular, in the transition interval, the vector inside this second integral is
\[
\begin{aligned}
R_\alpha-\frac{X_h-(a+2\ep+v)R_\alpha}{\Delta-4\ep}
&=\frac{(\Delta-4\ep+a+2\ep+v)R_\alpha-X_h}{\Delta-4\ep} =\frac{(b-2\ep+v)R_\alpha-X_h}{\Delta-4\ep},
\end{aligned}
\]
because $\Delta+a=b$. On the two constant regions, it is $R_\alpha$ and $0$, respectively.

Now, we apply these vectors to $f$. On the transition interval,
\[
\left|
\frac{X_hf-(a+2\ep+v)R_\alpha f}{\Delta-4\ep}
\right|
\le
\frac{|X_hf|+|a+2\ep+v|\,|R_\alpha f|}{\Delta-4\ep}.
\]
Since $|v|\le\ep$,
\[
|a+2\ep+v|\le |a|+3\ep.
\]
Similarly,
\[
\left|
\frac{(b-2\ep+v)R_\alpha f-X_hf}{\Delta-4\ep}
\right|
\le
\frac{|X_hf|+(|b|+3\ep)|R_\alpha f|}{\Delta-4\ep}.
\]
Thus the term $3\ep$ comes from two contributions: $2\ep$ from moving the cutoff endpoints, and another $\ep$ from the smoothing shift $v$. Meanwhile, on the constant regions, the corresponding derivatives are either $0$ or $R_\alpha f$. Therefore, every scalar term in either average has absolute value at most
\[
C_\ep:=
\max\left\{
\|R_\alpha f\|_\infty,\,
\frac{\|X_hf\|_\infty+(\max\{|a|,|b|\}+3\ep)\|R_\alpha f\|_\infty}{\Delta-4\ep}
\right\}.
\]
This implies that 
\[
\max\left\{
\|X_{g_\ep}f\|_\infty,\,
\|(R_\alpha-X_{g_\ep})f\|_\infty
\right\}
\le C_\ep.
\]

Finally, since $R_\alpha h=0$, we have $R_\alpha g_\ep=0$. Therefore, the four brackets in our definition of $pb_{4,\alpha}$ are $\pm X_{g_\ep}f$ and $\pm(R_\alpha-X_{g_\ep})f$. Also, $g_\ep=0$ on the neighborhood $\{h<a+\ep\}$ of $Y_0$, and $g_\ep=1$ on the neighborhood $\{h>b-\ep\}$ of $Y_1$. Thus $(f,g_\ep)$ is a smooth admissible pair, and the definition as an infimum gives
\[
p\le C_\ep.
\]
This inequality holds for every sufficiently small $\ep>0$. Letting $\ep$ go to $0$ gives
\[
p\le
\max\left\{
\|R_\alpha f\|_\infty,\,
\frac{\|X_hf\|_\infty+\max\{|a|,|b|\}\|R_\alpha f\|_\infty}{\Delta}
\right\}.
\]
This completes the proof. \end{proof}

Now, we are ready to give the proof of Theorem~\ref{thm:subset-chords}.

\begin{proof}[Proof of Theorem~\ref{thm:subset-chords}]
Assume first that \eqref{eq:main-condition} is strict, and suppose that no connection exists with $|t|\le T$ and $|s|\le L$. Lemma~\ref{lem:averaging} supplies an admissible separating function $f$ with the derivative bounds \eqref{eq:averaging-bounds}. Hence Lemma~\ref{lem:cutoff} gives
\[
p\le\max\left\{\frac1L,
\frac1\Delta\left(\frac1T+\frac{\max\{|a|,|b|\}}L\right)\right\}<p,
\]
a contradiction. Therefore, a connection exists whenever \eqref{eq:main-condition} is strict.

For the non-strict inequality, replace $T,L$ by $(1+1/j)T,(1+1/j)L$. The left-hand side of \eqref{eq:main-condition} is then divided by $1+1/j$, so the preceding argument gives $x_j\in X_0$ and times $t_j,s_j$ satisfying
\[
\varphi_1^{s_j}\varphi_h^{t_j}(x_j)\in X_1,
\qquad |t_j|\le(1+1/j)T,\qquad |s_j|\le(1+1/j)L.
\]
These triples lie in the compact set $X_0\times[-2T,2T]\times[-2L,2L]$. Passing to a subsequence gives a limit $(x,t,s)$ with $|t|\le T$ and $|s|\le L$. Since the flows are continuous and $X_1$ is closed, $\varphi_1^s\varphi_h^t(x)\in X_1$. Also, $t=s=0$ would imply $x\in X_0\cap X_1$, which is impossible. Finally, if the disjointness condition $\left(\bigcup_{|s|\le L}\varphi_1^s(X_0)\right)\cap X_1=\varnothing$ holds, then $t=0$ would imply $\varphi_1^s(x)\in X_1$ for some $|s|\le L$, again a contradiction. Hence we conclude $t\ne0$, as desired.
\end{proof}

\section{Proof of Theorem \ref{thm:contact-form-deformation}}\label{sec:proof-J}

\begin{proof}[Proof of Theorem \ref{thm:contact-form-deformation}] \rev{We begin with \keepmath{$d\alpha_s=d\alpha-s\,df\wedge dg$}.} For $n\ge1$, the identity $(df\wedge dg)\wedge(df\wedge dg)=0$ makes every term of order at least two in the binomial expansion of $(d\alpha_s)^n$ vanish. Therefore,
\[ (d\alpha_s)^n=(d\alpha)^n-ns\,df\wedge dg\wedge(d\alpha)^{n-1}. \]
Wedging with $\alpha_s$ gives
\[
\begin{aligned}
\alpha_s\wedge(d\alpha_s)^n =\alpha\wedge(d\alpha)^n -s\bigl(f\,dg\wedge(d\alpha)^n+n\,\alpha\wedge df\wedge dg\wedge(d\alpha)^{n-1}\bigr),
\end{aligned}
\]
where the possible quadratic term vanishes because it contains $dg$ twice. To identify the coefficient of $s$, recall that $\alpha(X_f)=f$ and $\iota_{X_f}d\alpha=(R_\alpha f)\alpha-df$. The form $\alpha\wedge dg\wedge(d\alpha)^n$ has degree $2n+2$, so it vanishes on $Y$. Contracting this zero form with $X_f$ gives
\[
\begin{aligned}
0&=\iota_{X_f}\bigl(\alpha\wedge dg\wedge(d\alpha)^n\bigr)\\
&=f\,dg\wedge(d\alpha)^n-(X_fg)\,\alpha\wedge(d\alpha)^n +n\,\alpha\wedge dg\wedge\bigl((R_\alpha f)\alpha-df\bigr)\wedge(d\alpha)^{n-1}.
\end{aligned}
\]
\rev{Here the term containing} $\alpha$ \rev{twice vanishes, while} $dg\wedge df=-df\wedge dg$\rev{. Hence the remaining terms yield}
\[ f\,dg\wedge(d\alpha)^n+n\,\alpha\wedge df\wedge dg\wedge(d\alpha)^{n-1}=(X_fg)\,\alpha\wedge(d\alpha)^n. \]
Substituting this identity into the expansion proves
\[ \alpha_s\wedge(d\alpha_s)^n=(1-sX_fg)\,\alpha\wedge(d\alpha)^n. \]
When $n=0$, we have $dg=(R_\alpha g)\alpha$ and $X_f=fR_\alpha$, so the same formula follows directly from $\alpha_s=\alpha-sf\,dg$. \rev{The contact condition now follows from this identity: since} $\alpha\wedge(d\alpha)^n$ is nowhere zero, $\alpha_s$ is contact whenever $1-sX_fg$ is nowhere zero. In particular, $|s|\norm{X_fg}_\infty<1$ implies $1-sX_fg>0$ everywhere (hence the corresponding deformation $\alpha_s$ is contact). If $X_fg=0$, all real parameters are allowed.

\rev{Next, we construct the desired contactomorphism using a variable-time Hamiltonian flow. For this purpose, recall that the contact Hamiltonian equations imply} $\mathcal L_{X_f}\alpha=(R_\alpha f)\alpha$ and $X_ff=fR_\alpha f$. Thus, for each fixed time $t$,
\[ (\varphi_f^t)^*\alpha=e^{\kappa(\varphi_f^t)}\alpha,\qquad f(\varphi_f^t(x))=e^{\kappa(\varphi_f^t)(x)}f(x), \]
where $\kappa(\varphi_f^t)(x)=\int_0^t R_\alpha f(\varphi_f^r(x))\,dr$. \rev{These identities also hold for negative} $t$\rev{. We may therefore define}
\[ \Psi_s(x):=\varphi_f^{-sg(x)}(x),\qquad \kappa_s(x):=\int_0^{-sg(x)}R_\alpha f(\varphi_f^r(x))\,dr. \]
Since the flow time depends on $x$, differentiating in a direction $v\in T_xY$ gives two terms:
\[ d\Psi_s(x)v=d\varphi_f^{-sg(x)}(x)v-s\,dg_x(v)X_f(\Psi_s(x)). \]
Evaluating $\alpha$ on this vector and using $f(\Psi_s(x))=e^{\kappa_s(x)}f(x)$, we obtain
\[
\begin{aligned}
(\Psi_s^*\alpha)_x(v)
&=e^{\kappa_s(x)}\alpha_x(v)-s\,dg_x(v)f(\Psi_s(x)) =e^{\kappa_s(x)}\bigl(\alpha_x(v)-sf(x)\,dg_x(v)\bigr).
\end{aligned}
\]
Therefore, $\Psi_s^*\alpha=e^{\kappa_s}\alpha_s$. Moreover, the integral defining $\kappa_s$ gives $|\kappa_s(x)|\le |s|\,|g(x)|\norm{R_\alpha f}_\infty,$ which proves the estimate in \eqref{eq:contact-deformation-pullback}.

\rev{The pullback formula and its estimate are now established. It remains to verify that} $\Psi_s$ \rev{is a diffeomorphism on the stated interval. We first prove local invertibility: the pullback identity implies}
\[ \Psi_s^*\bigl(\alpha\wedge(d\alpha)^n\bigr)=e^{(n+1)\kappa_s}\alpha_s\wedge(d\alpha_s)^n, \]
since every term containing $d\kappa_s\wedge\alpha_s$ vanishes after wedging with $\alpha_s$. For $|s|\norm{X_fg}_\infty<1$, the right-hand side is a positive volume form by the computation above. Hence $d\Psi_s$ is invertible at every point, and $\Psi_s$ is an orientation-preserving local diffeomorphism. \rev{To pass from local to global invertibility, note that} $Y$ \rev{is closed and each flow trajectory stays in its connected component. Hence} the restriction of $\Psi_s$ to each component is a covering map onto that component. \rev{Moreover, the homotopy} $r\mapsto\Psi_{rs}$, $0\le r\le1$, \rev{connects this map to the identity, so its degree is one.} Therefore, the covering has one sheet and $\Psi_s$ is a diffeomorphism. Smooth dependence on $s$ gives the required isotopy.
\end{proof}

\begin{proof}[Proof of Corollary \ref{cor:cover-deformation}] Since $T \cdot pb_\alpha(\mathcal U)<1$, the definition of the infimum gives a smooth subordinate partition $\vec f$ such that $T \cdot \Palpha(\vec f)<1$. Fix $a,b\in[0,1]^N$ and put $f=\sum_i a_if_i$ and $g=\sum_i b_if_i$. Then $0\le g\le1$ and $1-g=\sum_i(1-b_i)f_i$. \rev{Thus each coefficient vector in the brackets} $\{f,g\}_\alpha$ \rev{and} $\{f,1-g\}_\alpha$ \rev{belongs to} $[-1,1]^N$\rev{. The partition functional therefore bounds both brackets, giving}
\[ \norm{X_fg}_\infty\le\max\bigl\{\norm{\{f,g\}_\alpha}_\infty,\norm{\{f,1-g\}_\alpha}_\infty\bigr\}\le\Palpha(\vec f)<\frac1T. \]
Since the same bound holds for every choice of $a,b$, Theorem~\ref{thm:contact-form-deformation} applies for every $|s|\le T$, which completes the proof. \end{proof}

\bibliographystyle{amsplain} 
\bibliography{C_0_CPB_v44,cpb45-supplement}

\

\end{document}